\documentclass[12pt,a4paper]{article}
\usepackage{amsmath,amssymb,amsthm,amsfonts,amscd,euscript,verbatim}
\usepackage{t1enc}
\usepackage{hyperref}
\usepackage{newlfont}
 \theoremstyle{plain}
\newtheorem{theo}{Theorem}[subsection]

\newtheorem{pr}[theo]{Proposition}

 \newtheorem{lem}[theo]{Lemma}
 \newtheorem{coro}[theo]{Corollary}
  \newtheorem{conj}[theo]{Conjecture}
\theoremstyle{remark}
\newtheorem{rema}[theo]{Remark}

\theoremstyle{definition}
\newtheorem{defi}[theo]{Definition}
\newtheorem*{notat}{Notation}

\newcommand\ob{^{-1}}
\newcommand\smc{{SmCor}}
\newcommand\ssc{{Shv(SmCor)}}
\newcommand\ssca{{Shv(SmCor)_A}}
\newcommand\bmk{{B^-(Shv(SmCor))}}
\newcommand\dmk{{D^-(Shv(SmCor))}}
\newcommand\dmka{{D(Shv(SmCor)_A)}}
\newcommand\cmk{{C^-(Shv(SmCor))}}
\newcommand\kmk{{K^-(Shv(SmCor))}}
\newcommand\kmka{{K(Shv(SmCor)_A)}}
\newcommand\psc{{PreShv(SmCor)}}
\newcommand\dpk{{D^-(PreShv(SmCor))}}
\newcommand\dpe{DPM^{eff}}
\newcommand\hk{\mathfrak{H}}
\newcommand\hkgm{\mathfrak{H}_{gm}}

 \newcommand\zc{\mathcal{Z}{}}

\newcommand\dmge{DM^{eff}_{gm}}
\newcommand\dmgm{DM_{gm}}
\newcommand\dme{DM_-^{eff}}
\newcommand\dmea{DM_A^{eff}}
\newcommand\dms{DM^{s}{}}

\newcommand\mg{M_{gm}}
\newcommand\obj{Obj}

\newcommand\id{id}
\newcommand\cu{\underline{C}}

\newcommand\z{{\mathbb{Z}}}
\newcommand\q{{\mathbb{Q}}}
\newcommand\af{\mathbb{A}}
\newcommand\p{\mathbb{P}}
\newcommand\db{\mathcal{D}}
\newcommand\dbk{\mathcal{D}^{op}}
\newcommand\dbf{\mathcal{D}_{fin}^{op}{}}
\newcommand\pt{pt}
\newcommand\eps{\varepsilon}
\newcommand\si{\sigma}
\newcommand\ns{\{0\}}
\DeclareMathOperator\homm{\operatorname{Hom}}

\newcommand\ihom{{\underline{Hom}}}
\newcommand\chow{Chow^{eff}}
\newcommand\chown{Chow}
\newcommand\cho{Corr_{rat}}
\newcommand\ab{Ab}
\newcommand\var{Var}
\newcommand\sv{SmVar}
\newcommand\spv{SmPrVar}
\newcommand\sprv{Sch^{prop}}
\newcommand\spe{Spec}
 \DeclareMathOperator\ke{\operatorname{Ker}}
 \DeclareMathOperator\cok{\operatorname{Coker}}
\DeclareMathOperator\imm{\operatorname{Im}}
\DeclareMathOperator\co{\operatorname{Cone}}
\DeclareMathOperator\prt{Pre-Tr}
\DeclareMathOperator\prtp{Pre-Tr^+}
\newcommand\trp{Tr^+}
\newcommand\tr{Tr}

\DeclareMathOperator\inli{\varinjlim}
\DeclareMathOperator\en{\operatorname{End}}

\begin{document}

 \title{Differential graded motives:
 weight complex, weight filtrations and spectral sequences for
 realizations; Voevodsky vs. Hanamura}

 \author{M.V. Bondarko
  \thanks{ The author gratefully acknowledges the support from Deligne 2004
Balzan prize in mathematics. The work is also supported by RFBR
(grant no.
 08-01-00777a) and INTAS (grant no.
 05-1000008-8118).}}

\maketitle

\begin{abstract}

We describe  explicitly  the Voeovodsky's triangulated category of
motives $\dmge$ (and describe a 'differential graded enhancement' for
it). This enables us to able to verify that $\dmgm \q$ is
(anti)isomorphic to  Hanamura's $\db(k)$.

We obtain a  description of all subcategories (including those of
Tate motives) and of all localizations of $\dmge$.  We construct a
conservative {\it weight complex} functor $t:\dmge\to K^b(\chow)$;
$t$ induces an isomorphism $K_0(\dmge)\to K_0(\chow)$. A motif is
mixed Tate if and only if its weight complex is. Over finite fields the
Beilinson-Parshin conjecture holds iff $t \q$ is an equivalence.

  For a
realization $D$ of $\dmge$ we construct a spectral sequence $S$ (the
{\it spectral sequence of motivic descent}) converging to the
cohomology of  an arbitrary motif $X$.  $S$ is 'motivically
functorial'; it gives a canonical functorial weight filtration on
the cohomology of $D(X)$. 
For the 'standard'
realizations this filtration coincides with the usual one (up to a
shift of indices). For the motivic cohomology this weight filtration
is non-trivial and appears to be quite new.

 We define the (rational) {\it length} of a motif $M$; 
 modulo certain 'standard' conjectures this length
 coincides with the maximal length of the weight filtration of the
 singular cohomology of $M$.



 \end{abstract}

\tableofcontents

 \section*{Introduction}

We give  an explicit description of the category of effective
geometric motives of Voevodsky (see \cite{1}).  In what follows
$\dms$ is the full triangulated subcategory of the category $\dmge$
(defined in \cite{1})
 generated by motives of  smooth  varieties (we do not add the kernels of
projectors). It is proved that for any motivic complex $M$ (i.e.
an object of Voevodsky's $\dme$ that comes from $\dms$; in
particular, the Suslin complex of an arbitrary variety) there
exists a quasi-isomorphic complex $M'$ 'constructed from' the
Suslin complexes of smooth projective varieties; $M'$ is unique up
to a homotopy.  We prove that   Hanamura's  $\db(k)$ is
(anti)equivalent to Voeovodsky's $\dmgm\q$.

Our main category $\hk$ is defined as the category of {\it twisted
complexes} over a certain differential graded category whose
objects are cubical Suslin complexes; we construct an equivalence
$m:\hk \to \dms$.
 In terms of \cite{bk} our description of $\dms$ gives an
{\it enhancement} of this category. One should think of twisted
complexes as of the results of repetitive computation of cones of
morphisms in an 'enhanced' triangulated category. One can describe
any subcategory of $\hk$ that is generated by a fixed set of
objects; this method gives an integral description of the
triangulated category of Tate motives  similar to the rational
description of \cite{km}.
 Besides, any localization of $\hk$ can be described explicitly
using the construction of Drinfeld (see \cite{drinf}).

As an application we  consider the problem of constructing exact
functors from $\dms$ (i.e. realizations) in terms of cubical Suslin
motivic complexes.   The most simple and yet quite interesting of
functors constructed by our method are the {\it truncation functors}
$t_N$ that correspond to the canonical filtration of  Suslin
complexes.
 The target of $t_0$
is just the  category  $K^b(\cho)$ ($\cho$ is 'almost' the category
of effective Chow motives, see \ref{dvoev}). $t_0$ extends to
conservative weight complex functors $t:\dmge\to K^b(\chow)$ and
$t_{gm}:\dmgm\to K^b(\chown)$. We prove that $t$ induces an
isomorphism $K_0(\dmge)\to K_0(\chow)$ thus answering a question
from \cite{gs}.

We show that   if $W$ denotes the weight filtration on $H^i(X)$ for
a 'standard' realization $H$ then $W_{l+N} H^i(X)/W_{l-1} H^i(X)$
factors through $t_N$.

We prove that a motif $X$ belongs to a triangulated category
$M\subset\hk$ generated by motives of a given set of smooth
projective varieties $P_i$, if and only if the same is true for $t(X)$
(as a complex of Chow motives). In particular, the motif of a
smooth variety is a mixed Tate one iff its weight complex (as
defined by Gillet and Soul\'e) is.

For any realization $D$ of motives that belongs to a wide class of
'enhanceable' realizations (i.e. of realizations that admit a
differential graded 'enhancement') we construct a family of
'truncated realizations'.  In particular, this result can be applied to
'standard' realizations and motivic cohomology; an interesting new
family of realizations is obtained. This yields a canonical spectral
sequence $S$ converging to the cohomology of $D(X)$ of an arbitrary
motif $X$. $S$ could be called the {\it spectral sequence of motivic
descent}. The $E_1$-terms of $S$ are expressed in cohomology of smooth projective
varieties, the $E_{n}$-terms of $S$ have a nice description in terms
of $t_{2n-2}(X)$, $n\ge 1$. 
$S$ is canonical  and  'motivically functorial', it is also
functorial with respect to transformations of functors. $S$ gives a
canonical non-trivial weight filtration for 'differential graded'
realizations of motives; for the 'standard' realizations this
filtration coincides with the usual one (up to a shift of indices).

 The simplest case of $S$
 for motivic cohomology is the Bloch's long exact localization
sequence for higher Chow groups (see \cite{blo}).
 Our 'weight' filtration on motivic cohomology
is non-trivial and was not mentioned in the literature; it gives a
new filtration on the $K$-theory of a smooth variety $X$.

We also study motives  with compact support ($\mg^c$ in the notation
of Voevodsky). We give an explicit description of $\mg^c$ for  a
variety $Z$ and prove that the weight complex of Gillet and Soul\'e
can be described as $t_0(m\ob(\mg^c(Z)))$ (with reversed arrows);
the functor $h$ of Guillen and Navarro Aznar is (essentially)
$t_0(m\ob(\mg(Z)))$.

We define the 'length' of a motif ({\it stupid, fine or
rational}); this is a natural motivic analogue of the length of
weight filtration for a mixed Hodge structure. For a smooth
variety $X$ the length of $\mg(X)$ lies between the length of the
weight filtration of the singular cohomology of $X$ and the
dimension of $X$. If certain 'standard' conjectures are  valid
then the rational length of a motif coincides with the
(appropriately defined) length of the weight filtration of its
singular realization.

We note that in the current paper we apply several results of
\cite{1} that use resolution of singularities; yet applying de
Jong's alterations one can easily extend (most of) our results (at
least) to motives with rational coefficients over an arbitrary
perfect $k$. For a finite $k$ the Beilinson-Parshin conjecture (that
the only non-zero $H^i(X,\q(n))$ for smooth projective $X$ is
$H^{2n}$) holds iff $t \q:\dmge\q\to K^b(\chow \q)$ is an
equivalence (note that here $\dmge\q$ and $t \q$ denote the
appropriate idempotent completions). Much can also be proved with
integral coefficients.

The author would like to note that several interesting results of
this paper (in particular, the properties of $t_0$  and $t$) follow
just from the fact that $\hk=\tr(J)$ for $J$ being a negative
differential graded category (see subsection \ref{dhk}). Hence these
results are also valid for any other example of this situation.
Moreover, in  \cite{bws}  a set of axioms of so-called {\it weight
structures} for a triangulated category $C$ is introduced (see
Remark \ref{ass}). This (abstract) approach allows  extending (most
of) the results of the current paper to a wide class of triangulated
categories and realizations that do not necessarily have a
differential graded 'enhancement'; in particular, it can be applied
to the stable homotopy category (of spectra). The definition of the weight
filtration (and spectral sequences) for realizations is quite easy
in this context. This formalism is also closely (and nicely) related
to $t$-structures. Yet in this abstract setting it is difficult to
define truncation
functors (especially the 'higher' ones). 

Besides, recently Levine (in his very interesting paper \cite{lesm}) extended some of the results of the current paper to relative motives (i.e. motives over a regular base $S$, which is essentially of finite type over $k$). This paper also contains a nice exposition of the yoga of cubical objects
 (which is essential for the current paper) and of related tensor products.

Now we list the contents of the paper. More details can be found
in the beginnings of  sections.

 In the first section we recall some basic notation of \cite{1} (with minor modifications).
 Next we
describe cubical Suslin complexes and their properties. Most of
the proofs are postponed till \S\ref{prope} since they are not
important for the understanding of main results.

We start  \S2 by recalling the formalism of differential graded
categories (\ref{dg}) and twisted complexes (\ref{tw} and
\ref{prtp}).
 Next we use this formalism to construct our main objects of
study: a triangulated category $\hk$ and a functor $h$ from $\hk$
into the homotopy category of complexes of Nisnevich sheaves with
transfers (\ref{dhk}). 
 $\hk$ as a triangulated category is generated by motives
of smooth projective varieties. We also describe $\hk$ and $h$ more
explicitly (in \ref{expl}). We define 'stupid' filtration for
objects of $\hk$ that is similar to the 'stupid' filtration for
complexes over an additive category (\ref{sstufil}). Note that
(similarly to the case of complexes) the filtration for a motif also
depends on the choice of its 'lift' to a certain differential graded
category (as the stupid truncation object of $K(A)$ depends on its
lift to $C(A)$); yet this filtration is 'functorial enough' for our
purposes. In \ref{dop} we introduce more differential graded
definitions; they will be used in section \ref{new}.

In \S3 we prove (Theorem \ref{main}) that $h$ composed with the
natural functor from the homotopy category of sheaves with
transfers to the derived category  gives an equivalence  $m:\hk\to
\dms$. Lastly we give a certain differential graded description of
Voevodsky's $\dmgm$; to this end we define $\otimes \z(1)$ on
$\hk$ and formally invert this functor. 

In \S\ref{comphan} 
we  prove that Voevodsky's $\dmgm\q$ is (anti)isomorphic to
Hanamura's $\db(k)$.

In \S5 we verify the properties of cubical Suslin complexes. The
reader not interested in the proofs of  auxiliary results of \S1 may
skip this section.

 In \S6 using the canonical filtration of the (cubical) Suslin complex
  we define the 'truncation'
functors $t_N:\hk\to \hk_N$ (in \ref{tn}). $\hk_0$ (i.e. the target
of $t_0$) is just  $K^b(\cho)$. These functors are new though their
certain restrictions to varieties  were (essentially) considered in
\cite{gs} and \cite{gu} (and were shown to be quite important).
 We prove that  $t_0$ extends to $t:\dmge\to K^b(\chow)$ and  $t_{gm}:\dmgm\to K^b(\chown)$
 (in \ref{finlength}). All $t_N$ (see Theorem \ref{ttn}),  $t$ (Proposition
 \ref{funt}), and $t_{gm}$ (Remark \ref{nummot})
  are conservative. $t$ induces an isomorphism
$K_0(\dmge)\cong K_0(\chow)$ (Theorem \ref{kze}); $t_{gm}$ gives
$K_0(\dmgm)\cong K_0(\chown)$ . Certainly, this extends to an
isomorphism $K_0(\dmgm)\cong K_0(\chown)$ (Corollary \ref{newcor}).

We define the {\it length} of a motif (three types); the {\it stupid
length} (see \ref{motlen}) is not less than the {\it fine} one
(see \ref{finlength}), which is not less than the {\it rational
length}. We prove that motives of smooth varieties of dimension
$N$ have stupid length $\le N$ (parts 1 and 2 of Theorem
\ref{ttn}); besides $t_N(X)$ contains all information on motives
of stupid length $\le N$ (see part 3 of Theorem \ref{ttn} and
\ref{hodge}).

 At the end of the section we calculate $m\ob
(\mg^c(X))$ for a smooth $X$ explicitly (in \ref{expmgc}). Using
this result as well as $cdh$-descent we prove that the weight
complex of Gillet and Soul\'e for $X/k$ can be described as
$t_0(m\ob(\mg^c(X)))$ (with arrows reversed; see \ref{wecompl}).
Besides, $t_0(m\ob(\mg(X)))$ essentially coincides with the functor
$h$ described in Theorem 5.10 of \cite{gu}.

In \S \ref{new} we  study {\it realizations} of the category of
motives and their connections with (certain) weight filtrations.
 The differential graded categories formalism yields a general recipe of constructing
 realizations
(see \ref{enhreal}). It is quite easy to determine which of those
{\it enhanceable} realizations can be factored through $t_N$.

We verify that the \'etale and motivic cohomology are enhanceable
realization (see \S\ref{etdr}, \S\ref{sing}, and part 1 of Remark
\ref{rzn}). Very probably, this result could be extended to all
other 'standard' realizations. In \ref{hodge} for any enhanceable
realization $D$ we describe an interesting new family of
'truncated realizations'; they correspond to 'forgetting'
cohomology outside a given range of weights. Truncated
realizations give a filtration of the complex that computes the
given 'enhanced' realization of a motif $Y$. We obtain a spectral
sequence $S$ converging to $D^i(Y)$ (see (\ref{spectr})). Its
$E_{n}$-terms for $n\ge 2$ have a nice description in terms of
$t_{2n-2}(Y)$; in particular, $E_1$-terms are functorial in
$t_0(Y)$. $S$ is the {\it spectral sequence of motivic descent}.
$S$ gives a canonical integral weight filtration for 'enhanced'
realizations of motives; for the 'standard' realizations this
filtration coincides with the usual one (up to a shift of
indices). $S$ is 'motivically functorial', it is also functorial
with respect to ('enhanced') transformation of functors.

In \ref{clawe}  we prove that for a 'standard' $H$ the $N$-th
truncated realization computes $W_{l+N} H^l(X)/W_{l-1} H^l(X)$,
while $W$ equals the 'standard'  weight filtration.  A morphism
$f$ induces a zero morphism on cohomology if $t_0(f)$ is zero. We
also prove (modulo certain 'standard' conjectures) that the
rational length of a motif coincides with the 'range' of
difference of $l$ with the weights of $H^l$ for all $l$; see
Proposition \ref{procon}. We conclude the section with a discussion
of $qfh$-descent and motives of singular varieties (see
\ref{sing}).

In section \ref{last} we apply the general theory of \cite{bk} to
describe any subcategory of $\hk$ that is generated by a fixed set
of objects (see \ref{ssubcat}). In particular, this method can be
used to obtain an integral description of the triangulated category
of effective Tate motives (i.e. the full triangulated subcategory of
$\hk$ generated by $\z(n)$ for $n>0$).

In \ref{loc} we describe the construction of 'localization of
differential graded categories' (due to Drinfeld). This gives us a
description of localizations of $\hk$.  As an application, we
prove that the motif of a smooth $X/k$ is a mixed Tate one
if and only if the weight complex of $X$ (defined in \cite{gs}) is.

In \ref{motcharp} we verify that over an arbitrary perfect field
one can apply our theory (at least) with rational coefficients.
 Moreover, over
finite fields the Beilinson-Parshin conjecture holds iff $t
\q:\dmge\q\to K^b(\chow \q)$ is an equivalence. We also describe
an idea for constructing a certain 'infinite integral' weight
complex functor in finite characteristic.

In \ref{sindepl} we prove that traces of endomorphisms of cohomology
of motives induced by endomorphisms of motives do not depend on the
choice of a Weil cohomology theory. In particular, this statement can be
applied for the morphisms induced by 'open correspondences' (as
described in Definition 3.1 of \cite{bloesn}). We obtain a
generalization of Theorem 3.3. of loc. cit. to the case of varieties
which are not necessarily complements of  smooth projective
varieties by strict normal crossing divisors.

 In \ref{summand} we note  that one can modify the description
 of  $\dms$ so that $\z(n)$ will have stupid length $0$.  Lastly we
 describe certain functors $m_N:\hk_N\to \dme$ (see \ref{mn}). $t_N$
and $m_N$ could be related to the (yet conjectural) weight
filtration on $\dmge$.

The author is deeply grateful to prof. B. Kahn,  Dr. F. Deglise,
prof. A. Suslin, prof. M. Spiess, and to Dr. N. Dourov for  very
useful comments. The author is also deeply indebted to the
officers and guests of Max Planck Institut (Bonn) and of IHES
(Paris) for excellent working conditions during his stays in
these places.

\begin{notat}

In this paper all complexes will be cohomological i.e. the degree
of all differentials is $+1$.

We recall that for any triangulated $T$ there exists a unique
category $T'\supset T$ that is obtained from $T$ by 'adding the
kernels of all projectors'; $T'$ is called the idempotent
completion of $T$ (see \cite{ba}).

For an additive category $A$ we will denote by $A\otimes\q$ its
rational hull i.e. $\obj A\otimes\q=\obj A$  while morphisms are
tensored by $\q$.  $A\q$ will usually (except for some notation of
\S \ref{comphan})  denote the idempotent completion of
$A\otimes\q$.

We will call a realization of motives (usually of $\hk$) {\it
enhanceable} if it has a differential graded enhancement, see \S
\ref{new}.

Other notation will be more or less standard.  $k$ will denote the
ground field; we will assume (except in subsection \ref{motcharp})
that the characteristic of $k$ is zero. $\pt$ is a point, $\af^n$
is the $n$-dimensional affine space (over $k$), $\p^n$ is the
projective space of dimension $n$.

For an additive category $A$ we denote by $C^-(A)$  the category
of complexes over $A$ bounded from above; $C^b(A)\subset C^-(A)$
is the subcategory of bounded complexes;
  $K^-(A)$ is the homotopy category of $C^-(A)$ i.e.
  the morphisms of complexes are considered up to homotopy
equivalence; $K^b$ denotes the homotopic category  of bounded
complexes; sometimes we will also need the unbounded categories
$C(A)$ and $K(A)$; $\ab$ is the category of abelian groups.

 For a category $C,\ A,B\in\obj C$, we denote by
$C(A,B)$ the set of  $C$-morphisms from  $A$ to $B$.

For categories $C,D$ we write $C\subset D$ if $C$ is a full
subcategory of $D$.

We list the main definitions of this paper. Some basic motivic
definitions (mostly coming from \cite{1}) will be given in
\ref{dvoev}. $C(X)$ will be defined in \ref{cusus}; $g^l$ will be
defined in \ref{gl}; differential graded categories, $H(C)$ for a
differential graded category $C$, $S(A)$, $S_N(A)$, $B^-(A)$,
$B^b(A)$, $B(A)$, and $C(A)$ for an additive category $A$ will be
defined in \ref{dg}; the categories of twisted complexes
($\prt(C)$, $\tr(C)$, $\prtp(C)$, $\trp(C)$), arrows, $[P]$ and
$P[i]$ for $P\in\obj C$ will be defined in \ref{tw} and
\ref{prtp}; $\tr(F)$, $\prt(F)$, $\trp(F)$, and $\prtp(F)$ for a
differential graded functor $F$ will be defined in Remark
\ref{remf}; $J$, $\hk$, and $h$ will be defined in \ref{dhk};
$\hk'$, $h'$, $j$, and $J'$ will be defined in \ref{expl}; $C_-$
and different types of truncations of complexes ($\tau_{\le b}$,
$\tau_{[a,b]}$ and the canonical $[a,b]$-truncation) will be
defined in \ref{dop}; $m$ will be defined in \ref{ma}; $\dmgm$
will be described in \ref{stens}; $C_N(P)$, $\hk_N$, and  $t_N$
will be defined in \ref{tn}; $t$, $t\q$, and $\dmge{}'\cong \dmge$
will be defined in \ref{finlength}; truncated realizations will be
defined in \ref{hodge}.

\end{notat}

\section{Cubical Suslin complexes}

In this paper instead of the simplicial Suslin complex $\cu(L(P))$
we consider its cubical version $C(P)$. In this section we prepare
for the proof of the following fact: there exists a differential
graded category $J$ (it will be defined in \S2) whose objects are
the (cubical) Suslin complexes of smooth projective varieties, while
its morphisms are related to the morphisms between those complexes
in $\dme$.

First we recall basic definitions of Voevodsky (along with some
'classical' motivic definitions).

\subsection{Some definitions of Voevodsky: a reminder}\label{dvoev}

We use much of the notation from  \cite{1}. We recall (some of) it
here for the convenience of the reader. Those who remember
Voevodsky's notation well (and agree to identify certain equivalent
categories) could skip this subsection; note only that $\dms$ is the
smallest strict triangulated subcategory of $\dme$ containing all
motives of smooth varieties.


$\var\supset \sv\supset \spv$ will denote the class of all varieties
over $k$, resp. of smooth varieties, resp. of smooth projective
varieties.

 $\smc$ is the category of 'smooth
correspondences' i.e. $\obj \smc=\sv$, $\smc (X,Y)=\bigoplus_U\z$ for all
integral closed $U\subset X\times Y$  that are finite over $X$ and
surjective over a connected component of $X$.

$\ssc=\ssc_{Nis}$ is the abelian category  of  additive cofunctors
$\smc\to\ab$ that are sheaves in the Nisnevich topology (when
restricted to the category of smooth varieties); these sheaves are
usually called 'sheaves with transfers'. Moreover, by default all
sheaves will be sheaves in Nisnevich topology. By an abuse of
notation we will also denote by $\ssc$ the set of all Nisnevich
sheaves with transfers; $\dmk$ is the bounded above derived category of $\ssc$.

For $Y\in \sv$ (more generally, for $Y\in \var$, see \S4.1 of
\cite{1}) we consider $L(Y)=\smc(-,Y)\in \ssc$. For $X\in \sv$ we
denote by $L^c(X)(Y)\supset L(X)(Y)$  the group whose generators are
the same as for $L(X,Y)$ except that $U$ is only required to be
quasi-finite over $X$. $L(X)=L^c(X)$ for proper $X$. Note that
$L^c(X)$ is also a sheaf.

$\mg(X)=\cu(L(X))\cong C(L(X))$ is the Suslin complex of $L(X)$, see
subsection \ref{cusus} below for details; $\mg^c(X)=\cu(L^c(X))\cong
C(L^c(X))$.

$S\in \ssc$ is called homotopy invariant if for any $X\in \sv$ the
projection $\af^1\times X\to X$ gives an isomorphism $S(X)\to
S(\af^1\times X)$.

$\dme\subset \dmk$ is the subcategory of complexes whose cohomology
sheaves are homotopy invariant. It was proved in \cite{1} that for
any $F\in\ssc$ we have $\cu(F)\in \dme$.

The functor $RC:\dmk\to\dme$ is given by taking total complexes of
the Suslin bicomplex of a complex of sheaves, see \S3.2 of \cite{1}
for details.

 $\dms$ will denote the full strict triangulated subcategory of  $\dme$
 generated by
$\mg(X)$ for $X\in \sv$ (we do not add the kernels of projectors).
$\dms$ has a natural tensor structure that can be defined using the
relation $\mg(X)\otimes \mg(Y)=\mg(X\times Y)$; tensor
multiplication of morphisms is defined by means of a similar
relation.

In \cite{1} Voevodsky  defined $\dms$ as a certain localization of
$K^b(\smc)$ (note that he didn't introduce any notation for $\dms$);
then $\dmge$ was defined as the idempotent completion of $\dms$. Yet
Theorem 3.2.6 of \cite{1} (essentially) states that 'his' $\dms$ is
equivalent to those defined here. So we will denote by $\dmge$ the
idempotent completion of 'our' $\dms$.

By definition (cf. subsection 2.1 in \cite{1}) the Tate motif
$\z(1)[2]\in\dme$ can be represented as the cone of the natural map
$\mg(\pt)\to \mg(\p^1)$. Note also that $\mg(\pt)$ is a direct
summand of $\mg(\p^1)$, so $\z(1)[2]$ also is. In particular,
$\z(1)[2]$ is also a Chow motif; see below.

 $\dmgm$ in
\cite{1} was obtained from $\dmge$ (considered as an abstract
category i.e. not as a subcategory of $\dme$) by the formal
inversion of $\z(1)$ with respect to $\otimes$. We will use the same
definition; see \S\ref{stens} 
 below for details.  $\dmgm$ is a rigid
tensor triangulated category. We will also consider the idempotent
completion $\dmgm\q $ of $\dmgm\otimes\q$.


$\cho$ will denote the (homological) category of rational
correspondences. Its objects are smooth projective varieties; the
morphisms are morphisms in $\smc$ up to homotopy equivalence. The
category $\chow$ is the idempotent completion of $\cho$; it was
shown in Proposition 2.1.4 of \cite{1} that $\chow$ is naturally
isomorphic to the usual category of effective homological Chow
motives.

$\chown$ will denote the whole category of Chow motives i.e.
$\chow[\z(-1)[-2]]$. Note that the Chow motif that we denote by
$\z(1)[2]$ (this is compatible with the embedding $\chow\to\dmge$)
was sometimes  ('classically') denoted by $\z(1)$; a similar
convention was used by Hanamura. We will not pay much attention to
this distinction below.

Note  (as it is well known already from the works on Tate motives
that come from quadratic forms) that $\obj\chow\neq \obj\chow\q$
i.e. on the rational level one gets more idempotents in $\cho$.
Certainly, the same is true for $\dmge$ and $\dmge\q$.

We recall also that for categories of geometric origin (for example,
for $\cho$ and $\smc$) the addition of objects is induced by the
disjoint union of varieties operation.

\subsection{The definition of the cubical complex}\label{cusus}

For any    $P\in\sv$ we consider the sheaves
$$C'^i(P)(Y)=\smc(\af^{-i}\times Y,P),\ Y\in \sv;\ C^i=0{\text{
for }}i>0.$$ 
We will usually consider projective $P$.

 By Yoneda's lemma
$$C'^i(P)(Y)\cong \ssc(L(Y),L(P))= \ssc(C'^0(Y),C'^i(P)).$$ For
all $1\le j\le -i$, $x\in k$, we define $d_{ijx}=d_{jx}:C'^i\to
C'^{i+1}$ as $d_{jx}(f)=f\circ g_{jx}$, where
$g_{jx}:\af^{-i-1}\times Y\to\af^{-i}\times Y$ is induced by the map
$(x_1,\dots,x_{-1-i})\to (x_1,\dots,x_{j-1},x,x_j,\dots,x_{-1-i})$.
We define $C^i(P)(Y)$ as $\cap_{1\le j \le -i}\ke d_{j0}$. One may
say that $C^i(P)(Y)$ consists of correspondences that 'are zero if
one of the coordinates is zero'. The boundary maps $\delta^i:C^i\to
C^{i+1}$ are defined as $\sum_{1\le j\le -i}(-1)^jd_{j1}$. Again,
$C^i=0$ for positive $i$.

Since $C'^0=C^0$, we have $C^i(P)(Y)\cong \ssc(C^0(Y),C ^i(P))$.

\begin{rema}\label{sustot}

1. The definition of the cubical Suslin complex can be easily
extended to an arbitrary complex $D$ over $\ssc$ (or over a slightly
different abelian category). One should consider the total complex
of the double complex whose terms are
$$D^{ij}(X)=\cap_{1\le l \le -i}\ke g_{jl0}^*: D^j(\af^{-i}\times
X)\to D^j(\af^{-i-1}\times X),$$ the boundaries are induced by
$\delta^i$.

2. In the usual (simplicial) Suslin complex one defines
$\cu^i(F)(X)=F(D^{-i}\times X)$, where $D^{-i}\subset \af^{1-i}$ is
given by $\sum_{1\le l\le 1-i} x_l=1$; the boundaries come from
restrictions to $x_l=0$.

It is well known that cubical and simplicial complexes do not differ
much; the main advantage of cubical complexes is that descriptions
of (various) products become much nicer (see \S2.5 of \cite{le1}).
\end{rema}

We formulate the main property of $C$.

\begin{pr}\label{prop}
 For any $j\in\z$, $Y\in\sv$, and $P\in\spv$  there is
a natural isomorphism $H^jC(P)(Y)\cong A_{0,-j}(Y,P)\cong
\dme(\cu(Y),\cu(P)[j])$.
\end{pr}
\begin{proof}

 $A_{0,-j}(Y,P)\cong\dme(\cu(Y),\cu^c(P))$ by Proposition 4.2.3
\cite{1}.  Since $P$ is projective, by Proposition 4.1.5 \cite{1} we
have $\cu^c(P)=\cu(P)$.

The first isomorphism  will be described in  \S\ref{prope} below.

All these isomorphisms are natural. \end{proof}

In particular the cohomology presheaves of  $C(P)$ are homotopy
invariant.

We denote the initial object of  $\smc$ by $0$. We define $C^i(0)=0$
for all $i\in \z$.  We obtain
$$p(C(P))\in\obj\dme\subset \obj\dmk,$$ where $p:\kmk\to \dmk$ is
the natural projection.

\subsection{The assignment $g\to (g^l)$}\label{gl}

Let $P,Y$ belong to $\sv$. We construct a family of morphisms $C(Y)\to
C(P)[i]$.

 For any  $f\in C'^i(P)(Y),\ l\le 0$,  we define
$f^l:C'^l(Y)\to C'^{l+i}(P)$ as follows. To the element  $h\in
\smc(Z\times \af^{-l},Y),\ Z\in\sv$, we assign $(-1)^{li} f\circ
(\id_{\af^{-i}}\otimes h) $. It is easily seen that the same formula
also defines  the maps $f^l:C^l(Y)\to C^{l+i}(P)$ for $f\in
C^i(P)(Y)$.

\begin{pr}\label{dmc}

1. The assignment $g\to G=(g^l)$ defines a homomorphism  $\ke
\delta^i(P)(Y)\to\kmk(C(Y),C(P)[i])$.

2. The assignment $g\to G=(g^l)$ induces an isomorphism
$H^i(C(P)(Y))\cong \dme(C(Y),C(P)[i])$.
\end{pr}
\begin{proof}

1. For any $f\in C'^i(P)(Y), h\in C'^l(Y)(Z),\ Z\in\sv$ we have an
equality
\begin{equation}\label{fu}
\delta^{i+l}f^l(h)=(-1)^if^{l+1}\delta^l(h)+(\delta^i f)^l(h).
\end{equation}
 Hence if
$\delta^i g=0,\ g\in C^i(P)(Y)$ then $G$ defines a morphism of
complexes $C(Y)\to C(P)[i]$.

2. Using (\ref{fu}) we obtain that the elements of
$\delta^{i+1}(C^{i+1}(P)(Y))$ give homomorphisms $C(Y)\to C(P)[i]$
that are homotopy equivalent to $0$. Hence we obtain a homomorphism
$H^iC(P)(Y)\to \dme(p(C(Y)),p(C(P)[i]))$. The bijectivity of this
homomorphism will be proved in  \S\ref{prope} below.
\end{proof}

\section{Differential graded categories; the description of $\hk$
and $h:\hk\to K^-(\ssc)$}

Categories of {\it twisted complexes} (defined in subsections
\ref{tw} and \ref{prtp}) were first considered in \cite{bk}. Yet our
notation differs slightly from that of \cite{bk}; some of the signs
are also different.

In subsections \ref{dhk} and \ref{expl} we define and describe our
main categories: $J,\hk,J'$, and $\hk'$.

In subsection \ref{sstufil} we define a natural 'stupid' filtration
on $\hk'$ that is 'close' to those on $C^b(\cho)$; we prove its
natural properties.

We will not need the formalism of \S \ref{dop} till \S \ref{new}. 

\subsection{The definition of differential graded
categories}\label{dg}

Recall that an additive category $C$ is called graded if for any
$P,Q\in\obj C$ there is a canonical decomposition $C(P,Q)\cong
\oplus_i C^i(P,Q)$ defined;
 this decomposition satisfies $C^i(*,*)\circ C^j(*,*)\subset C^{i+j}(*,*)$.
 A differential graded category (cf. \cite{bk} or \cite{drinf}) is a graded category
 endowed with an additive operator
$\delta:C^i(P,Q)\to C^{i+1}(P,Q)$ for all $ i\in \z, P,Q\in\obj C$.
$\delta$ should satisfy  the equalities $\delta^2=0$ (so $C(P,Q)$ is
a complex of abelian groups); $\delta(f\circ g)=\delta f\circ
g+(-1)^i f\circ \delta g$ for any $P,Q,R\in\obj C$, $g\in C^i(P,Q)$,
$f\in C(Q,R)$. In particular, $\delta (\id_P)=0$.

We denote $\delta$ restricted to morphisms of degree $i$ by
$\delta^i$.

For an additive category $A$ one can construct the following
differential graded categories. The notation introduced below will
be used throughout the paper.

We denote the first one by $S(A)$. We set $\obj S(A)=\obj A;\
S(A)^i(P,Q)=A(P,Q)$ for $i=0$;
 $S(A)^i(P,Q)=0$ for $i\neq 0$.
We take $\delta=0$.

We also consider the category $B^-(A)$ whose objects are the same as
for $C^-(A)$, whereas for $P=(P^i)$, $Q=(Q^i)$ we define
$B^-(A)^i(P,Q)=\prod_{j\in \z} A(P^j, Q^{i+j})$.  Obviously $B^-(A)$
is a graded category.

We denote by $B^b(A)$ the full subcategory of $B^-(A)$ whose objects
are  bounded complexes. $B(A)$ and $C(A)$ will denote the
corresponding categories whose objects are unbounded complexes.

We set $\delta f=d_Q\circ f-(-1)^i f \circ d_P$, where $f\in
B^i(P,Q)$, $d_P$ and $d_Q$ are the differentials in $P$ and $Q$.
Note that the kernel of $\delta^0(P,Q)$ coincides with $C(A)(P,Q)$
(the morphisms of complexes); the image of $\delta^{-1}$ are the
morphisms homotopic to $0$.

For any $N\ge 0$ one can define a full subcategory  $S_N(A)$ of
$B^b(A)$ whose objects are complexes concentrated in degrees
$[0,N]$. We have $S(A)=S_0(A)$.

$B^b(A)$ can be obtained from $S(A)$ (or any $S_N(A)$) by means of
the category functor $\prt$ described below.


For any differential graded $C$ we define a category $H(C)$; its
objects are the same as for $C$; its morphisms are defined as
$$H(C)(P,Q)=\ke \delta^0_C(P,Q)/\imm \delta^{-1}_C(P,Q).$$

\subsection{Categories of twisted complexes ($\prt(C)$ and
$\tr(C)$)}\label{tw}

Having a differential graded category $C$ one can construct two
other differential graded categories $\prt(C)$ and $\prtp(C)$ as
well as  triangulated categories $\tr(C)$ and $\trp(C)$. The
simplest example of these constructions is $\prt(S(A))=B^b(A)$.

\begin{defi}
1. The objects of $\prt(C)$ are  $$\{(P^i),\ P^i\in\obj C, i\in\z,
q_{ij}\in C^{i-j+1}(P^i,P^j)\};$$
 here almost all $P^i$ are $0$;
 for any $i,j\in \z$ we have
$\delta q_{ij}+\sum_l q_{lj}\circ q_{il}=0$. We call $q_{ij}$ {\it
arrows} of degree $i-j+1$. For $P=\{(P^i),q_{ij}\}$,
$P'=\{(P^i{}'),q'_{ij}\}$ we set
$$\prt(C)^l(P,P')=\bigoplus_{i,j\in\z}C^{l+i-j}(P^i,P'^j).$$ For
$f\in C^{l+i-j}(P^i,P'^j)$ (an arrow of degree $l+i-j$) we define
the differential of the corresponding morphism in $\prt(C)$ as
$$\delta_{\prt(C)}f=\delta_C f+\sum_m(q'_{jm}\circ
f-(-1)^{(i-m)l}f\circ q_{mi}).$$

2. $\tr(C)=H(\prt(C))$.
\end{defi}

It can be easily seen that $\prt(C)$ is a differential graded
category (see \cite{bk}). There is also an obvious translation
functor on $\prt(C)$. Note also that the terms of the complex
$\prt(C)(P,P')$ do not depend on $q_{ij}$ and $q'_{ij}$, whereas the
differentials certainly do.

We denote by $Q[j]$ the object of $\prt(C)$ that is obtained by
putting $P^i=Q$ for $i=-j$, all other $P^j=0$, all $q_{ij}=0$. We
will write $[Q]$ instead of $Q[0]$ (i.e., $Q[i]$ is the translation
of $[Q]$ by $[i]$).

Immediately from definition we have $\prt(S(A))\cong B^b(A)$.

A morphism $h\in\ke \delta^0$ (a closed morphism of degree $0$) is
called a {\it twisted morphism}. For a twisted morphism
$h=(h_{ij})\in \prt((P^i,q_{ij}),(P^i{}',q'_{ij}))$, $h_{ij}\in
C(P^i,P^j{}')$ we define $\co(h)=(P^i{}'',q_{ij}'')$, where
$P^i{}''=P^{i+1}\oplus P^i{}'$,
$$q"_{ij}=\begin{pmatrix}q_{i+1,j+1}&0\\
h_{i+1,j}&q'_{ij}\end{pmatrix}$$

We have a natural triangle of twisted morphisms
 \begin{equation}\label{dgtri}
P\stackrel{f}{\to} P'\to \co(f)\to P[1],
\end{equation}
the components of the second map are $(0, id_{P'^i})$ for $i=j$ and
$0$ otherwise. This triangle induces a triangle in the category
$H(\prt(C))$.

\begin{defi}
For distinguished triangles in $\tr(C)$ we take  the triangles
isomorphic to those that come from (\ref{dgtri}) for any
$P,P'\in \prt(C)$, $f$ being twisted.
\end{defi}

We summarize the properties of $\prt$ and $\tr$ of \cite{bk} that
are most relevant for the current paper. We have to replace
bounded complexes by complexes bounded from above. Part II4 of the following proposition is 
new.

\begin{pr}\label{mdg}

I $\tr(C)$ is a triangulated category.

II For any additive category $A$ there are natural isomorphisms

1. $\prt (B^-(A))\cong B^-(A)$.

2. $\tr(B^-(A))\cong K^-(A)$.

3. $\prt(B(A))\cong B(A)$.

4. $\tr (S_N(A))\cong B^b(A)$

\end{pr}
\begin{proof}
I See Proposition 1 \S1 of \cite{bk}.

II 1,3. See Lemma II.II.1.2.10 of \cite{le3}.

2. Immediate from assertion II1.

4. We have natural full embeddings $S_0(A)\subset S_N(A)\subset
B^b(A)$. Since $\tr (S_0(A))\cong \tr (B^b(A))\cong B^b(A)$, we
obtain the assertion.
\end{proof}

\subsection{The categories  $\prtp(C)$ and $\trp(C)$}\label{prtp}

We define certain categories $\prtp(C)$ and $\trp(C)$; our constructions yield categories equivalent to those defined in \S4 of \cite{bk}. 
The definitions of $\prtp(C)$ and $\trp(C)$ coulds also be found in
\S2.3 of  \cite{drinf}; there these categories were denoted
by $C^{pre-tr}$ and $C^{tr}$.

\begin{defi}\label{dprtp}
$\trp(C)$ is defined as the full (strict) triangulated category of $\tr(C)$ generated by 
$[P]:\ P\in \obj C$; we denote
the corresponding full subcategory of $\prt(C)$ by $\prtp(C)$.  
\end{defi}

The following statement is an easy consequence of the definitions above.

\begin{pr}\label{gentrp}

1. There are natural embeddings of categories $i:C\to \prtp(C)$ and
$H(C)\to \trp(C)$ sending $P$ to $[P]$.

2. $\prt(i)$, $\tr(i)$,  $\prtp(i)$, and $\trp(i)$ are equivalences
of categories.
\end{pr}

\begin{proof}

1. By definition of  $\prtp(C)$ (resp.  of $\trp(C)$) there exists a
canonical isomorphism of bifunctors $C(-,-)\cong \prtp(C)([-],[-])$
(resp. $HC(-,-)\cong \trp(C)([-],[-])$). In remains to note that
both of these isomorphisms respect addition and composition of
morphisms; the first one also respects differentials.

2. The proof was given in \S3 and \S4 of \cite{bk}.
\end{proof}

\begin{rema}\label{remf}
1. Since $\prt$, $\prtp$, $\tr$, and $\trp$ are functors on the
category of differential graded categories, any differential
category functor $F:C\to C'$ naturally induces functors
 ${\prt}F$, ${\prtp}F$, ${\tr}F$, and ${\trp}F$.
We will use this fact throughout the paper.

For example, for $X=(P^i,q_{ij})\in \obj \prt(C)$ we have
${\prt}F(X)=(F(P^i),F(q_{ij}))$; for a morphism $h=(h_{ij})$ of
$\prt(C)$ we have ${\prt}F(h)= (F(h_{ij}))$. Note that the
definition of ${\prt}F$ on morphisms does not involve $q_{ij}$; yet
${\prt}F$ certainly respects differentials for morphisms.

2. Let $F: \prtp(C)\to D$ be a differential graded functor. Then the
restriction of $F$ to $C\subset \prtp(C)$ (see Proposition
\ref{gentrp}(1)) gives a differential graded functor $FC:C\to D$.
Moreover, since $FC=F\circ i$, we have $\prtp(FC)= \prtp(F)\circ
\prtp(i)$; therefore $\prtp(FC)\cong \prtp(F)$.
\end{rema}

\subsection{Definition of $\hk$ and $h$}\label{dhk}

For   $X,Y,Z\in\spv$, $i,j,l\le 0$, $f\in C^i(X)(Y)$, $g\in
C^j(Y)(Z)$ we have the equality
\begin{equation}\label{e1} (f^j(g))^l=f^{j+l}(g^l).\end{equation}

Hence we can define a  (non-full!) subcategory $J$ of $\bmk$ whose
objects are $[P]=C(P)$, $P\in\spv$, the morphisms are defined as
$$J^i(C(P),C(Q))=\{\bigoplus_{l\le 0}(g^l):\ g\in C^{i}(Q)(P)\},$$
the composition of morphisms and the boundary operators are the same
as for $\bmk$. There is an obvious addition defined for morphisms;
the operation of disjoint union of varieties gives us the addition
on objects. It follows immediately from (\ref{fu}) that $J$ is a
differential graded subcategory of $\bmk$.

 Note that
$J^i(-,-)=0$ for $i>0$; this is a very important property! In
particular, for any $i<0$, $X,Y,Z\in \obj J$, it implies that
$dJ^i(Y,Z)\circ J^0(X,Y)\subset dJ^i(X,Z)$ and $J^0(Y,Z)\circ
dJ^i(X,Y)\subset dJ^i(X,Z)$. This is crucial for the construction of
truncation functors $t_N$ (see subsection \ref{tn} below).  We call
categories that have no morphisms of positive degree {\it negative}
differential graded categories; this property will be discussed in
subsection \ref{dop} below.

 We define $\hk$ as $\tr(J)$. Since $C^l=0$ for $l>0$, we have $\hk=\trp(J)$ (see Remark 6.2.2(2) of \cite{bws}).
Now Proposition \ref{gentrp} implies the following statement
immediately.

\begin{pr}\label{genhk}
$\hk$ is generated by $[P]$, $P\in\spv$, as a triangulated category.
Here $[P]$ denotes the object of $\hk$ that corresponds to
$[P]=C(P)\in \obj J$.
\end{pr}

We consider the functor $h:\hk\to\kmk$ that is induced by the
inclusion $J\to\bmk$.

We also note that any differential graded functor $J\to A$ induces a
functor $\hk\to \trp(A)$.

The definition of $\hk$ implies immediately that $\hk([P],Q[i])=H^i
(C(Q)(P))$ for $P,Q\in\spv$.

\subsection{An explicit description of $\hk$ and $h$}\label{expl}

For the convenience of the reader we describe $\hk$ and $h$
explicitly. Since in this subsection we just describe the category
of twisted complexes over $J$ explicitly, we don't need any proofs
here. 

We define $J'=\prtp(J)$. $J'$ is an {\it enhancement} of $\hk$ (in
the sense of \cite{bk}). The idea is that taking cones of (twisted)
morphisms becomes a well-defined operation in $J'$ (in $\hk$ it is
only defined up to a non-canonical isomorphism).

We describe an auxiliary category $\hk'$. We set $\obj \hk'=\obj J'=\obj
\hk$, whereas $\hk'(X,Y)=\ke \delta^0_{J'}(X,Y)$ for $X,Y\in \obj
\hk'$.

Hence the objects of $\hk'$ are $(P^i,\ i\in\z,f_{ij},\ i<j)$, where
$(P^i)$ is a finite sequence of  (not necessarily connected) smooth
projective varieties (we assume that almost all $P^i$ are $0$),
 $f_{ij}\in C^{i-j+1}(P^j)(P^i)$ for all $m,n\in \z$ satisfy
 the condition
\begin{equation}\label{dif}
\delta^{m-n+1}(P^n)(f_{mn})+\sum_{m<l<n}f^{m-l+1}_{ln}(f_{ml})=0.
\end{equation}
Morphisms  $g: A=(P^i,f_{ij})\to B=(P'^i,f'_{ij})$ can be described
as sets $(g_{ij})\in C^{i-j}(P'^j)(P^i),i\le j $, where the $g_{ij}$
satisfy
\begin{equation}\label{mor}
\delta_{P^j{}'}^{i-j}(g_{ij}) +\sum_{j\ge l\ge i}
f'{}^{i-l}_{lj}(g_{il})=\sum_{j\ge l\ge i}g_{lj}^{i-l+1}(f_{il})\
\forall i,j\in\z.\end{equation} We will assume that $g_{ij}=0$ for
$i>j$.

Note that $g_{ij}=0$ if $P^i=0$ or $P^j=0$.
 Hence the morphisms for any pair of objects in $\hk'$ are
defined by means of a finite set of equalities.

The composition of $g=(g_{ij}):A\to B$ with $h=(h_{ij}):B\to
C=(P'^i,f'_{ij})$ is defined as $$l_{ij}=\sum_{i\le r\le j}
h^{i-r}_{rj}(g_{ir}).$$ 

$\hk'$ has a natural structure of an additive category. The direct
sum of objects is defined by means of a disjoint union of varieties.

The morphisms $g,h: A=(P^i,f_{ij})\to B=(P'^i,f'_{ij})$ are called
homotopic  ($g\sim h$) if there exist  $l_{ij} \in
C^{i-j-1}(P'^j)(P^i), i-1\le j$, such that
\begin{equation}\label{mothom} g_{ij}-h_{ij}=\delta^{i-j-1}_{P'^{j}}l_{ij}+
 \sum_{i-1\le r\le j}
f{'}^{i-r-1}_{rj}(l_{ir})+ \sum_{i\le r\le j+1}
l^{i-r+1}_{rj}(f_{ir}).\end{equation}

Now $\hk$ can be described as a category whose objects are the same
as for $\hk'$, whereas $\hk(A,B)=\hk'(A,B)/\sim$. The translation on
$\hk$ is defined by shifts of indices (for $P^i, f_{ij}$). For
$g=(g_{ij})\in\hk'(A,B)$ its cone  is defined as
$C=Cone(g)\in\obj\hk'$, the $i$-th term of $C$ is equal to
$P''{}^i=P^{i+1}\oplus P'^i$, whereas $$h_{ij}\in
C^{i-j+1}(P^{j+1}\oplus P'^j)(P^{i+1}\oplus
P'^i)=\begin{pmatrix}f_{i+1,j+1}&0\\
g_{i+1,j}&f'_{ij}\end{pmatrix}, i<j;$$ we have obvious natural maps
$B\to C\to A[1]$.

It is easily seen that $\hk$ coincides with the category defined in
\ref{dhk}. We denote the projection $\hk'\to \hk$ by $j$.

Moreover, as in Theorem  4.6 of \cite{hb} one can check (without
using the formalism described above) that $\hk$ with the
structures defined is a triangulated category. 
 One can also check
directly that $P[0]$ for $P\in\spv$ generate $\hk$ as a triangulated
category.

 For $A=(P^i,f_{ij})\in\obj(\hk')$ we define $h'(A)\in \cmk$
as $(C^j_{A},\delta_A^{j}:C^j_{A}\to C^{j+1}_A)\in \cmk$. Here
$C^{j}_A=\sum_{i\le j} C^{i-j}(P^j)$, the component of
$\delta_A^{j}$ that corresponds to the morphism of $C^{i-j}(P^j)$
into $C^{i-j'+1}(P^j{}')$ equals $\delta^{i-j}_{P^j}$ for $j=j'$ and
equals $f^{i-j}_{jj'}$ for $j'\neq j$.

Note that the condition  (\ref{dif}) implies $d_{h'(A)}^2=0$.

Now we define $h'$ on morphisms. For $(l_{ij}):A\to B,\ s\in\z$, we
set $h'(l)_s=\oplus_{i,j}l^{s-i}_{ij}$.

 One can check explicitly that $h'$ induces an exact functor $h:\hk\to\kmk$.

By abuse of notation we denote by $h'$ also the functor $J'\to
\bmk$.

\begin{rema} $P^i$ should be thought about as being 'stratification
pieces' of the motif $A=(P^i,f_{ij})$. In particular, let $Z$ be
closed in $X$, $Z,X\in\spv$, $Y=X-Z$; suppose that $Z$ is everywhere
of codimension $c$ in $X$. If we adjoin $Z(c)[2c]$ to $\obj J$ (see
subsection \ref{summand}) then $\mg(Y)$ can be presented in $\hk$
as $((X,Z(c)[2c]),g_Z)$, where $g_Z$ is the Gysin morphism (see
Proposition 3.5.4 of \cite{1}). See also Proposition \ref{mgc} for a
nice explicit description of the motif with compact support of any
smooth quasi-projective $X$.
\end{rema}

The main distinction of $\hk$ from the motivic category $\db$
defined by Hanamura (see \cite{h}, \cite{hb}, and \S\ref{comphan})
is that
the Bloch cycle complexes (used in the definition of $\db$) are
replaced by the Suslin complexes; we never have to choose
distinguished subcomplexes for our constructions (in contrast with
\cite{h}). Note also that our definition works on the integral level
in contrast with those of \cite{h}.

\subsection{'Stupid filtration' for motives}\label{sstufil}

As we will see several times below, the category $\hk$ (hence also
$\dms$, cf. Theorem \ref{main}) is very close to $K^b(\cho)$. In a
certain sense, $\hk$ has 'stupid filtration' related to those of
$K^b(\cho)$. Certainly, this filtration is only defined on the level
of $\hk'$ (note that the stupid truncation of an object of
$K^b(\cho)$ depens on its lift to $C^b(\cho)$!).

\begin{pr}\label{stufil}

Let $X=(P^i,f_{ij})$ (as in subsection \ref{expl}).

1. For any $a\le b\in \z$ the set $(P^i,f_{ij}:\ a\le i,j\le b)$
gives an object $X_{[a,b]}$ of $\hk'$ (and so also of $\hk$).

2. If $P^i=0$ for $i<a$ then $(\id_{P^i},\ i\le b)$ gives a morphism
$X\to X_{[a,b]}$ (in $\hk'$ and $\hk$).

3. If $P^i=0$ for $i>b$ then $(\id_{P^i},\ i\le b)$ gives a morphism
$ X_{[a,b]}\to X$.

4. If $P^i=0$ for $i<a$ and for $i>c$, $a<b<c$, then we have a
distinguished triangle $ X_{[a,b]}\to X \to X_{[b,c]}$

\end{pr}
\begin{proof}

1. We have to check that the equality (\ref{dif}) is valid for
$X_{[a,b]}$.  Yet all terms of (\ref{dif}) are zero unless $a\le
i\le j \le b$. Moreover, in the case $a\le i\le j \le b$ the terms
of (\ref{dif}) are the same as for $X$. Both of these facts follow
immediately from the negativity of $J$.

2,3: we have to check the condition (\ref{mor}) for these cases;
again this is  obvious by the negativity of $J$.

4. We should check that $X\to X_{[b,c]}$ is homotopy equivalent to
the second morphism of the triangle corresponding to $X_{[a,b]}\to
X$; this easily follows from (\ref{dgtri}) (see also the
corresponding part of subsection \ref{expl}).

\end{proof}

The definition of the stupid filtration and its properties are quite
similar to those described in \S1 of \cite{ha3} (see property (6) in
the end of that section). Note that we only used the fact
that there are no morphisms of positive degree between objects of
$J$. See \cite{bws} for a vast generalization of this observation.

\subsection{Other generalities on differential graded
categories}\label{dop}

We describe some new differential graded categories and differential
graded functors. We will need them in section \ref{new} below.

\subsubsection{Differential graded categories of morphisms}

For an additive category $A$ we denote by $MS(A)$ the category of
morphisms of $S(A)$. Its objects are $\{(X,Y,f):\ X,Y\in \obj A,
f\in A(X,Y) \}$; $$MS_0((X,Y,f),(X',Y',f'))=\{(g,h):\ g\in A(X,X'),\
h\in A(Y,Y'),\ f'\circ g=h\circ g\}.$$ As for $S(A)$, there are no
morphisms of non-zero degree in $MS(A)$; hence the differential for
morphisms is zero.

We denote $\prt (MS(A))$ by $MB^b(A)$. We recall that a twisted
morphism is a closed morphism of degree $0$, i.e., an element of the
kernel of $\delta^0$.

\begin{pr}
1. $MB^b(A)$ is the category of closed morphisms of $B^b(A)$ (cf. \S2.10 of \cite{drinf}). That
means that its objects are $\{(X,Y,f):\ X,Y\in \obj B(A), f\in \ke
\delta^0 (B(A)(X,Y)) \}$, $$MB^b{}^i((X,Y,f),(X',Y',f'))=\{(g,h):\
g\in B(A)^i(X,X'),\ h\in B(A)^i(Y,Y'),\ f'\circ g=h\circ f\}.$$

2. Let $MB(A)$ denote the unbounded analogue of $MB^b(A)$. Then
$\prt(MB(A))\cong MB(A)$.

3. Let $Cone:MB(A)\to B^b(A)$ denote the natural cone functor. Then
the functor $\prt(Cone)$ in naturally isomorphic to $Cone$.
\end{pr}
\begin{proof}

1. Easy direct verification.

2. The proof is very similar to those of part II1 of Proposition
\ref{mdg}. First we note that $\prt(MB^bB(A))\cong MB^b(A)$, then
extend this to the unbounded analogue.

3. Obviously, $Cone$ is a differential graded functor. Hence it
remains to apply part 2 of Remark \ref{remf}.
\end{proof}

We have obvious  differential graded functors $p_1,p_2:MB(A)\to
B(A)$: $p_1(X,Y,f)=X$, $p_2(X,Y,f)=Y$.

\begin{coro}\label{cqua}
 Let $F:J\to MB(A)$ be a differential graded functor.

1. $\prt(F)$ gives a functorial system of closed morphisms
$\prt(p_1\circ F)(X)\to \prt(p_2\circ F)(X)$ in $B(A)$ for $X\in\obj
J'=\obj \hk$.

2. Let $A$ be an abelian category. Suppose that for any $P\in\spv$
the complex $F([P])$ is exact. Then there exists a natural
quasi-isomorphism $\trp(p_1\circ F)(X)\sim \trp(p_2\circ F)(X)$ for
$X\in\obj \hk$.
\end{coro}
\begin{proof}
1. Obvious.

2. We have to show that $\prtp(Cone(F))(X)$ is quasi-isomorphic to
$0$ for any $X\in \obj J'$. We consider the exact functor
$G=\trp(Cone(F))(X)$; it suffices to show that $G=0$. Recall that
$[P],\ P\in \spv$, generate $\hk$ as a triangulated category. Hence
$G([P])= 0$ for any $P\in\spv$  implies that $G=0$.
\end{proof}

\subsubsection{Negative differential graded categories; truncation functors}

We recall that a differential graded category $C$ is called {\it
negative} if $C^i(X,Y)=0$ for any $i>0$, $X,Y\in \obj C$.

Certainly in this case all morphisms of degree zero are closed (i.e.,
satisfy $\delta f=0$). This notion is very important for us since
$J$ is negative.

For any differential graded $C$ there exist a unique 'maximal'
negative subcategory $C_-$ (it is not full unless $C$ is negative
itself!). The objects of $C_-$ are the same as for $C$, whereas
$C_{-,i}(X,Y)=0$ for $i>0$, $=C^i(X,Y)$ for $i<0$, $=\ke
\delta^0(C(X,Y))$ for $i=0$.

Obviously, if $F:D\to C$ is a differential graded functor, $D$ is
negative, then $F$ factors through the faithful embedding $C_-\to
C$.

Suppose that $A$ is an abelian category.

Then zeroth (or any other) cohomology defines a functor $B_-(A)\to
S(A)$.

More generally, we define two versions of the canonical truncation
functor  for $B_-(A)$. We will need these functors in subsection
\ref{hodge} below.

Let $X$ be a complex over  $A$, $a,b\in\z$, $a\le b$. We define
$\tau_{\le b}$ as the complex $$\dots\to X^{b-2}\to X^{b-1}\to \ke
(X^b\to X^{b+1}), $$ here $\ke (X^b\to X^{b+1})$ is put in degree
$b$. $\tau_{[a,b]}(X)$ is defined as $\tau_{\le b}(X)/\tau_{\le
a-1}(X)$ i.e., it is the complex $$X^{a-1}/\ke(X^{a-1}\to X^{a})\to
X^a \to X^{a+1}\to \dots \to X^{b-1}\to \ke (X^b\to X^{b+1}).$$

The canonical $[a,b]$-truncation of $X$ for $a<b$ is defined as
$$X^{[a,b]}=X^a/dX^{a-1} \to X^{a+1}\to \dots \to X^{b-1}\to \ke
(X^b\to X^{b+1}), $$ again $\ke (X^b\to X^{b+1})$ is put in degree
$b$; for $a=b$ we take $H^a(X)$. Recall that truncations preserve
homotopy equivalence of complexes.

\begin{pr}\label{trunc}
1. $\tau_{\le b}$, $\tau_{[a,b]}$, and the canonical
$[a,b]$-truncation define differential graded functors $B_-(A)\to
B_-(A)$.

2. Let $F:J\to B(A)$  be a differential graded functor; we can
assume that its target is $B_-(A)$. We consider the functors
$\tau_{[a,b]}F$ and $F_{[a,b]}$ that are obtained from $F$ by
composing it with the corresponding truncations. Then there exists a
functorial family of quasi-isomorphisms $\trp (\tau_{[a,b]}F)(X)\to
\trp (F_{[a,b]})(X)$ for $X\in\obj \hk$.
\end{pr}
\begin{proof}
1. Note that all truncations give idempotent endofunctors on $C(A)$.

Hence it suffices extend truncations to all morphisms of  $B_-(A)$
and prove that truncations respect $\delta$.

The definition of truncations on morphisms of negative degree is
very easy. The only morphisms in $B_-(A)$ of degree zero are twisted
ones, i.e., morphisms coming from $C(A)$.

It remains to verify that if a given truncation $\tau$ of a morphism
$f=(f_{i}):(X^i)\to (Y^i)$ in $B_-(A)$ is zero then $\tau(\delta
f)=0$.

First we check this for $\tau =\tau_{\le b}$. $\tau_{\le b}f=0$
means that $f(\tau_{\le b}X)=0$ (i.e., the corresponding restrictions
of $f_i$ are zero). Since the boundary maps $\tau_{\le b}X$ into
itself, the definition of $\delta$ for $B(A)$ gives the result.

Now we consider the case $\tau =\tau_{[a, b]}$. $\tau_{[a, b]}(f)=0$
means that $f(\tau_{\le b}X)\subset \tau_{\le a}Y$. Again it
suffices to note that the boundary maps $\tau_{\le b}X$ and
$\tau_{\le a}Y$ into themselves.

The case of canonical truncation can be treated similarly.

 2. The natural morphism $m([P]):\tau_{[a,b]}F([P])\to
F_{[a,b]}([P])$ gives a functor $H:J\to MB(A)$ such that
$p_1(H)=\tau_{[a,b]}F$ and $p_1(H)=F_{[a,b]}$. It remains to note
that $m([P])$ is a quasi-isomorphism for any $P\in\spv$ and apply
part 2 of Corollary \ref{cqua}.
\end{proof}

\begin{rema}\label{prod}
1. Another way to obtain new differential graded categories is to
take 'tensor products' of differential
 graded categories. In particular, one can consider the categories
  $J\otimes J$ and $\tr(J\otimes J)$ (which could be denoted by $\hk\otimes\hk$).

  2. If $H^i(C(X,Y))=0$ for some differential graded category $C$,
  any $X,Y\in\obj C$, $i>0$, then $\trp(C_-)\sim \trp(C)$. Indeed,
  the embedding $C_-\to C$ gives a functor $F:\trp(C_-)\to \trp(C)$.
  Then an easy argument
  shows  that $F$ is a full embedding; see part 3 of Remark \ref{risom} below. Lastly, since $\obj C$
  generates $\trp(C)$ is a triangulated category, one can easily
  prove (by induction) that $F$ gives an equivalence, i.e., that any
  $X\in\obj \trp(C)$ is isomorphic to some $Y\in F_*(\obj
  \trp(C_-))$.

  3. One can also define {\it positive} differential graded
  categories in a natural way. Positive differential graded
  categories 
seem to be related with $t$-structures. Yet we
  will not study this issue in the current paper.

\end{rema}

\section{The main classification result}

In this section we prove the equivalence of $\hk$ and $\dms$. It
follows that  the presentation of a motif as $m(X)$ for $X\in\hk$
could be thought about as being a 'motivic injective resolution'.

Unfortunately, we don't know how to compare Voevodsky's $\dmge$
(or just $\sv$) with $\hk$ 'directly'.

In \S\ref{stens} we define the Tate twist $\otimes \z(1)$ on
$\hk$. We formally invert it; this provides us with a certain
differential graded description of Voevodsky's $\dmgm$.  

\subsection{The equivalence of categories $m:\hk \to \dms$}\label{ma}

We denote the natural functor  $\kmk\to\dmk$ by $p$, denote $p\circ
h$ by $m$.

\begin{theo}\label{main}

$m$ is a full exact embedding of triangulated categories; its
essential image is $\dms$.
\end{theo}
\begin{proof}

Since  $h$ is an exact  functor, so is $m$. Now we check that $m$ is
a full embedding. By part 2 of Proposition \ref{dmc}, $m$ induces an
isomorphism $\dme(m([P]),m(Q[i]))\cong\hk([P],Q[i])$ for
$P,Q\in\spv$, $i\in\z$. Since $[R],\ R\in\spv$, generate $\hk$ as a
triangulated category (see Proposition \ref{genhk}), the same is
true for any pair of objects of $\hk$ (cf. part 3 of Remark
\ref{risom} below).

We explain this argument in more detail.

First we verify that for any smooth projective $P/k$ and arbitrary
$B\in \hk$ the functor $m$ gives an isomorphism
\begin{equation}\label{mors}
\dme(m([P]),m(B))\cong\hk([P],B).
\end{equation}
By Proposition \ref{dmc},  (\ref{mors}) is fulfilled for $B=P'[j]$,
$j\in\z$,
 $P'\in\spv$.
For any distinguished triangle  $X\to Y\to Z\to X[1]$ in $\hk$ the
functor  $m$ defines a morphism of long exact sequences
\begin{equation}\label{hmor1}\begin{CD}
 \to \hk([P],Y)@>>>\hk([P],Z)@>>>
 \hk([P],X[1])\to
 \\@VV{}V @VV{}V@VV{}V\\
\to \dme(m([P]),m(Y))@>>>\dme(m([P]),m(Z))@>>>
\dme(m(P),m(X)[1])\to\end{CD}
\end{equation}
Thus if  $m$ gives an isomorphism in (\ref{mors}) for $B=X[i]$ and
$B=Y[i]$ for $i=0,1$, then $m$ gives  an isomorphism for $B=Z$.
Since objects of the form  $B=[P']$ generate $\hk$ as a triangulated
category,
 (\ref{mors}) is fulfilled for any $B\in \hk$. Hence for all
 $i\in\z$, $B\in \hk$,
we have $\dme(m(P[i]),m(B))\cong\hk(P[i],B)$. For any distinguished
triangle $X\to Y\to Z\to X[1]$ in $\hk$ the functor $m$ defines a
morphism of long exact sequences
 $(\dots\to
 \hk(Z,B)\to\dots)\to (\dots\to
 \dme(m(Z),m(B))\to\dots)$ similar to (\ref{hmor1}). Now the same argument as above
proves that $m$ is a full embedding.

It remains to calculate the 'essential image' $M$ of the map that is
induced by $m$ on $\obj (\hk)$ (we adjoin to $M$ all objects that
are isomorphic to those coming from $\hk$). According to Proposition
\ref{isom} below we have $m([P])=C(P)\cong \mg(P)$.

Since $\hk$ is generated by $[P]$ for $P\in\spv$ as a triangulated
category,
  $M$ is the strict triangulated subcategory of $\dme$ that is generated by all
$\mg(P)$. Since the tensor structure on $\dme$ is defined by means
of the relation $\mg(X)\otimes \mg(Y)=\mg(X\times Y)$ for
$X,Y\in\sv$,
 $M$ is a tensor subcategory of $\dme$. Since
$\mg(P)\in \obj \dms$ for any $P\in\spv$, we have $M\subset \dms$.
It remains to prove that $M$ contains $\dms$. Now, $\z(1)[2]\in\dme$
can be represented as a cone of $\mg(\pt)\to \mg(\p^1)$ (we will
identify $\z(1)$ with $\mg(\z(1))$). 
 Hence $\z(n)\in \obj M$ for any
$n>0$. Since $\mg$ is a tensor functor, if $\mg(Z)\in \obj M$ for
$Z\in\sv$, then $Z(c)[2c]\in \obj M$ for all $c>0$. Now we apply
Proposition 3.5.4 of \cite{1}  as well as the Mayer-Vietoris
triangle for motives (\S2 of \cite{1}); similarly to Corollary 3.5.5
of \cite{1} we conclude that $\obj M$ contains all $\mg(X)$ for
$X\in\sv$ (cf. the remark in \cite{1} that precedes Definition
2.1.1). A more detailed version of this argument will be used in the
proof of Theorem \ref{ttn} below.

Since $\dms$ is the smallest triangulated subcategory of
 $\dme$ containing motives of all smooth varieties, we prove the claim.
\end{proof}

\begin{rema}\label{risom}
1. In order to calculate $\dme (M,M')$ for $M,M'\in \dms$ (using the
explicit description of $\hk$ given in subsection \ref{expl}) in
terms of cycles one needs to know $m\ob(M)$ and $m\ob(M')$ (or the
preimages of their duals).
 See subsection \ref{expmgc} for a nice result in this direction.

 2. Note that we do not construct any comparison functor $\dms\to
 \hk$ (or $\sv\to \hk$) explicitly. Yet if we denote by $\smc_{pr}$ the full
 subcategory of $\smc$ whose objects are smooth projective
 varieties, then we have obvious functors $\smc_{pr}\to K^b(\smc_{pr})\to
 \hk$. Note that the Voevodsky's description of $\dmge$ also gives
 a canonical functor $K^b(\smc_{pr})\to
 \dms\subset\dmge$.

3. Exactly the same argument as used in the (first part of the)
proof in fact immediately proves the following (general and
well-known) statement. Let $r:C\to D$ be an exact functor (i.e.,
$C,D$ are triangulated categories); let $S$ be a class of objects of
$C$ that generate it as a triangulated category. Then $r$ is a full
embedding if and only if it induces isomorphism $C(X,Y[i])\to
D(r(X),r(Y[i]))$ for all $i\in \z$, $X,Y\in S$.

We will use this statement for $S=\{[P],\ P\in\spv\}$, $C=\hk$ or
$\hk\otimes \q$, below.

\end{rema}

\subsection{On motives of (possibly) singular varieties}

Recall (see \S4.1 of \cite{1}) that for any $X\in\var$ (not
necessarily smooth!) there were certain objects $\mg(X)$ and
$\mg^c(X)$ of $\dme$ defined. $\mg(X)$ was called the motif of $X$;
$\mg^c(X)$ was called the motif of $X$ with compact support.

 The following statement follows easily.

\begin{coro}\label{motvar}
For any (not necessarily smooth) variety $X/k$ there exist $Z,Z'\in
\hk$ such that $m(Z)\cong \mg(X)$, $m(Z')\cong \mg^c(X)$.
\end{coro}
\begin{proof}

It is sufficient to verify that $\mg(X), \mg^c(X)\in \dms$. The
proof of this fact is the same as for Corollaries 4.1.4 and 4.1.6 in
\cite{1}. Indeed, for the
 proofs in \cite{1} one doesn't need to add the kernels of projectors.

\end{proof}

\begin{rema}
We obtain that for the Suslin complex of an arbitrary variety $X$
there exists a quasi-isomorphic complex $M$ 'constructed from' the
Suslin complexes of smooth projective varieties; $M\in h_*(\obj
\hk)\subset \obj \kmk$ is unique up to a homotopy.

Moreover, as we have noted in the proof of Lemma \ref{cohy}, the
cohomology of $C(P)$ as a complex of presheaves for any $P\in \spv$
coincides with its hypercohomology (the corresponding fact for
$\cu(P)$ was proved in \cite{1}). Hence the same is true for any
$M\in h_*(\obj \hk)$. Hence for $X\in\sv$ the quasi-isomorphism
$C(X)\to M$ is given by an element of $H^0(M)(X)$.

This result shows that the presentation of a motif as $m(X)$ for
$X\in\hk$ could be thought about as being a 'motivic' analogue of
taking an injective resolution; here the Suslin complexes of smooth
projective varieties play the role of injective objects.
\end{rema}

\subsection{Tensor products, direct summands and non-effective motives
($\dmgm$)}\label{stens}

For any $P\in \spv$ one can in the obvious way define a differential
graded functor $-\otimes [P]$ that maps $[Q]\in \obj J$ to $[Q\times
P]$. Applying $\tr$ to it, one gets an exact functor $-\otimes
[P]:\hk\to \hk$.

More generally, for any $C\in C^b(J^0)$ one can easily define an
exact  functor $-\otimes [C]:\hk\to \hk$ (see Remark 1.14 of
\cite{lesm}; the construction is very similar to the construction of
a tensor product on $K^B(A)$ for a tensor additive category $A$).
This tensor product  will also map $[Q],\ Q\in\spv$, to the
naturally defined complex $[Q\times C]$. We will present the motif
$\z(1)$ as $([\pt]\to [\p^1])[-2]$ and define a functor $\otimes
\z(1)$ using this construction; we will denote this functor by
$-(1)$. For $n\ge 0$ we will denote by $X(n)$ the $n$-th iteration
of $-(1)$ applied to $X\in \obj \hk$. Note that there exists a
functorial isomorphism $X(1)[2]\oplus X\cong X\otimes \p^1$.

Now we define $\hkgm$ as $\inli_{-(1)}\hk$, i.e., it is the category
whose objects are $(X,m):\ X\in \obj \hk,\ m\in \z$;
$$\hkgm((X,m),(X',m')=\inli_i \hk(X(n+i),X'(n'+i))$$ (note that
these groups are only defined in the case $n+i,n'+i\ge 0$). Note
that $-(1)$ can be extended to $\hkgm$; $\hkgm$ is a
triangulated category.

Next, we recall that $\dmge$ is a tensor triangulated category,
and the corresponding functor $\otimes \z(1):\dmge\to \dmge$ is
fully faithful. Following Voevodsky we define $\dmgm$ in a similar
way as $\inli_{-(1)}\dmge$. Besides, recall that $-(1)$ is a full
embedding of $\dmge$ (and hence also of $\dmgm$) into itself: this
is the Cancellation Theorem (see Theorem 4.3.1 of \cite{1} or
\cite{voevc} for the characteristic $p$ case).

Unfortunately, the author does not know how to compare the twist
functors (i.e., $-(1)$) on $\hk$ and $\dmge$ in the general case.
Still, they 'take on the same values' on $P[i],\ P\in \spv,
i\in\z$; this turns to suffice for comparison of $\hkgm$ with
$\dmgm$.

\begin{pr}\label{prgm}
1. The idempotent completion of $\hkgm$ is equivalent to $\dmgm$.

2. The functor $-\otimes \z(1)$ is an autoequivalence of $\hkgm$.
\end{pr}
\begin{proof}

1. First we verify that $-(1)$ is a full embedding of $\hk$ into
itself.

By part 3 of Remark \ref{risom}, to this end it suffices to prove
that for any $P,Q\in\spv$, $i\in\z$ the map $\hk([P],Q[i])\to
\hk(P(1),Q[i](1))$ induced by $-(1)$ is an isomorphism. The latter
fact easily follows from the results of \S8 of \cite{4} (that  the
suspension map is an isomorphism).

Now we recall that by the definition of the tensor product on
$\dmge$ (see \S3 of \cite{1}; note that the product can be
extended to $\dme$) for any $P\in \spv$ we have $m([P](1))\cong
\mg(P)\otimes \z(1)$. This fact along with the Cancellation
Theorem (see above) easily implies that $\obj \hk(i)\approx \obj
\dms(i)$ for any $i\ge 0$. Here by $\approx$ we means that  (via
the isomorphism $m$ of Theorem \ref{main}) any object of $\hk(i)$
is isomorphic to an object of $\dms(i)$ (here the twist is defined
as in $\dmge$) and vice versa. Indeed, both categorical images in
question are generated (as triangulated categories) by $[P](i)$,
$P\in \spv$.

Now, by Lemma \ref{leqlim} below we have $\inli_{-(1)}\dms \cong
\hkgm$. Since idempotent completions commute with the operation
$C\to \inli_{-(1)} C$ in the case when $-(1)$ is a full embedding,
we obtain the claim.

2. Since $-(1)$ is a full embedding of $\hk$ into itself, the same
is true on $\hkgm$. Now, we consider a full (triangulated)
subcategory $H$ of $\hk$ that is obtained by leaving exactly one
element in each class of isomorphic objects in $\hk$ and define
$H_{gm}= \inli_{-(1)}H$. Clearly, $H_{gm}$ is isomorphic to
$\hkgm$, whereas ${-(1)}$ is invertible on $H_{gm}$.

\end{proof}


\begin{rema}\label{rtens}
1. Certainly, all the arguments above can be applied to the
rational hulls of $\hk$ and $\dmge$.

2. Note that $\hkgm$ has a differential graded enhancement: one
can define a differential graded category $J'_{gm}$ with $\obj
J'_{gm}=\obj \hkgm$ and $J'_{gm}((X,m),(X',m')=\inli_i
J'(X(n+i),X'(n'+i))$. Then $\trp(J'_{gm})\cong \hkgm$ (similarly
to Proposition \ref{gentrp}(2)).

3. In order to simplify the definition of $-(1)$ on $\hk$ one can
modify it (i.e., replace it by an equivalent category). To this end
one could add direct summands of  objects to $J$; see
\S\ref{summand} below for more detail.

\end{rema}

Now we prove the lemma used in the proof.

\begin{lem}\label{leqlim}

Let $C$ be a category, let $F,F'$ be  full embeddings of $C$ into
itself. Suppose that for any $i\ge 0$ each object of $F^i(C)$ is
isomorphic to some object of $F'{}^i(C)$ and vice versa; here
$F^i$, $F'{}^i$ denote the corresponding  $i$-th iteration
functors. Then $\inli _F C\cong \inli _{F'}C$.

\end{lem}
\begin{proof}

We can assume that each isomorphism class of objects in $C$
consists of a single object. Then $\obj F^i(C)=\obj F'{}^i(C)$.

It follows that the following functors are well-defined
(invertible) auto-equivalences of $C$: for $i>0$ we consider
$G_i=F'{}^{-i}\circ F^{i}$, for $i\le 0$ we define $G_i=F^{i}\circ
F'{}^{- i}$. Here $F^{l},\ F'{}^{l}$ for $l\le 0$
denote the composition inverses of $F^{-l},\ F'{}^{ -l}$. 

Now, it remains to note that sending $(X,m)\in \inli _F C$ to
$(G_m(X),m)\in \inli _{F'} C$ for all $m\in\z$ defines an
equivalence of the categories in question (one can easily construct
the inverse functor).

\end{proof}


\section{Comparison of $\dmgm\q$  with Hanamura's category of
motives}\label{comphan}

We  prove that $\dmgm\q$  is anti-isomorphic to the Hanamura's
triangulated category of motives.

Prof. Hanamura has kindly informed the author that he has
(independently) obtained an alternative proof of the
anti-equivalence of motivic categories; his proof uses the extension
of the functor constructed in \cite{ha2}. Unfortunately, Hanamura's
proof is not
available for the public (in any form).

\subsection{The plan}

 We will not recall Hanamura's definitions in detail here
(they are rather long) so the reader should consult \S2 of \cite{h}
for the definition of $\db_{fin(k)}$ and $\db(k)$ (see also \S4 of
\cite{hb} and \S4.5 \cite{le1}). One of
 unpleasant properties of Hanamura's construction is that it uses
a certain composition operation for Bloch's complexes (see below)
which is only partially defined. Yet on the target of our comparison
functor the composition is always defined.

Note that Hanamura's category is cohomological (i.e., the functor
$\spv\to \db_{fin(k)}$ is contravariant); so it's natural to
consider contravariant functors $\hk\to \db_{fin(k)}$. For that
reason we will consider the categories $\dbf\subset\dbk$ that equal
$\db_{fin}(k)^{op}\subset \db(k)^{op}$.

 $\dbf$ could
'almost' be described as $\tr(I)$ for a certain differential graded
$I$. The problem is that composition of morphisms in $I$ is only
partially defined; overcoming this difficulty makes the construction
rather complicated.

Moreover, in our definition of $\hk$ we use not necessarily
connected $P\in\spv$. In Hanamura's notation $P$ should be replaced
by $\oplus P^i$, where $P^i$ are connected components of $P$; yet we
will ignore this distinction below.

The proof consists of three parts.

I Construction of a functor $F:\hk\otimes\q \to \dbf$.

II Proof that $F$ is a full embedding.

III Proof that $F$ extends to an equivalence $\dmgm\q\to \dbk$.

\subsection {Construction of a comparison functor
$F$}\label{alter}

First we modify $\hk\otimes \q$ slightly. We define $J\q$ as the
category whose objects are the same as for $J$ while the morphisms
are given by rational alternated cubical Suslin complexes.

These are defined similarly to the alternated cubical Bloch
complexes.

Let $\Sigma_n$ for $n\ge 0$ denote the group of permutations of $n$
elements; it acts on $\af^n$ by permuting coordinates.

We define $$J\q^i([P],[Q])=\q\otimes \{a_i
 (f):\,f\in J^i([P],[Q])\}\subset
\q\otimes J^i([P],[Q]) $$ for all $i\in\z,\ P,Q\in\spv$. Here $a_i$
is the idempotent
$(\sum_{\si\in\Sigma_{-i}}{\operatorname{sgn}}(\si)\si)/(-i)!$,
${\operatorname{sgn}}(\si)$ is the sign of a permutation. One can
easily see that $(a_i)$ is an endomorphism of the complex of
$\q\otimes J([P],[Q])$, i.e., we can consider the boundaries of $J\q$
induced by those of $J$ (note that $J\q$ is a direct summand of
$J\otimes\q$).

We define the composition of morphisms in $J\q$ as in a
factor-category of $J\otimes\q$. It is possible since for any $i,j\in
\z$, $P,Q,R\in\spv$, $f\in J^i(P,Q),$ $g\in J^j(Q,R)$ we have
$a_{i+j}(g\circ f)=a_{i+j}(a_j(g)\circ a_i(f))$ (an easy direct
verification).

An easy standard argument (see Lemma 2.28 of
\cite{le1}) 
shows immediately that the alteration procedure  gives a
quasi-isomorphism $J([P],[Q])\otimes\q\to J\q([P],[Q])$ for any
$P,Q\in\spv$. Hence (by part 3 of Remark \ref{risom}) there exists
an equivalence of categories $G:\hk\otimes\q\to \hk\q$, where
$\hk\q=\tr(J\q)$.


Hence for our purposes it suffices to construct a natural embedding
$H:\hk\q\to \dbf$. As we have said above, $\dbf$ is 'almost equal'
to $\tr(I)$ for a certain differential graded 'almost category' $I$.
Hence it suffices to define a certain differential graded functor
$G:J\q\to I$ and define $H=\tr G$. Indeed, the 'image' of $H$ will
be a subcategory of $\dbf$ (not necessarily full)  which composition
of morphisms is compatible with those of $\hk\q$; we don't have to
care of morphisms and object of $\dbf$ that do not come from
$\hk\q$.

We describe $I$.  Let $\zc^r(X)$ for $X\in\sv$, $r\ge 0$, denote the
Bloch's alternated cubical cycle complex (we use the version of the
definition described in \S1 of \cite{h}). It is defined in the
following way (note that we consider cohomological complexes; that
forces us to reverse arrows in the usual notation).

One first defines a sequence of groups $\zc'^r(X)$ with
$$\zc{}'{}^r{}^i(X)=\sum_{U\subset X\times\af^{-i}}\q,$$ where $U$
runs through all  integral closed subschemes of $X\times\af^{-i}$
 of codimension $r-i$ that intersect faces properly;
here 'faces' mean subvarieties of $X\times \af^{-i}$ defined by
putting some of the last coordinates of $X\times\af^{-i}$ equal to
$0$ or $1$. $\zc^r$ is obtained form $\zc{}'{}^r$ by alteration, see
the beginning of this subsection.

The boundaries $\zc^r{}^i(X)\to \zc^r{}^{i+1}(X)$ are defined as
$\sum_{0\le j\le -i} (g_{j0}{}_*-g_{j1}{}_*) $, $g_{jx}$ are defined
as in subsection \ref{cusus}.

We have a natural map $\zc^r(X)\otimes \zc^s(Y)\to \zc^{r+s}$ (as
complexes) defined by applying $a_{i+j}$ to the corresponding tensor
product of cycles ($i,j$ are the indices of terms; see the text
after Lemma 2.28 in \cite{le1}).

The objects of $I$ are pairs $(P,r):\ P\in\spv,r\in\z$.

Now $I((X,r),(Y,s))=\zc^{\dim Y+r-s}(X\times Y)$. Composition of
morphisms $f\in I^i(X,Y)$, $g\in I^j(Y,Z)$ is defined by the formula
$$g\circ f=pr_{X\times Z\times \af^{-i-j} }(X\times l(Y)\times
Z\times \af^{-i-j})\cap (g\otimes f), $$ where $l:Y\to Y\times Y$ is
the diagonal embedding, $pr$ is the natural
projection. Note that the composition 
is defined only if $g\otimes f$ intersects
$l(Y)\times X\times Z\times \af^{-i-j}$ properly.

A more detailed description of $\zc^r$ (and the discussion of
several other questions relevant for the results of this section)
could be found in \cite{le1}. The reader is strongly recommended to
look at \S2.5 and \S4.3 loc. cit.

The construction of $H$ uses two facts.

(i) For any connected $P,Q\in \spv$, $\dim Q=r$, we have a natural
embedding of complexes $J\q([P],[Q])\to \zc^r (P\times Q)$.

This is obvious from the definition of $\zc^r$.

(ii) For any connected $P,Q,R\in \spv$, $\dim R=s$, the partially
defined composition $\zc^r(P,Q)\times \zc^s(Q,R)\to \zc^s(P,R)$ (cf.
Proposition 4.3 of \cite{hb}) for Bloch complexes restricted to
$J\q(P,Q)\times J\q(Q,R)$ coincides with the map induced by the
composition in $J\q$. This is very easy since the  composition of
morphisms in $J\q$ is exactly the one induced from the composition
of   Bloch complexes (described above). Note  that the composition
of morphisms coming from $J\q$ is always well-defined.

Hence sending $[P]\in\obj J\q$ to $(P,0)\in \obj I$ and embedding
$J\q([P],[Q])\to \zc^r (P\times Q)$ we obtain a differential graded
functor $G:J\q\to I$. This gives $H=\tr G$.

We note that in Hanamura's construction one often has to choose
distinguished subcomplexes (and modify the choice of elements of
Bloch complexes) in order to 'compute'  compositions of arrows in
his categories.  Yet in the 'image' of $H$ this problem never
occurs.

\subsection {$F$ is a full embedding}\label{embed}

We check that $H$ (and so also $F$) is fully faithful.

By part 3 of Remark \ref{risom}, it suffices to check that $H$
induces an isomorphism
$$\hk\q([P],[Q][i])\to\dbf((P,0),(Q,0)[i])$$ for all $P,Q\in\spv$,
$i\in\z$. Since in $\db(k)$ there is Poincare duality (by
definition, see also (4.5) of \cite{hb}), for $Q$ of pure dimension
$n$ we have $$\dbf((P,0),(Q,0)[i])=\dbf((P\times
Q,0),(\pt,n)[i])=Ch^r(P\times Q,-i).$$
Here $Ch^r(P\times Q,-i)$ denotes the rational higher Chow group,
the last equality follows from theorems (4.10) and (1.1) of
\cite{hb}. Since the same is true in $\hk\q$ (see propositions 4.2.3
and 4.2.9 of \cite{1}; cf. also Proposition 12.1 of \cite{frisus}),
the isomorphism is compatible with the maps of corresponding
complexes, we get the claim.

Another way to express the same argument is to say that $G:J\q\to I$
induces a quasi-isomorphism on morphisms (considered as complexes of
abelian groups). This easily implies that $\tr G$ is a full
embedding (even in our situation when $I$ is 'not quite a
category').

\subsection {Conclusion of the proof}

First we check that any object of the full triangulated subcategory
$\dbf^+\subset\dbf$ whose objects are 'positive' diagrams (i.e.,
those  that contain only symbols of the form $(P,r)$, $r\ge 0$) is
equivalent to $\hk\otimes\q$. We should check that any object of
$\dbf^+$ is isomorphic to an object of the form $H(X),\ X\in \obj
\hk\otimes\q$.

We define the 'zeroth part' of $\dbf$ as the full triangulated
subcategory of $\dbf$  generated by $(P,0)$. We note  $\dbf^+$ is
generated by $(P,r),\ r\ge 0$, as a triangulated category. This is
easy to see directly and  also follows immediately from property (6)
of \S1 of \cite{ha3}. Hence it suffices to verify that any $(P,r),\
r\ge 0,\ P\in\spv$ is isomorphic to a certain object of the 'zeroth
part' of $\dbf$.

In fact, $(P,r)$ belongs to the triangulated subcategory of $\dbf$
generated by $(\p^l)\times P$ for $0\le l\le r$.  It suffices to
check that $(P,r)$ is a direct summand of $(P\times\p^r,0)$, its
complement is $(P\times\p^{r-1},0)$.   First we note that morphisms
between 'pure' objects (i.e., between the objects of $\dbf$ of the type $(P,r)$)
by definition equal to the morphisms between corresponding Chow
motives. The last statement is just a well-known property of Chow
motives; it follows easily from the fact that $\chow$ is a tensor
category (note that $\dbf$ also is and the tensor multiplications
are compatible) and $[\p^r]=[\p^{r-1}]\oplus
[{\operatorname{pt}}](r)[2r]$ in $\chow$ (in the covariant
notation).

 A similar direct sum statement  is verified in \S 2.3 of \cite{hb}. The main
 difference is that Hanamura presented $(P,r)$ as a direct summand of $(P\times(\p^1)^r,0)$
 and does not care about its complement. Since we will idempotent
 complete $\dmgm\otimes\q$ in the end of the section, this version also fits  our
 purposes.

Hence we get the equivalence of $\hk\otimes\q$ with  $\dbf^+$.

Note that $\q(1)\in \obj \hk\otimes \q$ differs from
${\operatorname{pt}}(1)$ in $\dbf^+$ by a shift by $[2]$; cf. also the remark
at the bottom of page 139 of \cite{hb}.

Now we prove that this equivalence can be extended to $\dmgm\q$. 

The definition of morphisms in $\dbf$ immediately implies the
cancellation theorem, i.e., $\dbf(X,Y)=\dbf(X(1),Y(1))$. Moreover,
our definition of $-(1)$ on $\hk$ (in \S\ref{stens}) is easily
seen to be compatible the Tate twist on $\dbf$ (again, they differ
by a shift by $[2]$).

 Hence we can extend $F$ to a functor $F':\hkgm
 \otimes\q\to \inli_{-(1)}\dbf^+$; $F'$ is also a full
 embedding. It
remains to recall that $\dmgm\q$ is isomorphic to the idempotent
completion of $\hkgm
 \otimes\q$ (see part 1 of Remark \ref{rtens}),  whereas $\dbk$ is the
idempotent completion of $\inli_{-(1)}\dbf^+=\cup_{n>0}
\dbf^+(-n)$.

The proof is finished.


\begin{rema}\label{rchan}

 Hanamura defined $\db(k)$ as a triangulated subcategory of a
certain 'infinite' analogue of $\db_{fin}(k)$. Yet we do not need
this definition for the proof of equivalence (above) since $\db(k)$
was defined as the idempotent completion of $\db_{fin}(k)$; recall
that the idempotent completion of a triangulated category is
canonical (see \cite{ba}).


\end{rema}

\section{The properties of  cubical Suslin complexes}\label{prope}

The main result of this section is that the cubical complex $C(X)$
is quasi-isomorphic (as a complex of presheaves) to the simplicial
complex $\cu(X)$ that was used in \cite{1}. This fact was mentioned
by M. Levine (Theorem 2.25 of \cite{le1}) yet no complete proof was
given. One of the possible  methods of the proof (proposed in
\cite{le1}) is the use of a bicomplex method. Recall that this
method yielded a similar result for motivic cohomology in \S4 of
\cite{le2} and can be adjusted to yield the proof of the
statement desired. Yet  we use another method here. The reader not
interested in the details of the proof should skip this section.

\subsection{A certain adjoint functor for the derived category of presheaves}

We denote by  $\psc$ the category of presheaves (of abelian groups)
on $\smc$, by $\dpk$ the derived category of $\psc$ (complexes are
bounded from above), by $\dpe$ we denote a full subcategory of
$\dpk$ whose objects are complexes with homotopy invariant
cohomology.

\begin{lem}\label{inv}
$C(P)\in \dpe$.
\end{lem}

\begin{proof}
The scheme of the proof is the same as for Proposition 3.6 in
\cite{3}. First we check that for  $Y\in \sv$ and  $P\in\spv$ the
maps $i_0^*,i_1^*:C(P)(Y\times \af)\to C(P)(Y)$ are homotopic; here
$i_0^*,i_1^*$ are induced by the embeddings $i_x:Y\times\{x\}\to
Y\times{\af},\ x=0,1$. We consider the maps $pr_i:C'(P)^{i}(Y)\to
C'(P)^{i}(Y\times \af)$ induced by the projections $Y\times\af\to
Y$. Also consider the maps $h_i':C'(P)^i(Y\times \af)\to
C'(P)^{i-1}(Y)$ induced by isomorphisms  $Y\times \af\times
\af^{-i}\cong Y\times \af^{-i+1}$ (we put the $\af$ multiplier at the first coordinate of the right hand side), and  $h_i:C(P)^i(Y\times
\af)\to C(P)^{i-1}(Y)$, $h_i=h_i'-pr_{i-1}\circ i_0^*$. We have
$$\delta^{i-1}_* h_i+h_{i+1}\delta^{i}_*=(i_1^*-i_0^*)_i,$$ i.e.,
$h_i$ gives the homotopy needed. Then $i_0^*,i_1^*$ induce
coinciding maps on cohomology. Let  $Y=U\times \af$, $U\in\sv$. We
consider the morphism $H=\id_U\times\mu:U\times\af^2\to U\times
\af$, where $\mu$ is given by multiplication. Since the
maps induced by $H\circ i_0$ and $H\circ i_1$ on the cohomology
of $C(P)(Y)$ coincide, the composition
 $U\times\af\to
U\stackrel{\id_U\times i_0}{\to}U\times \af$ induces an isomorphism
on the cohomology of $C(P)(U\times \af)$ for any $U\in\sv$. Hence  the
cohomology presheaves of $C(P)$ are homotopy invariant.
\end{proof}

Now we formulate an analogue of Proposition 3.2.3 in \cite{1}. By $(F)$
we denote a complex concentrated in degree $0$ whose non-zero term
is $F$.

\begin{pr}\label{fun}

1. There exists an exact functor $R:\dpk\to \dpe$ right adjoint to
the embedding $\dpe\to \dpk$. Besides $R((F))\cong \cu(F)$ (see the
definition of $\cu(F)$ in \S3.2 \cite{1}).

2. In $\dpk$ we have $R(\cu(L(P)))\cong C(P)\cong R(L(P))$.
\end{pr}

\begin{proof}
1. The proof is similar to the proof of existence of the projection
 $RC:\dmk\to\dme$
in 3.2 of \cite{1}. We consider the localising subcategory $\cal{A}$
in $\dpk$ that is generated by all complexes  $L(X\times \af)\to
L(X)$ for $X\in\sv$. As in \cite{1} we have
 $\dpk/{\cal A}\cong  \dpe$ (cf. Theorem
9.32 
of  \cite{vbook}).

Now as in the proof of Proposition 3.2.3 in \cite{1} we should
verify the following statements.

1. For any $F\in \psc$ the natural morphism  $\cu(F)\to (F)$ is an
isomorphism in $\dpk/\cal{A}$.

2. For all $T\in \dpe$ and $B\in \cal{A}$ we have $\dpk(B,T)=0$.

The proof of the first assertion may be copied word for word from
the similar statement in 3.2.3 of \cite{1}.
This was noted in the proof of Theorem 3.2.6 of \cite{1}.

As in \cite{1}, for the second assertion we should check for any
$X\in\sv$ the bijectivity of the map $$\dpk((L(X)), T)\to
\dpk((L(X\times\af)), T)$$ induced by the projection $X\times
\af\to X$. Since
 representable presheaves are projective in  $\psc$ (obvious from Yoneda's lemma, cf.  2.7
in  \cite{vbook}), this follows immediately from the homotopy
invariance of the cohomology of  $\cu(F)$.

2.  From part  2 of Lemma \ref{skl} below we obtain that the
morphism $(L(P))\to C(P)$ induces an isomorphism $R((L(P)))\cong
R(C(P))$ in the category $\dpk$. Using assertion 1 we obtain that
the map $\cu(P)\to (L(P))$ induces an isomorphism $R((L(P)))\cong
R(\cu(P))$. Since $R$ is right adjoint to an embedding of
categories, it remains to note that $\cu(L(P)),C(P)\in\dpe$.
\end{proof}

\begin{lem}\label{skl}
1. $R(C^j(P))=0$ for $j<0$.

2. The morphism $i_P:(L(P))\to C(P)$ induces an isomorphism
$R((L(P)))\cong R(C(P))$ in $\dpk$.
\end{lem}
\begin{proof}

1. We consider the same maps $h_i:C(P)^{i}(Y\times \af)\to
C(P)^{i-1}(Y)$ as in the proof of Lemma \ref{inv}. Obviously $h_i$
is epimorphic, besides $\ke h_i\cong C(P)^i(Y)$. We obtain an exact
sequence
\begin{equation}\label{ex}0\to   C(P)^i(Y)\to C(P)^{i}(Y\times \af)\to
C(P)^{i-1}(Y)\to 0.
\end{equation}

We prove the assertion by induction on $j$. The case $j=-1$ follows
immediately from (\ref{ex}) applied for the case $i=0$. If
$R(C^j(P))=0$ for $j=m$, then $R$ maps $C(P)^m$ and $C''(P)^m$ to
$0$, where $C''(P)^m(Y)=C(P)^m(Y\times \af)$. Applying (\ref{ex})
for $i=m$ we obtain $R(C(P)^{m-1})=0$ (recall that  $R$ is an exact
functor).

2. Follows from assertion 1 immediately.
\end{proof}

Now we recall  (Theorem 8.1 of \cite{4}) that the cohomology groups
of $\cu(L(P))(Y)$ are exactly $A_{0,-i}(Y,P)$. Hence we completed
the proof of Proposition \ref{prop}.

\subsection{Proof of Proposition {\ref{dmc}}}

\begin{lem}\label{cohy}
For all $i\in \z$, $P,Y\in\spv$, the obvious homomorphism
$$\kmk((L(Y)),C(P)[i])\to \dmk((L(Y)),C(P)[i])$$ is bijective.
\end{lem}
\begin{proof}

By definition the homomorphism considered in the map from the
cohomology of $C(P)=C^c(P)$ into its hypercohomology. By Theorem 8.1
of \cite{4} for  $\cu(P)(=\cu^c(P))$ the corresponding map is
bijective. Hence the assertion follows from $\cu(P)\cong C(P)$ in
$\dpk$.
\end{proof}

\begin{pr}\label{isom}

 $i_P:(L(P))\to C(P)$ induces an isomorphism
$ RC((L(P)))\cong C(P)$ in $\dme$.
\end{pr}
\begin{proof}
Literally repeating the argument of the proof of part  1 of Lemma
\ref{skl} we obtain $RC(C^j(P))=0$ for $j<0$. Therefore
$RC((L(P)))\cong RC(C(P))$. It remains to note that
$C(P)\in\obj\dme$.
\end{proof}

Now we finish the proof of Proposition \ref{dmc}. The assignment
$g\to G=(g^l)$ defines a homomorphism $$\kmk((L(Y)),C(P)[i])\to
\kmk(C(Y),C(P)[i]).$$ Hence it is sufficient to verify that the map
$$\dme(RC((L(Y))),C(P)[i])\to \dme(C(Y),C(P)[i])$$ induced by this
homomorphism coincides with the homomorphism induced by the map
$i_P{}_*:RC(L(Y))\cong RC((C(Y)))$. Since $G\circ i_P{}_*=g$, we are
done.

\section{Truncation functors, the length of motives, and $K_0(\dmge)$}

  In  subsection \ref{tn}
  using the canonical filtration of the (cubical) Suslin complex
  we define the {\it truncation
functors} $t_N$.  These functors are new though certain very partial
cases were (essentially) considered in \cite{gs} and \cite{gu}
(there another approaches were used).

 The target of $t_0$
is just  $K^b(\cho)$ (complexes of rational correspondences, see
\ref{dvoev}). In subsection \ref{finlength} we prove that $t_0$
extends to $t:\,\dmge\to K^b(\chow)$ and to $t_{gm}:\dmgm\to
K^b(\chown)$. In \ref{skz} we prove that $t$ induces an
isomorphism $K_0(\dmge)\cong K_0(\chow)$ thus answering the
question of 3.2.4 of \cite{gs}.  In Corollary \ref{newcor} we
extend this result to an isomorphism $K_0(\dmgm)\cong
K_0(\chown)$.

The functors $t$ (Proposition \ref{funt}), $t_{gm}$ (Remark
\ref{nummot}),
 and all $t_N$ (see Theorem \ref{ttn}) are conservative. $t$
induces a natural functor $t_{num}:\dmge\to K^b(Mot_{num}^{eff})$.
Over a finite field $t_{num}\q$ is (conjecturally) an equivalence,
cf. Remark \ref{nummotp}.

We define the {\it length} of a motif: {\it stupid} length in
\ref{motlen}, {\it fine} and {\it rational} length in
\ref{finlength}. The stupid length is not less than the fine length,
the fine length is not less than the rational one. We prove (Theorem
\ref{ttn}) that motives of smooth varieties of dimension $N$ have
stupid length $\le N$; besides $t_N(X)$ contains all information on
motives of length $\le N$. The length of a motif is a natural
motivic analogue of the length of weight filtration for a mixed
Hodge structure.

For a smooth quasi-projective variety $X$ we calculate $m\ob
(\mg^c(X))$ explicitly (in \ref{expmgc}). Using this result we prove
that the weight complex  of Gillet and Soul\'e  for a smooth
quasi-projective variety $X$ can be described as
$t_0(m\ob(\mg^c(X)))$. Next we recall the $cdh$-topology of
Voevodsky and prove this statement for arbitrary  $X\in\var$ (see
\ref{wecompl}). Besides, $t_0(m\ob(\mg(X)))$ essentially coincides
with the functor $h$ described in Theorem 5.10 of \cite{gu}.

In the next section  we will verify that the  weight
 filtration of 'standard'  realizations is closely related to $t_N$;
 the rational length of a motif coincides with the (appropriately defined) length of
 the weight filtration of its singular realization.

\subsection{Truncation functors of level $N$}\label{tn}

For $N\ge 0$ we denote the $-N$-th canonical filtration of $C(P)$ as
a complex of presheaves (i.e., $\dots\to 0\to
C^{-N}(P)/d_PC^{-N-1}(P)\to C^{-N+1}(P)\to\dots \to C^0(P)\to
0\to\dots$) by $C_N(P)$.

We denote by $J_N$ the following differential graded category. Its
objects are the symbols $[P]$ for $P\in\spv$, whereas
$J_N([P],[Q])^i=C_{N}^{i}(Q)(P)$.   The composition of morphisms is
defined  similarly to those in $J$. For morphisms in $J_N$ presented
by $g\in C^{i}(Q)(P)$, $h\in C^{j}(R)(Q)$, we define their
composition as the morphism represented by $h^i(g)$ for $i+j\ge -N$
and $0$ for $i+j<-N$. Note that for $i+j=-N$ we take the class of
$h^{i}(g) \mod d_RC^{-N-1}(R)(P)$; for $i=-N$, $j=0$, and vice
versa, $g$ is only defined up to an element of $d_QC^{-N-1}(Q)(P)$
(resp. $h$ is defined up to an element of $d_RC^{-N-1}(R)(Q)$) yet
the composition is well-defined. The boundary on morphisms is also
defined as in $J$, i.e., for $g\in J_N(P,Q)$ we define $\delta g= d_Q
g$. Certainly, all $J_N$ are negative (i.e., there are no morphisms
of degree $>0$).

We have an obvious functor $J\to J_N$. As noted in Remark
\ref{remf}, this gives canonically a functor $t_N: \hk\to \tr(J_N)$.
We denote $\tr(J_N)=\trp(J_N)$ by $\hk_N$; note that $\hk_0$ is
precisely $K^b(\cho)$.

For any $m\le N$ we also have an obvious functor $J_N\to J_m$. It
induces a functor $t_{Nm}:\hk_N\to \hk_m$ such that $t_m=t_{Nm}\circ
t_N$.

Certainly, one can give a description of $\hk_N$ that is similar to
the description of $\hk$ given in subsection \ref{expl}. Hence
 objects of $\hk_N$ can be represented as certain $(P^i,
f_{ij}\in C^N{}^{i-j+1}(P^j)(P^i), i<j\le i+N+1)$, the morphisms
between $(P^i, f_{ij})$ and $(P^i{}', f'_{ij})$ are represented by
certain $g_{ij}\in C_N^{i-j}(P^j{}',P^i),\ i\le j\le i+N$, etc. The
functor $t_N$ 'forgets' all elements of $C^m([P],[Q])$ for $P,Q\in
\spv$, $m<-N$, and factorizes $C^{-N}([P],[Q])$ modulo coboundaries.
In particular, for $N=0$ we get ordinary complexes over $\cho$.

\subsection{The 'stupid' length of motives (in $\dms$); conservativity of $t_0$}\label{motlen}

It was proved in \cite{1} that the functor $\mg$ gives a full
embedding of $\cho\to\dme$. In this subsection we prove a natural
generalization of this statement.

We will say that $P=(P^i,f_{ij})\in\obj\hk'$  is concentrated in
degrees $[l,m]$, $l,m\in\z$, if  $P^i=0$ for $i<l$ and $i>m$. We
denote the corresponding additive set of objects of $\hk'$ by
$\hk'_{[a,b]}$. We denote  by $\hk_{[a,b]}$ the objects of $\hk$
that are isomorphic to those coming from $\hk'_{[a,b]}$.

Obviously (from the description of distinguished triangles in
$\tr(C)$ for any differential graded $D$) if $A\to B\to C\to A[1]$
is a distinguished triangle, $A,C\in \hk_{[a,b]}$, then $B\in
\hk_{[a,b]}$.

\begin{theo}\label{ttn}

1. For any smooth variety $Y/k$ of dimension $\le N$ we have
$m\ob(\mg(Y))\in \hk_{[0,N]}$.

2. For any smooth variety $Y/k$ of dimension $\le N$ we have
$m\ob(\mg^c(Y))\in  \hk_{[-N,0]}$.

3. If $A\in \hk_{[a,b]}$,
 $B\in \hk_{[c,d]}$, $N\ge d-a$, $N\ge
0$, then $\hk(A,B)\cong \hk_{N}(t_N(A),t_N(B))$.

4. If $s\in \hk(X,X)$ for $X\in \obj \hk$ is an idempotent and
$t_0(s)=0$ then $s=0$.

 5. $t_0$ is conservative, i.e., for $Y \in \hk$ we have $Y=0\iff t_0(Y)=0$.

6. $f:A\to B$ is an isomorphism if and only if $t_0(f)$ is.

\end{theo}
 \begin{proof}

1. Obviously, the statement is valid for smooth projective  $Y$. We
prove the general statement by induction on dimension.

 By the projective bundle theorem (see Proposition 3.5.3 of
\cite{1}) for any $c\ge 0$ we have a canonical isomorphism
$\p^c\cong \oplus_{0\le i\le c}\z(i)[2i]$.  Hence $\z(c)[2c]$ can be
represented as a cone of the natural map
$\mg(\p^{c-1})\to\mg(\p^c)$. Therefore $\z(c)[2c]\in
m(\hk_{[-1,0]})$.

One can easily show that for any $X\in \hk_{[e,f]}$, $e,f\in\z$,
$c>0$ we have $X(c)[2c]\in \hk_{[e-1,f]}$ (here we consider
Voevodsky's Tate twist). Indeed, if $X$ can be decomposed into a
Postnikov tower whose 'factors' are $P^i$ (the terms of the weight
complex $t(X)$) then $X(c)[2c]$ can be decomposed into a
Postnikov tower whose 'factors' are $P^i(c)[2c]$; since $t$ is an
exact functor, we obtain that $t(X(c)[2c])$ is homotopy equivalent
to a complex concentrated in degrees $[e-1,f]$. 

  We recall the Gysin distinguished triangle (see Proposition 3.5.4 of
\cite{1}). For a closed embedding $Z\to X$, $Z$ is everywhere of
codimension $c$, it has the form
\begin{equation}\label{gys}
\mg(X-Z)\to \mg(X)\to \mg(Z)(c)[2c]\to \mg(X-Z)[1].
\end{equation}

Suppose that the assertion is always fulfilled for $\dim Y=N'<N$.

Let $X/k$ be a smooth quasi-projective variety. Since $k$ admits resolution of
singularities, $X$ can be represented as a complement to a
$P\in\spv$ of a divisor with normal crossings $\cup_{i\ge 0}Q^i$.
Then using (\ref{gys}) one proves by induction on $j$ that the
assertion is valid for all $Y^j=P\setminus(\cup_{0\le i\le
j}Q^i)$. To this end we check  by the inductive assumption for $j\ge
0$ that $$\mg (P\setminus(\cup_{0\le i\le j}Q^i)\setminus
(P\setminus(\cup_{0\le i\le
j+1}Q^i)))=\mg(Q^{i+1}\setminus(\cup_{0\le i\le j}Q^j))  \in
m(\hk_{[0,N-1]}).$$ Hence  $\mg(X)\in m(\hk_{[0,N]}).$

If $X$ is not quasi-projective we can still choose closed $Z\subset
X$ (of codimension $>0$) such that $X-Z$ is quasi-projective. Hence
the assertion follows from the inductive assumption by applying
(\ref{gys}).

2. The proof is similar to those of the previous part. The
difference is that we don't have to twist and should use the
distinguished triangle of  Proposition 4.1.5 of \cite{1}:
\begin{equation}
 \mg^c(Z)\to \mg^c(X)\to \mg^c(X-Z)\to \mg^c(Z)[1]
\end{equation}
 instead of (\ref{gys}).

3. We can assume (by increasing $d$ if needed) that $N=d-a$.

Let $A=(P^i,f_{ij})$, $B=(P'^i,f'_{ij})$. As we have seen in
subsection \ref{expl}, any $g\in \hk(A,B)$ is given by a certain set
of $g_{ij}\in C^{i-j}(P'^j)(P^i),i\le j $. The same is valid for
$h=(h_{ij})\in \hk_N(t_N(A),t_N(B))$; the only difference is that
$h_{ad}$ is given modulo $d_{P'_d}C^{-N-1}(P'^d)(P^a)$. Both
$(g_{ij})$ and $(h_{ij})$ should satisfy the conditions
\begin{equation}\label{morn}
\delta_{P^j{}'}^{i-j}(m_{ij}) +\sum_{j\ge l\ge i}
f'{}^{i-l}_{lj}(m_{il})=\sum_{j\ge l\ge i}m_{lj}^{i-l+1}(f_{il})\
\forall i,j\in\z.\end{equation}

First we check surjectivity. We recall that the conditions
(\ref{morn}) for $g$ depend only on $g_{ij}$ for $(i,j)\neq (a,d)$
and on $d_{P'^d}g_{ad}$. Hence if $(h_{ij})$  satisfies the
conditions (\ref{morn}) then $h=t_N(r)$, where $r_{ij}=h_{ij}$ for
all $(i,j)\neq (a,d)$, $r_{ad}$ is an arbitrary element of
$C^{-N}(P'^d)(P^a)$ satisfying $r_{ad}\mod
d_{P'^d}C^{-N-1}(P'^d)(P^a)=h_{ad}$.

Now we check injectivity. Let $t_N(g)=0$ for $g=(g_{ij})\in
\hk'(A,B)$. Note that $C_N(P)$ is a factor-complex of $C(P)$ for any
$P\in \spv$. Hence similarly to subsection \ref{expl} one can easily
check that there exist $l_{ij} \in C^{i-j-1}(P'^j)(P^i), i\le j$,
such that
\begin{equation}\label{homotop}
g_{ij}=\delta^{i-j-1}_{P'^{j}}l_{ij}+
 \sum_{i\le r\le j}
(f{'}^{i-r-1}_{rj}(l_{ir})+ l^{i-r+1}_{rj}(f_{ir}))
\end{equation}
 for all $(i,j)\neq (a,d)$, for $i=a$, $j=d$ the equality
 (\ref{homotop}) is fulfilled modulo $d_{P'^d}q$ for some $q\in
 C^{-N-1}(P'^d)(P^a)$. We consider $(l'_{ij})$, where $l'_{ij}=l_{ij}$
for all $(i,j)\neq (a,d)$, $l'_{ad}=l_{ad}+q$. Obviously, if we
replace $(l_{ij})$ by $(l'_{ij})$ then (\ref{homotop}) would be
fulfilled for all $i,j$. Therefore $g=0$ in $\hk(A,B)$.

 4. Let $X=(P^i,f_{ij})$ be as in subsection \ref{expl}; let $s$ be
given by a set of $s_{ij}\in J^{i-j}([P^i],[P^j])$, $i\le j$;
$(s_{ij})$ are defined up to a homotopy of the sort described in
\ref{expl}.
 $t_0(s)$ is homotopic to zero. Since this homotopy can be represented by a set of
$m_i\in \smc(P^i,P^{i-1})=J^{0}([P^i],[P^{i-1}])$, we can lift this
homotopy to $\hk$. This means that we take $l_{ij}=m_i$ for $j=i-1$,
$l_{ij}=0$ for $j\neq i-1$, where $l$ is as in (\ref{mothom}); this
allows us to assume that $s_{ii}=0$.  Next, since $s^2=s$ in $\hk$,
we have $s^n=s$ for any $n>0$. Now note that all degrees of
components of $s^n$ are $\le -n$. Hence $s^r=0$ if $X\in
\hk_{[a,b]}$ for $r>b-a$.

5. Immediate from assertion 4 applied to $X=Y$ and $s=\id_Y$.

6. Follows immediately from assertion 5 (recall that a morphism is
an isomorphism if and only if its cone is zero).

\end{proof}

In fact,  for a smooth quasi-projective $X$ one can compute
$\mg^c(X)$ explicitly (see Proposition \ref{mgc} below).

We say that $X\in \hk$ has {\it stupid length} $\le N$ if for some
$l\in\z,m\le l+N$, the motif $X\in \hk_{[l,m]}$. We will define the
{\it fine length} of a motif below.

\begin{rema}

1. In fact, surjectvity (but not injectivity) in part 3 is also
valid for $d-a=N+1$. The proof is similar to those for the case
$d-a=N$. We should choose $r_{ad}$, $r_{a+1,d}$, and $r_{a,d-1}$;
the classes of $r_{a+1,d}$ and $r_{a,d-1}$ modulo coboundaries are
fixed. This choice affects the equality (\ref{morn}) only for $i=a,\
j=d$. Note also that this equality only depends on $d_{P'^d}r_{ad}$.
One can choose arbitrary values of $r_{a+1,d}$ and $r_{a,d-1}$ in
the corresponding classes. Then the equality (\ref{morn}) with
$r_{ad}=0$ will be satisfied modulo $d_{P'_d}q$ for some $q\in
 C^{-N-1}(P'^d)(P^a)$. Therefore if we take $r_{ad}=q$ then $t_N(r)=h$.

2.  Let $A\ \in \hk_{[a,b]}$, $B\in \hk_{[c , d]}$, $N+1\ge d-a$,
$N\ge 0$. Then one can check that $A\cong B$ iff $t_N(A)\cong
t_N(B)$.

Indeed, if $f:A\to B$ is an isomorphism then $t_N(f)$ also is.

Conversely, let $f_N:t_N(A)\to t_N(B)$ be an isomorphism. Then, as
was noted above,  there exists an $f\in \hk(A,B)$ such that
$f_N=t_N(f)$ ($f$ is not necessarily unique). Since $t_N(F)$ is an
isomorphism, $t_0(f)=t_{N0}(t_N(f))$ also is. From part 6 of Theorem
\ref{ttn} we obtain that $f$ gives an isomorphism $A\cong B$.

It follows immediately that two objects $A,B\in\hk $ of stupid
length $\le N+1$ are isomorphic if and only if $t_N(A)\cong t_N(B)$.

3. One could define $\hk_{N,[0,N]}\subset \hk_N$ similarly to
$\hk_{[0,N]}$. Then $t_N$ would give an equivalence of additive
categories $\hk_{[0,N]}\to \hk_{N,[0,N]}$. Indeed, this restriction
of $t_N$ is surjective on objects; it is an embedding of categories
by part 3 of Theorem \ref{ttn}.

\end{rema}

\subsection{{\it Fine} length of a motif (in $\dmge$); conservativity of the weight
complex functor $t:\dmge\to K^b(\chow)$}\label{finlength}

One can check (using the method of the proof of Proposition
\ref{idempdm} below) that the stupid length of a motif $M\in\obj
\dms$ coincides with the length of $t_0(M)\in K^b(\cho)$. Yet
replacing $K^b(\cho)$ by $K^b(\chow)$ one can obtain a more
interesting invariant which will be defined on the whole $\obj
\dmge$.

\begin{pr}\label{funt}
1. $t_0$ can be extended to an exact functor $t:\dmge{}'\to
K^b(\chow)$, where $\dmge{}'$ is the idempotent completion of $\hk$.

2. $t$ is conservative. \end{pr}

\begin{proof}
1. 
By the main result of \cite{ba}, $t_0$ can be canonically extended
to an exact functor from $\dmge{}'$ to the idempotent completion of
$K^b(\cho)$. It remains to note that the idempotent completion of
$K^b(\cho)$ is exactly $K^b(\chow)$ (see, for example, Corollary
2.12 of \cite{ba}).

2. Immediate from part 4 of Theorem \ref{ttn}.

\end{proof}

\begin{rema}\label{nummot}

1. Exactly the same arguments yield that there exists a functor
$t_0\otimes\q:\hk\otimes\q\to K^b(\cho\otimes \q)$ that can be
extended to a conservative $t\q:\dmge{}'\q\to K^b(\chow \q)$.

2. Since $\hk\cong \dms$, we have $\dmge{}'\cong \dmge$. For this
reason we will usually identify these categories.

3. As was shown in \S\ref{stens}, $\dmgm$ is isomorphic to the
idempotent completion of a certain  triangulated category $\hkgm'$
that has a differential graded enhancement. Now, note that the
complex $J'_{gm}(P(i)[2i],Q(j)[2j])$ (in the notation of part 2 of
Remark \ref{rtens}) is acyclic in positive degrees for any
$P,Q\in\sv$, $i,j\in \z$ (indeed, since the Tate twist is
invertible, we can assume $i,j\ge 0$). Hence (replacing $\hkgm$ by
an equivalent category again) we obtain that $\dmgm$ is equivalent
to $\tr(I)$ for a certain negative differential graded $I$ (see
part 2 of Remark \ref{prod}).

Now recall that $\chow$ is a tensor category and the functor
$\otimes \z(1)[2]:\chow\to \chow$ is fully faithful. $\chown$ is
usually  defined as the 'union' of $\chow(-i)[-2i]$, whereas each
$\chow(-i)[-2i]$ is isomorphic to $\chow$, i.e., $\chown
=\inli_{-\otimes \z(1)[2]}\chow$.
 We obtain that $HI\subset\chown$. It
follows that one can  (essentially) extend $t$ to an exact
functor $t_{gm}:\dmgm\to K^b(\chown)$. Since $t$ is conservative,
$t_{gm}$ also is.

4. Composing $t$ (or $t\q$) with the natural functor $\chow\to
Mot_{num}^{eff}$ (resp. $\chow\q\to Mot_{num}^{eff}\q$)
 one gets a
functor $t_{num}:\dmge\to K^b(Mot_{num}^{eff})$ (resp.
$t_{num}\q:\dmge\q\to K^b(Mot_{num}^{eff}\q)$). These functors can
be easily extended to functors $\dmgm\to K^b(Mot_{num})$ and
$\dmge\q\to K^b(Mot_{num}\q)$. Here $Mot_{num}$ (and
$Mot^{eff}_{num}$) denotes the category of (effective) numerical
motives; note that $Mot_{num}\q$ is an abelian category (see
\cite{ja}). Yet in order to obtain the 'correct' structure of
$K^b(Mot_{num}\q)$ (i.e., those compatible with 'standard'
realizations) one needs  Kunneth projectors for numerical motives;
currently they are known to exist only over a finite field (see
\S\ref{motcharp}). Moreover, one cannot prove that
$t_{num}\q$ is conservative without assuming certain 'standard'
conjectures (cf. Proposition \ref{becon}). See also part 2 of Remark
\ref{nummotp} below for further discussion of the case of a finite $k$.

Note also: if numerical equivalence of cycles coincides with
homological equivalence (a standard conjecture!) then the cohomology
of $t_{num}\q(X)$ computes the (pure) weight factors of \'etale (and
singular) cohomology of $X$; see subsection \ref{clawe} below. This
generalizes to (weight complexes of) motives Remark 3.1.6 of \cite{gs}.

\end{rema}

\begin{defi}
1. For $M\in \obj \dmge$ we write  $M\in \dmge{}_{[a,b]}$ if
$t(M)\cong W$ for a complex $W$ of Chow motives concentrated in
degrees $[a,b]$.

2. We define the {\it fine} length of $M$ as  the smallest
difference $b-a$ such that $M\in \dmge{}_{[a,b]}$.

3. The {\it rational length} of a motif $M$ is the length of
$t(M)\otimes\q\in K^b(\chow \q)$ ($\chow\q$ is the idempotent
completion of $\chow\otimes\q$!).

\end{defi}

Obviously, the fine length of $M\in\obj \dms$ is not greater than
its stupid length and not less than its rational length. Note also
that the fine length of $\z(n)$ is $0$.  Besides, we also  have the
sets $\dmge{}_{[a,b],\q}\subset\obj \dmge\q$ for each $a\le b\in
\z$.

 The length of a motif is a natural motivic analogue of the
length of weight filtration of a mixed Hodge structure or of a
geometric representation (i.e., of a representation coming from the
\'etale cohomology of a variety). Note that even  the stupid length
of the motif of a smooth variety is not larger than its dimension;
this is a motivic analogue of the corresponding statement for the
weights of singular and \'etale cohomology.

 The results of the next section along with Proposition \ref{idempdm}
 easily imply (see subsection \ref{clawe}) that the length of the weight filtration of
the singular  or \'etale cohomology of a motif $M$ is not greater
than the fine length of $M$. Moreover, the standard conjectures
imply that $M\in \dmge{}_{[a,b],\q}$ iff for any $l$ the weights of
$H^l(M)$ (singular or \'etale cohomology) lie between $l+a$ and
$l+b$ (see \S \ref{clawe} below).

\subsection{The study of $K_0(\dmge)$ and $K_0(\dmgm)$}\label{skz}

We recall some standard definitions (cf. 3.2.1 of \cite{gs}). We
define the Grothendieck group $K_0(\chow)$ as the Abelian group whose
generators are of the form $[A],\ A\in\obj \chow$; the relations are
$[C]=[A]+[B]$ for $C\cong A\oplus B\in\obj \chow$. Note that
$A\oplus 0\cong A$ implies $[A]=[B]$ if $A\cong B$. We use the same
definition for $K_0(\chown)$.

The $K_0$-group of a triangulated category $T$ is defined as the Abelian
group whose generators are $[t],\ t\in \obj T$;  if $A\to B\to C\to
A[1]$ is a distinguished triangle then $[B]=[A]+[C]$. Note that this
also immediately implies $[A]=[B]$ if $A\cong B$.

The existence of $t$ allows us to calculate $K_0(\dmge)$ easily. To
this end we prove the following statement.

\begin{pr}\label{idempdm}
For $X\in\obj\dmge$ if $X\in \dmge{}_{[a,b]}$ then $X$ is a direct
summand of some $X'\in\hk_{[a,b]}$.
\end{pr}
\begin{proof}
The proof is similar (and generalizes) those of part 4 of Theorem
 \ref{ttn}.

Suppose that  $X$ is a direct summand of an object $Y=(P^i,f_{ij})$
of $\hk$ of length $r$; $X$ is given by an idempotent endomorphism
$s=(s_{ij})$ of $Y$.
  It suffices to verify that $Y$ is a direct
summand of $X'=Y_{[a,b]}$ (see Proposition \ref{stufil}).

Similarly to the proof of part 4 of Theorem
 \ref{ttn}, we can assume that $s_{ii}=0$ for $i>b$ and
$i<a$. Indeed, $s_{ii}=0$ is 'homotopic to $0$ outside of $[a,b]$',
whereas we can alter $s_{ii}$ by $ f^{0}_{i-1i}(l_{ii-1})+
l^{0}_{i+1i}(f_{ii+1})$ for any set of $l_{kk-1}\in
C^0(P^k,P^{k-1})$ (see (\ref{mothom})).

Again we have $s=s^n$ for any $n$. If we replace $s$ by $s^{r+1}$ we
easily obtain that $s_{ij}=0$ for $i>b,j>b$ and $i<a,j<a$.

Obviously, it suffices to check that (for the new choice of $s$):

1.) Arrows $(s_{ij},\ a\le i\le b)$ give a morphism $s'\in
\hk(Y_{[a,b]},Y)$.

2.) $(s_{ij},\ a\le j\le b)$ give a morphism $s''\in
\hk(Y,Y_{[a,b]})$.

3.) $s=s''s'$ (in $\hk$).

For 1.): We have to check all equalities (\ref{mor}) for $s'$.
 Both sides of (\ref{mor}) belong to $J^{j-i+1}([P^i],[P^j])$
  for some $i,j, a\le i\le b$. Since the components of $s'$ are taken
from $s$, we only have to compare the differences for both sides for
$i,j\ge a$ (for other values of $i,j$ both sides are zero for $s'$).
The only summands in (\ref{mor}) for $i,j$ that distinguish $s$ from
$s'$ are those of the form $s_{lj}^{i-l+1}(f_{il})$ for $l>b$. Yet
in these cases $s_{lj}=0$ by our assumption.

2.) This is proved similarly.

3.) In $J'$ we have $s^2=s's''+(s''+s')d+d(s''+s')+ d^2$, where the
components of $d$ are morphisms in $J([P^i],[P^j])$ for $i<a,j>b$.
It remains to note that all (possible) degrees of arrows in
$(s''+s')d+(s'+s'')d$ are positive and $d^2=0$.
\end{proof}

\begin{theo}\label{kze}
 $t$ induces an isomorphism $K_0(\dmge)\cong
K_0(\chow)$.\end{theo}

\begin{proof}
 Since $t$ is an exact functor, it gives an abelian group
homomorphism $a:K_0(\dmge)\to K_0(K^b(\chow))$. By Lemma 3 of 3.2.1
of \cite{gs}, there is a natural isomorphism $b:K_0(K^b(\chow))\to
K_0(\chow)$. The embedding $i:\chow\to \dmge$ (see Proposition 2.1.4
of \cite{1}) gives a homomorphism $c:K_0(\chow)\to K_0(\dmge)$. The
definitions of $a,b,c$ imply immediately that $b\circ a\circ
c=\id_{K_0(\chow)}$. Hence $a$ is surjective, $c$ is injective.

It remains to verify that $c$ is surjective.

We claim that if  $t(X)=P^i\to P^{i+1}\to\dots \to P^j$, $P^l\in
\obj \chow$,
 then the class $[X]\in K_0(\dmge)$ equals $\sum (-1)^l[P^l]$.

We prove this fact by induction on the (fine) length of $X$. The
length one case follows immediately from Proposition \ref{idempdm}
(and the conservativity of $t$ can be considered as a partial case of this statement).

To make the inductive step it suffices  to show the existence of a
morpihsm $l:P^j[-j]\to X$ in $\dmge$ that gives the obvious morphism
of complexes after we apply $t$ to it. Indeed, then the length of
the cone of $l$ would be less than $j-i$ (cf. part 4 of Proposition
\ref{stufil}).

By Proposition \ref{idempdm}, $X$ is a direct summand of some $X'\in
\hk_{[i,j]}$. Hence the existence of $l$ follows from part 3 of
Theorem \ref{ttn}.

\end{proof}

\begin{coro}\label{newcor}
$t_{gm}$ (defined in part 3 of Remark \ref{nummot}) induces an
isomorphism $ K_0(\dmgm)\to K_0(\chown)$.

\end{coro}
\begin{proof}

Note that the definitions of the Tate twist on $\hk$ and $\chow$
is compatible. Since $-(1)$ is an autoequivalence on $\dmgm=\cup
\dmge(-i)[-2i]$ and on $K^b(\chown)=\cup K^b(\chow)(-i)[-2i]$, we
obtain the claim (from Theorem \ref{kze}).

\end{proof}


\begin{rema}\label{rkring}
1. Note that the categories $\chow\subset \dmge$ have compatible
tensor categories structures. Hence their $K_0$-groups are
actually rings, whereas the isomorphism constructed is an
isomorphism of rings.

The same is true for $\chown\subset \dmgm$.

2. Alternatively, one can prove Corollary \ref{newcor} by noting
that $K_0(\chown)=K_0(\chow)[\z(1)]\ob$ and
$K_0(\dmgm)=K_0(\dmge)[\z(1)]\ob$.

\end{rema}
\subsection{Explicit calculation of  $m\ob (\mg^c(X))$;
the weight complex  of smooth quasi-projective
varieties}\label{expmgc}

Let $\mg^c(X)$ for $X\in\smc$ denote the motif of $X$ with compact
support (cf. \S2.2 or \S4.1 of \cite{1}).

\begin{pr}\label{mgc}
For a smooth quasi-projective $X/k$ let $j:X\to P$ be an embedding
for $P\in\spv$, let $P\setminus X=\cup Y_i, 1\le i\le m,$ be a
smooth normal crossing divisor. Let $U_i=\sqcup_{(i_j)}Y_{i_1}\cap
Y_{i_2}\cap\dots \cap Y_{i_r}$ for all $1\le i_1<\dots < i_r\le
m$, $U_0=P$. We have $r$ natural maps $U_r\to U_{r-1}$. We denote by
$d_r$ their alternated sum (as a finite correspondence). We consider
$Q=(Q^i, f_{ij})$, where $Q^i=U_{-i}$ for $0\le i\le -m$, $P^i=0$
for all other $i$; $f_{ij}=d_i$ for $0>i\ge -m$, $j=i+1$, and
$f_{ij}=0$ for all other $(i,j)$. Then $\mg^c(X)\cong  m(Q)$.
\end{pr}
\begin{proof}

Let $j^*$ denote the natural morphism $\mg^c(P)=\mg(P)\to
\mg^c(X)$ (see 4.1 of \cite{1}). Then by Proposition 4.1.5 loc.
cit. the cone of $j^*$  is naturally isomorphic to
$\mg^c(R_0)=\mg(R_0)$, where $R_0=\cup Y_i$ (note that $R_0$ is
proper). Hence our assertion is equivalent to the statement that
$C(R)\cong 0$, where $C$ denotes the Suslin complex of $R$,
 $$R=L(U_m)\stackrel{d_{m*}}{\to}L(U_{m-1})\stackrel{d_{m-1*}}{\to}\dots
 \stackrel{d_{2*}}{\to} L(U_1)\to L(R_0).$$
 Now, it remains to note that $R$ is acyclic as a complex of presheaves. The latter fact easily follows from
  the definition
 of $\smc$, since obviously: for any
 $Z\in\sv$ if $T\subset Z\times R_0$ is closed integral, then $T\subset
 Z\times
 Y_i$ for some $i$.

See also Lemma 7.1 of \cite{grf}.
\end{proof}

We also get an explicit presentation of $\mg^c(X)$ as a complex over
$\smc$ (this corresponds to the first description of $\dmge$ in
\cite{1}). The terms of the complex  are (motives of) smooth
projective varieties.

\begin{rema} \label{rmgc}

1. Applying Proposition \ref{mgc} along with the statements of
subsection \ref{hodge} below we get a nice machinery for computing
cohomology with compact support. Moreover, Proposition \ref{mgc}
appears to be connected with the Deligne's definition of (mixed)
Hodge cohomology of smooth quasi-projective $X$. This is no surprise
(cf. subsection \ref{clawe} below); yet a deeper understanding of
this matter could improve our understanding of cohomology.

2. Using Proposition \ref{mgc} along with Theorem \ref{main} one can
write an explicit formula for $\dme(\mg^c(X),\mg^c(Y))$ for smooth
quasi-projective $X,Y/k$.

Using section 4.3 of \cite{1} one can also calculate
$\dme(\mg(X),\mg(Y)[s])$ for any $X,Y\in \sv$, $s\in \z$. Indeed,
if $\dim X=m$, $\dim Y=n$, $X,Y$ are smooth equidimensional, then
(in the category of geometric motives $\dmgm$) $$\begin{gathered}
\dme(\mg(X),\mg(Y)[s])=\dmgm(\mg(Y)^*, \mg(X)^*[s])\\=\dmgm
(\mg^c(Y)(-n)[-2n],\mg^c(X)(-m)[-2m+s])\\ =\dmge
(\mg^c(Y)(m)[2m],\mg^c(X)(n)[2n+s]);\end{gathered} $$ one can also
use $\hkgm$ for these calculations.
\end{rema}

 In the \S2 of \cite{gs} for any $X/k$ a certain {\it weight complex}
  $W(X)$ of Chow motives was defined.
In order to make the notation of \cite{gs} compatible with ours we
reverse the arrows in the category of Chow motives. Thus we consider
homological Chow motives instead of cohomological ones considered in
\cite{gs}. We have $W(X)\in K^-(\chow)$.

Let $m\ob$ denote the equivalence of $\dms\subset \dme$ with $\hk$
inverse to $m$.

\begin{coro}
For any smooth quasi-projective $X/k$ we have
$t_0(m\ob(\mg^c(X)))\in K^b(\cho)\cong W(X)$ (in $K^-(\chow)$).  
\end{coro}
\begin{proof}
By Proposition 2.8 of \cite{gs} the weight complex of $X$ is
isomorphic to the image in $K^b(\chow)\subset K^-(\chow)$ of the
complex $U$ defined in Proposition \ref{mgc} (with arrows reversed).
\end{proof}

\subsection{$cdh$-hypercoverings; the weight complex of Gillet and
Soul\'e for arbitrary varieties}\label{wecompl}

We recall one of the main tools of \cite{1} (cf. Definition 4.1.9);
it allows  computations with motives of  non-smooth varieties.

\begin{defi}
$cdh$-topology is the smallest Grothendieck's topology such that
both Nisnevich coverings and coverings of the form $X'\coprod Z\to
X$ are $cdh$-coverings; here $p:X'\to Z$ is a proper morphism,
$i:Z\to X$ is a closed embedding, and the morphism $p\ob(X-i(Z))\to
X-i(Z)$ is an isomorphism.
\end{defi}

By Lemma 12.26 of \cite{vbook}, proper $cdh$-coverings are exactly
envelopes in the sense of 1.4.1 of \cite{gs}. Therefore, a
hyperenvelope in the sense of  \cite{gs} is exactly the same thing
as a proper $cdh$-hypercovering. We recall that a
$cdh$-hypercovering is an augmented simplicial variety $X.$ such
that each $X_i\to (cosk_{i-1}sk_{i-1}(X))_i$ is a $cdh$-covering.

We introduce the category $\sprv$. Its objects are varieties over
$k$, its morphisms are proper morphisms of varieties.

 In  \cite{gs} the weight complex functor $W:\sprv\to
K^b(\chow)$ was defined in the following way. The weight complex
 for a simplicial smooth projective variety $T$ was defined
(up to the reversion of arrows) as $T_0\to T_1\to T_2\to\dots$; the
boundary maps were given by alternated sums of face maps. Recall
that we reverse arrows in $W(T)$! 

For $X\in\var$ a proper $Y\supset X$ was chosen; $Z=Y-X$. It was
shown in \cite{gs} that there exist hyperenvelopes $Z.$ of $Z$, $Y.$
of $Y$, and a simplicial closed embedding   $Z.\to Y.$ extending the
map $Z\to Y$, whereas the terms of $Z.$ and $Y.$ are smooth
projective varieties. Then $W(X)$ was defined as the cone of
$W(Z.)\to W(Y.)$ (if we reverse arrows). By means of comparing
different hyperenvelopes
 Gillet and
Soul\'e showed that $W(X)$ is well-defined as an object of $K^b(\cho)$
and gives a functor $\sprv\to K^b(\cho)$.

\begin{pr} \label{cws}
The functor $t_0(m\ob(\mg^c(X))):\sprv\to K^b(\cho)$ is equivalent
(after we reverse all arrows) 
to the functor $W$.
\end{pr}
\begin{proof}
We recall that to compare $t_0\circ m\ob\circ \mg^c$ with the
functor of Gillet and Soul\'e we should fix some choice of $m\ob \circ
\mg^c$ ({\it a priori} the latter one is only defined up to an isomorphism,
cf. Remark \ref{risom}). So first we should check that
$t_0(m\ob(\mg^c(X)))$ is isomorphic to the weight complex of $X$
defined in \cite{gs}; this doesn't require any choices.

Since $Y.\to Y$ is a $cdh$-hypercovering, the $cdh$-sheafification of
the  corresponding complex $L(Y.)\to L(Y)$ is quasi-isomorphic to
$0$. Then Theorem 5.5 of \cite{4} shows that $C(Y.)\cong C(Y)$.
Hence $C(Y.)$ calculates $\mg(Y)$. The same is true for $Z$.

 By
Proposition 4.1.5 of \cite{1} there exists a distinguished triangle
$$\mg(Z)(=\mg^c(Z))\to \mg(Y)(=\mg^c(Y))\stackrel{j}{\to}
\mg^c(X)\to \mg(Z)[1]$$ in $\dmge$. Hence we obtain that
$t_0(m\ob(\mg^c(X)))\cong W(X)$ (in $K^-(\chow)$).

The definitions  imply that  $W$ and $t_0(m\ob(\mg^c(X))$ coincide
as functors on the category of proper smooth varieties (cf. Remark
\ref{risom}). In order to compare the functors in general we can
define $m\ob \circ \mg^c$ using the method of \S2 \cite{gs} (see the
description above). It can be easily seen that this method allows 
lifting $W(X)$ to a functor $W':\sprv\to \hk$ which (as we have just
proved) can be identified with $m\ob \circ \mg^c$. Certainly, in
order to prove that $W'$ is well-defined one should replace the
usage of the Gersten acyclicity (i.e., of Proposition 2 of 1.4.3 of
\cite{gs}) in the proof of Theorem 2 of \cite{gs} by the usage of
Theorem 5.5 of \cite{4}.

\end{proof}

\begin{rema}\label{imp}

1. In Theorem 5.10 of \cite{gu} also a certain functor $h:Sch_k\to
K^b(\chow)$ was constructed ($Sch_k$ is the category of varieties
over $k$). It can be shown that $h$ is equivalent to  the
restriction of $t:\dmge\to K^b(\chow)$ to motives of varieties (see
subsection \ref{skz} for the definition of $t$).

2.  In \S2 of \cite{gs} it was shown that  any two different
representatives $W_i$ of $W(X)$
  (considered as complexes over $\smc$)
  can be connected by a
 chain of certain homomorphisms $h_i$ of complexes of smooth projective
 varieties. Gillet and Soul\'e proved that $h_i$ induce isomorphisms
 on the level of $K^b(\cho)$. The main technical tools were
 Proposition 2 and Theorem 1 of \S1 of \cite{gs} showing that
 hyperenvelopes give quasi-isomorphisms of complexes of Chow
 motives.

 To any such $W_i$ we can associate an object of $\hk$. Since
$t_0(h_i)$ is an isomorphism, the corresponding map of motives will
be an isomorphism too, see part 6 of Theorem \ref{ttn}.

Hence one can prove that the method of \cite{gs} gives a
well-defined motif 
without using the $cdh$-descent reasoning above.

3. More generally, one can easily define Voevodsky's motives of
Deligne-Mumford stacks (i.e., for stacks coming from quotients of
varieties by finite groups) over $k$. For a finite $G,\ \#G=n$,
acting on a variety $X/k$ one can take $\mg(X/G)_\q=a_{G*}
\mg(X_\q)\in \dmge\q$ and $\mg^c(X/G)_\q=a_{G*} \mg^c(X_\q)\in
\obj\dmge\q$. Here $a_G$ is the idempotent correspondence
$\frac{\sum_{g\in G}g}{n}:X\to X$. As we have noted above there
exists a conservative exact weight complex functor $t_\q: \dmge\q\to
\chow\q$ (with properties similar to those of $t$). Certainly, for
$G=\{e\}$ we will have $t_\q(\mg(X/G))=t(\mg(X))$ and
$t_\q(\mg^c(X/G))=t(\mg^c(X))$. Besides, it is most probable that
$t_\q(\mg^c(X/G))$ would coincide with the weight complex for $X/G$
defined by H. Gillet and C. Soul\'e (see Proposition 14 of
\cite{sg}). Indeed, both of these weight complexes can be
calculated using proper hypercoverings of $X/G$. So it seems that
the method of the proof of Proposition \ref{cws} could be extended
to this case. Still the details of the proof of the isomorphness (as
well as of \cite{sg}) have to be written out.
\end{rema}

 \section{Realizations of motives; weight filtration;
  the spectral sequence of motivic descent}\label{new}

One of the main parts of the theory of motives is the problem of
constructing and studying different {\it realizations}, i.e., exact
functors $\dms\to T$ for $T$ being a triangulated category. Some
authors consider functors from the category of (smooth) varieties to
$T$, yet usually those functors can be factored through $\dms$
(cf. \cite{hu} and \cite{le3}).

  In subsection \ref{enhreal} we recall that any differential graded functor from $J$
  gives a realization of $\hk$ (and $\dms$). This method of constructing
  realizations
  is a vast generalization of the method  described
  in 3.1.1 of \cite{gs}. We call  realizations that can be constructed
  from differential graded functors
  {\it enhanceable} realizations; this class
  seems to contain all 'standard' realizations as well as all
  representable functors for the category of motives
  (cf. subsections \ref{etdr} and part 1 of Remark \ref{rzn}).
  In particular, in \ref{etdr} we verify that the \'etale
  cohomology realization is enhanceable; a reader who believes in
  this fact could skip this subsection.

  For any enhanced realization $D$ in \ref{hodge} we define a family of
  {\it truncated realizations}. One could say that
  truncated realizations correspond to 'forgetting
cohomology outside a given range of weights'.  In particular, for
'standard' realizations and motivic cohomology one obtains an
interesting new family of realizations this way.

In subsection \ref{hodge} we also prove that truncated realizations
of {\it length} $N$ can be factored through $t_N$; they give a
filtration on the natural complex that computes $D$.
  We obtain a spectral sequence $S$ converging to $D(Y)$ for a
motif $Y$.  $S$ could be called the {\it spectral sequence of
motivic descent} (note that the usual cohomological descent spectral
sequences compute cohomology of varieties only). For the cohomology
with compact support of a variety $S$ is very similar to the
spectral sequence considered in 3.1.2 of \cite{gs}; yet the origin
of $S$ is substantially different from those of the mentioned one.
Besides we don't need the sheaves to be torsion as one does for
\'etale cohomology. $E_{n}(S)$ can be expressed in terms of $t_{2n-2}(Y)$ (see (\ref{filt}));
in particular, $E_1$-terms depend only on $t_0(Y)$ and have a nice
description in terms of cohomology of smooth projective varieties. 
    $S$ gives a canonical weight
filtration on a wide class of cohomological functors; for the
'standard' realizations this filtration coincides with the usual one
(with indices shifted).

 We note that (as an easy partial case of our results) we get
  a canonical 'weight' filtration on the motivic cohomology of any
variety and the corresponding 'weight' spectral sequence for it.
 $S$ (and the filtration)
are 'motivically functorial'; they are also functorial with respect
to 'enhanced' transformations of functors (this includes regulator
maps).

In subsection \ref{clawe} we  verify that the our definition of
weights gives classical weights for 'standard' realizations (at
least,  rationally). Moreover,  if $W$ denotes the weight filtration
on $H^i(X)$ then $W_{l+N} H^i(X)/W_{l-1} H^i(X)$ is exactly the
corresponding truncated realization; hence it factors through
$t_N$. A morphism $f$ induces a zero morphism on cohomology if
$t_0(f)$ is zero. We also prove that
 the rational length of $X\in\dmge$
coincides with the 'range' of difference of $l$ with the weights of
$H^l(X)$ for all $l$ (cf. Proposition \ref{procon}). If we assume
certain  'standard' conjectures then there would be an equality, see
Proposition \ref{procon}.

We conclude the section by the discussion of $qfh$-descent
cohomology theories and $qfh$-motives of (possibly) singular
varieties. It turns out that a wide class of realizations (including
'standard' ones) are '$qfh$-representable' (hence they are
enhanceable realizations of $\dms$). Moreover, the $qfh$-motif of a
(not necessarily smooth) variety gives 'right values of standard
realizations'.

\subsection{Realizations coming from differential graded functors
('enhanceable' realizations)}\label{enhreal}

We consider the problem  of constructing and studying different {\it
realizations}  of motives, i.e., exact functors $\dms\to T$ for $T$
being a triangulated category.
 Our
description of $\dms$ gives us a simple recipe for constructing
realizations. Any differential graded functor $F:J\to X$ for a
differential graded category $X$ gives an exact functor
$\trp(F):\hk\to \trp(X)$ (and hence also a functor $\hk\to \tr(X)$),
cf. Remark \ref{remf}. It can be easily seen that $\trp(F)$ can be
factored through $t_N$ if $t(J^{l}([Y],[Z]))=0$ for any $Y,Z\in
\spv,\ l< -N$. This is always true if $X^l=0$ for $l<- N$. One can
also note that all functors factoring through $t_N$ can be
reduced (in a certain sense)  to functors of such sort.

We will say that $F$ gives an {\it enhancement} of the realization
$\trp(F)$; a realization that possesses an enhancement can be
called {\it enhanceable}. Obviously, any differential graded
transformation of enhancement induces an exact transformation of
realizations. We will mostly consider contravariant functors $F$.

 Note that for $N=0$, $X$ being equal to $S(A)$ for $A$ an
abelian category (see the definition of $S(A)$ in \ref{dg}), our
construction of $\trp(F)$ essentially generalizes to motives the
recipe proposed in 3.1.1 of \cite{gs} for cohomology of varieties
with compact support (also cf. \cite{gu}).

Note lastly that any exact functor $\dms\to T$ can be uniquely
extended to an exact functor from $\dmge$ to the idempotent
completion of $T$.

\subsection{'Representable' contravariant realizations; \'etale 
cohomology}
\label{etdr}

Now we verify that a large class of realizations are enhanceable;
this includes \'etale cohomology.

To this end we describe a  recipe for constructing a rich family of
contravariant differential graded functors from $J$.  Let $A$ be a
Grothendieck topology stronger than Nisnevich topology (for example,
\'etale topology). We consider the category $\ssc_A$ (i.e., the
morphisms are those of $\smc$, coverings are those of $A$); let
$C(\ssca)$ denote the category of (unbounded) complexes over
$\ssca$. We suppose that for any $X\in \smc$ the representable
presheaf $L(X)=\smc(-,X)$ is a sheaf. We denote by $\dmka$ and
$\dmea$ the categories of unbounded complexes over $\ssca$ that are
similar to the corresponding categories of \cite{1} (i.e., derived
category of complexes of sheaves, resp. derived category of
complexes of sheaves with homotopy invariant cohomology). We also
consider the categories $K(\ssca)$ and $B(\ssca)$ that are unbounded
analogues of $K^-(\ssca)$ and $B^-(\ssca)$ (cf. \ref{dg})
respectively.

Now we verify that the \'etale cohomology realization is enhanceable
in the case when $K$ has finite \'etale cohomological dimension. Let
$Y\in C^+(\ssca)$ be a complex of injective sheaves with transfers
bounded from below with homotopy invariant hypercohomology (we need
the hypercohomology condition if $A\neq Nis$). Now we consider
$C(L(X))$ for $X\in\sv$. Since $C(L(X))$ is quasi-isomorphic to
$\cu(L(X))$, for any $i\in \z$ we have
$$\dmka(C(L(X)),Y[i])=\dmka(\cu(L(X)),Y[i]).$$
 Since the
correspondence $(F)\to \cu(F)$ defines a functor $RC_A$ which is
left adjoint to the embedding $\dmea\to \dmka$ (cf. Proposition
3.2.3 of \cite{1}), $Y\in \dmka$, we have
$$\dmka(C(L(X)),Y[i])=\dmka(L(X),Y[i]).$$ Let $Z\in \obj \hk=\obj
J'$ satisfy $m(Z)\cong \mg(X)=\cu(L(X))$ (in $\dmk$ and so also in
$\dmka$, cf. Corollary \ref{motvar}). Since the terms of $Y$ are
injective sheaves, we conclude that $$H^{-i}(Y)(X)=\kmka
(L(X),Y[i])\cong \kmka (h(Z),Y).$$ Moreover, the complex
$B(\ssca)(h'(Z'),Y)$ computes the complex $Y(X)$ up to a
quasi-isomorphism (see the definitions of \ref{expl}).

Now we describe how the formalism of \S2 can be applied to the
computation of $B(\ssca)(h'(Z),Y)$. We have a contravariant functor
$Y^*:J\to C(\ab)$ that maps $[P]\in \obj J$ to $B(\ssca)(C(P),Y)$.
Let $a:J'=\prt J\to B(\ssca)$ denote the differential graded functor
induced by the embedding $J\to B^-(\ssca)$, cf. Remark \ref{remf}
and Proposition \ref{mdg}. Since $Y^*=B(\ssca)(-,Y)\circ a$, we
obtain that
$$B(\ssca)(h'(Z),Y)\cong \prt (Y^*)(Z).$$ Here $\prt (Y^*)$
denotes the extension of $Y^*$ to $J'$, cf. Remark \ref{remf}.

For example, we can take $A$ being the \'etale site. One can take $Y$ being an
injective resolution of $\z/n\z$  (or a resolution of any other
\'etale complex $C$ with transfers with homotopy invariant
hypercohomology) by means of \'etale sheaves with transfers. By
Proposition 3.1.8 and Remark 2 preceding Theorem 3.1.4 in \cite{1}
the cohomology of $Y(L(X))$ for $X\in\sv$ will compute the 'usual'
\'etale hypercohomology  (i.e., in the category of sheaves without
transfers) of $C$ restricted to $X$. Hence $\tr(Y^*)$ gives the
corresponding realization of motives.
 We obtain that in order to
compute the \'etale realization of motives (with coefficients in
$\z/n\z(r)$ for any $n>0,\ r\ge 0$) it suffices to know the
restriction of the corresponding 'representable functor' to the
subcategory of $\ssc_A$ consisting of sheaves of the form $C^i(P)$
for $P\in\spv$. Note also that we can compute morphisms in the
category of presheaves with transfers.

If $C$ is a complex  of $cdh$-sheaves (see subsection \ref{wecompl})
with transfers then the cohomology of $Y(L(X))$ for $X\in\sv$ will
compute the $cdh$-hyper\-cohomology of $C$ restricted to $X$ for any
$X/k$. Yet $cdh$-topology is not subcanonical; $cdh$-hypercohomology
does not necessarily coincides with 'usual' hypercohomology of $C$.
For  the
 computation of 'usual'
cohomology for singular varieties the $qfh$-topology seems to be
more useful; see subsection \ref{sing} below.

One can check  that the Galois action on $H^i_{et}(X\times_{\spe\,
k}\spe\, \overline{k},\z/n\z)$, where $\overline{k}$ is the
algebraic closure of $k$, $n$ is prime to the characteristic of $k$,
can be expressed in terms of our formalism. Indeed, isomorphisms
in the derived category corresponding to the Galois action can be
extended to an injective resolution of $\z/n\z$. Therefore all
statements of this section are valid for \'etale realization with
values in Galois modules.

A similar method for constructing the \'etale realization (without
using the formalism of differential graded categories) was described
in \cite{hu} (see the reasoning following Proposition 2.1.2).
 In the next subsection we describe a general method of obtaining
 weight filtrations for realizations.

It seems that that the same method can be applied to other
'classical' realizations including the  'mixed realization' one.
Yet filling out the details is rather hard. Fortunately, one can
avoid this by applying the {\it weight structure} formalism of
\cite{bws}; see  Remark \ref{ass} below.

\subsection{The spectral sequence of motivic descent;
weight filtration of realizations; the connection with
$t_N$}\label{hodge}

Now we consider a  contravariant functor $F:J\to B(A)$ for an
abelian $A$. We denote the functor $\prt(F):J'\to B(A)$ by $G$,
denote $\tr(F):\hk\to K(A)$ by $E$.

It seems very probable that one can 'enhance' all 'classical'
realizations this way (possibly, for a 'large' $A$). Besides, it
practice it usually suffices to consider functors whose targets are
categories of complexes bounded (at least) from one side.

The constructions of this subsection use the results of subsection
\ref{dop} heavily.

We recall that for a complex $X$ over  $A$, $a,b\in\z$, $a\le b$,
 its canonical $[a,b]$-truncation is the complex $$X^a/dX^{a-1}
\to X^{a+1}\to \dots X^{b-1}\to \ke (X^b\to X^{b+1}), $$ here $\ke
(X^b\to X^{b+1})$ is put in degree $b$; for $a=b$ we take $H^a(X)$.
 We also consider truncations of the type $\tau_{\le b}$
(i.e., truncations from above).

For any $b\ge a\in \z$ we consider the following functors (see
subsection \ref{dop}). By $F_{\tau_{\le b}}$ we denote the functor
that that sends $[P]$ to
 $\tau_{\le b}(F([P]))$. By $F_{\tau_{[a,b]}}$ we denote the functor
 that that sends $[P]$ to
 $\tau_{\le b}(F([P]))/\tau_{\le a-1}(F([P]))$. For $N=a-b$
  we consider the functor $F_{b,N}$
  that sends $[P]$ to  the $[a,b]$-th canonical
truncation of $F([P])$. These functors are differential graded;
hence they extend to $G_b=\prt(F_{\tau_{\le b}}):J'\to B^-(A)$,
$G_N^b=\prt(F_{b,N}):J'\to B^b(A)$, and
$G_{a,b}=\prt(F_{\tau_{[a,b]}}):J'\to B^b(A)$. We recall that
$G_{a,b}$ and $G_N^b$ are connected by a canonical functorial
quasi-isomorphism, see part 2 of Proposition \ref{trunc}. The reason
for considering both of them is that the functors $G_{a,b}$ are more
closely related to the spectral sequence (\ref{spectr}) below,
whereas $G_N^b$ behave better with respect to $t_N$.

 We denote $\tr(F_{b,N}):\hk\to K^b(A)$ by $F_N^b$. Since
$F_{b,N}$ is concentrated in degrees $[b-N,b]$, $F_{b,N}$ maps all
$J^m(X,Y)$ for $X,Y\in\obj J$, $m< -N$, to $0$. Hence one can
present $F^b_N$ as $b_{b,N}:\hk_N\to K(A)$ for a unique
 $b_{b,N}:\hk_N\to K(A)$.
 The set of
 $F^b_N$ could be called  {\it truncated realizations} for the realization $E$;
 $N$ is the {\it length} of the realization.
 These realizations appear to be new even in the case when $E$ is
 the \'etale cohomology.
Note that for $X=[P],\ P\in\sv$,
 the truncated realizations give exactly the corresponding
 truncations of $E([P])$ (i.e., of the corresponding 'cohomology'
 of $P$); that is what one usually expects from the weight filtration.

 The complexes $G_b(X)$ give a
 filtration of $G(X)$ for any $X\in\obj J'$; moreover $G_{a,b}(X)=G_b(X)/G_{a-1}(X)$.

Let $X=(P^{i}, q_{ij})\in \obj J'=\obj\hk$.
 We obtain the spectral sequence of a filtered complex (see \S III.7.5 of\cite{gelman})
  \begin{equation} S:E_1^{ij}(S)\implies H^{i+j}(G(X))\label{spectr}\end{equation}
 we call it the {\it spectral sequence of motivic descent}.
 Here $E_1^{ij}(S)=H^{i+j}G^{j,j}(X)=H^{i+j}(F^{j}_0(X))$.
Note that  $H^{i+j}(G(X))=
 H^{i+j}(E(X))$, in the right hand side we consider $X$ as an
 object of $\hk$.
It is easily seen that $S$ is $\hk'$-contravariantly functorial in
$X$; in fact, they are also  functorial (for example, see
(\ref{filt})) below. In fact, one could also define $E_0(S)$; this
would be $\hk'$-functorial also, yet difficult to compute.

Moreover, if $h:F\to F'$ is a differential graded transformation of
functors then the corresponding map of spectral sequences depends
only on $\tr(h)$. In particular, for the \'etale 
realization the spectral sequence does not depend on the choice of
an injective resolution for the corresponding complex (see the
previous subsection).

The spectral sequence $S$ is similar to those   coming from
hypercoverings (and hyperenvelopes). Yet its terms are 'much more
functorial'; it computes cohomology of any motif (not necessarily of
a motif of a variety).

By definition, $E_1^{ij}(S)=H^{i+j}(F^{j}_0(Y))$ is the $i$-th
cohomology group of the chain complex $A_l=H^j(P^{-l})$.
  Hence the $E_1$-terms are functorial in the complex
  $(P^l)\in K^b(\cho)$, i.e., in $t_0(Y)$. $S$
 is convergent:  if for $X$ of (stupid) length $N$ we choose a
 representative in $\hk'$ of length $N$
 then only $N+1$ rows of $E_1(S)$ would be non-zero.
 Besides if all $F^j_0(Y)$ are acyclic then $E(Y)$
 is acyclic.
 We denote the  filtration on $H^{s}(E(Y))$ given by $S$ by
 $W_l$; we call it the weight filtration of $H^s$.

 For any $b,N$ we also have a 'spectral subsequence'
  $$S_N^b:
 E_1^{ij}(S_N^b) \implies H^{i+j}(G_{a,b}(X))= H^{i+j}(G_N^b(X))=H^{i+j}(F^b_N(Y)).$$
 Its $E_1$-terms form a subset of the $E_1$-terms of $S$, the (non-zero)
 boundary maps are the same. We also have weight filtration on
$H^{s}(F^b_N(Y))$.

For any $0\le l\le N$ we have an obvious spectral sequence morphism
$S_{N-l}^{b-l}\to S_N^b$. It induces an epimorphism
$$\alpha^s_{l,b,N}: H^s(F^{b-l}_{N-l}(Y))\to W_{b-l}
(H^{s}(F^b_N(Y))).$$

Now we use an argument that can be applied to any filtered
complex. For any $N\ge 0$, $n\ge 1$, one easily sees that 
$E_n^{ij}(S_N^b)=E_{\infty}^{ij}(S_N^b)$ if $b-n<j<b-N+n$. Moreover,
if $b+1-n\ge j\ge b-N+n-1$ then $E_n^{ij}(S_N^b)=E_n^{ij}(S)$.
Therefore we have
\begin{equation}\begin{gathered}
E_n^{ij}(S)=Gr_{n-1}^WH^{i+j}(F^{j+n-1}_{2n-2}(Y)) =
W_{n-1}(H^{i+j}(F^{j+n-1}_{2n-2}(Y)))/W_{n-2}(H^{i+j}(F^{j+n-1}_{2n-2}(Y)))\\
=\imm \alpha^{i+j}_{n-1,j+n-1,2n-2}/\imm \alpha^{i+j}_{n,j+n-1,2n-2},
\end{gathered}\label{filt}\end{equation}
i.e., it is the middle factor of the weight filtration of
$H^{i+j}(F^{j+n-1}_{2n-2}(Y))$.
 A similar equality can be written for
$E_n^{ij}(S_N^b)$ for any $b\in\z, N\ge 0, n\ge 1$. Hence for any
$n\ge 1$ the $E_n$-terms of $S$ and all $S_N^b$ depend only on
$t_{2n-2}(Y)$.

Suppose that $X\in \hk'_{[c,d]}$. It can be easily verified (for
example, using the spectral sequence (\ref{spectr})) that for any
$j\in\z,b-N-c\le j \le b-d$, the $j$-th cohomology group of
$F_N^b(Y)$ coincides with $H^j(E(Y))$. Besides for any $j\in\z$ the
weights of $H^j(E(Y))$ lie between $j+c$ and $j+d$. In particular,
by part 1 Theorem \ref{ttn} the weights of $H^j(X)$ for $X\in\sv$ of
dimension $N$ lie between $j$ and $j+N$, the weights of $H^j_c(X)$
(the cohomology with compact support) lie between $j-N$ and $j$.

Hence all 'cohomological information' of a motif of length $\le N$
can be factored through $t_N$ (in order to prove this for fine
length $\le N$ one should use Proposition \ref{idempdm}). This
statement can be considered as the 'realization version' of Theorem
\ref{ttn} (parts 1 and 2).

\begin{rema}\label{rzn}

1. For an object $U\in \obj J'=\obj \hk'=\obj\hk$ (recall that
$J'=\prt(J)$) one can consider the differential (contravariant)
graded functor $J_U:J\to B(\ab)$ that maps $[P]$ to $J'([P],U)$.
Then $H^i(\tr(U))(X)=\hk(X,U[i])$; this means that representable
realizations are enhanceable.

In particular, we can take $U=\z(n)$ for $n\ge 0$. Hence we obtain
canonical 'weight' filtration on the motivic cohomology of any
variety and the corresponding 'weight' spectral sequence for it. A
simple example of this spectral sequence is closely related with the
Bloch's long exact localization sequence for higher Chow groups (see
\cite{blo}).

Indeed, let $Z\subset X \in \spv$, $Z$ is everywhere of codimension
$c$, let $Y=X-Z$. Then (cf. the proof of part 1 of Theorem
\ref{ttn}) the motif of $Y$ is a cone of $\mg(X)\to \mg(Z)(c)[2c]$.
So the length of $\mg(Y)$ is $1$; hence for any realization the
formula for $E_1(S)$ relates the cohomology of $Y$ with those of $X$
and of $Z(c)[2c]$. If the realization is motivic cohomology then
this relation is exactly the same as given by  the exact sequence of
Bloch. See also the end of Remark 2.4.2 of \cite{bws} and Remark
\ref{ass} below.

 This example shows that the weight filtration obtained this way is non-trivial in general;
 it appears  not to be mentioned in the literature.
The filtration is compatible with the regulator maps (whose targets
are 'classical' cohomology theories).

 Using the
spectral sequence relating algebraic $K$-theory to the motivic
cohomology (see \cite{frisus} and \cite{grf}), one can also obtain a
new filtration on the $K$-theory of a smooth variety $X$.

2. The spectral sequence functor  described is additive (here we fix
$F$ and consider $S$ as a functor from $\hk$). Hence one can easily
see that this functor can be uniquely extended to the whole $\dmge$.

3. It seems to be interesting to study the truncations and the
weight spectral sequence for the cohomology of the sheaf $Y\to
G_m(X\times Y)$ for a fixed variety $Y$ ($G_m$ is the multiplicative
group). These things seem to be related with the Deligne's
one-motives of varieties as they were described in \cite{bv}.

More generally, for any $X\in \dmge$ the functor
$\ihom(-,X):\dmge\to \dme$ is enhanceable.

4.
 In Remark 6.4.1 of \cite{bws} we verify that the truncated realizations coming
from representable realizations are representable also. The
truncated realizations are represented by $t$-truncations of the
objects representing the original realizations with respect to a
certain {\it Chow} $t$-structure. Note that this is the case for
motivic cohomology and for 'classical' realizations of motives (with
values in $\ab$).

More generally, one can define a $t$-structure on the category
$\tr(DG-Fun(J,B(A)))$ (differential graded functors) that
corresponds to the canonical truncation of $A$-complexes. Then the
realizations of the type considered here correspond to some objects
of this category; truncations of a realization with respect to this
$t$-structure would be exactly its truncated realizations.

5. Another important source of differential graded functors from $J$
(generalizing representable functors considered in part 1) are those
coming from localizations of $\hk$ (or of $\dms$ which is the same
thing). It will be discussed in subsection \ref{loc} below.
\end{rema}

\subsection{Comparison with 'classical weights'; comparison of the rational length
with the 'Hodge length'}\label{clawe}

Suppose now that there are no maps between different weights, i.e.,
for any $P,P',Q,Q'\in\spv,f\in \smc(P,P'),\ g\in \smc (Q,Q'),\ i\neq
j,$ we have $$A(\ke(H^i(E([P']))\stackrel{f^*}{\to}
H^i(E([P]))),(\cok( H^j(E([Q'])) \stackrel{g^*}{\to}
H^j(E([Q])))=0.$$ Recall that this condition is fulfilled for the
\'etale and Hodge realizations with rational coefficients. Then $S$
and all $S_N^b$ degenerate at $E_1$. Therefore $H^jF^b_N(X)=
W_{b}(H^j(E(Y)))/ W_{b-N-1}(H^j(E(Y)))$.

Besides, if there are no maps between different weights (in some
variation of the target category), then for any $H^l(X)$ there
cannot exist more than one filtration $W_j$ on $H^l(X)$ such that
$W_j/W_{j-1}$ is of weight $j$.  Hence for the \'etale and Hodge
realizations our weight filtration coincides with the usual one
(with indices shifted). Therefore, $t_N$ can be called the weight
functors.

The existence of the weight spectral sequence easily implies that
the length of the weight filtration of $H^l$ is not larger than the
stupid length of $M\in \obj \dms$. Moreover, Proposition
\ref{idempdm} gives the same inequality for the fine length of any
object of $\dmge$ and also for the rational length (on the rational
level). Besides, if $M\in \dmge{}_{[a,b],\q}$ then the weights of
$H^l\otimes \q(M)$ lie between $l+a$ and $l+b$.

\begin{conj}\label{lecon}
The converse implication is true also, i.e., if for all $l$ the
weights of $H^l\otimes \q(M)$ (here $H^*$ is the singular
cohomology) lie between $l+a$ and $l+b$ then $M\in
\dmge{}_{[a,b],\q}$.
\end{conj}

Now we show that Conjecture \ref{lecon} follows from certain
'standard' conjectures.

\begin{pr}\label{procon}

Suppose that following statements are valid:

1. The Hodge 
conjecture.

2. Any morphism of  Chow motives that induces an isomorphism on
singular cohomology is an isomorphism.

Then Conjecture \ref{lecon} is also valid.
\end{pr}

\begin{proof}  Indeed, we should check that if
$M\not \in \dmge{}_{[a,b],\q}$ then at least for one $l$ the weights
of $H^l\otimes \q(M)$ do not lie between $l+a$ and $l+b$. Using
Proposition \ref{idempdm} we obtain that we have the motivic descent
spectral sequence for the singular realization of objects of
$\dmge$; its $E_1^{ij}$-term is the $i$-th cohomology group of the
chain complex $A_l=H^j(P^{-l})$. Therefore one should check

(1) If the 'first term' of $t_\q(M)=(P^i)$ is at $c$-th place and
the map  $g:P^c\to P^{c+1}$ is not a projection onto a direct
summand then for some $l$ the $(l+c)$-th weight component of
$H^l\otimes \q(M)$ is non-zero.

(2) The dual statement.

We verify (1); (2) is similar (and follows from (1) by duality).
Suppose that for any $l$ the $l+c$-th weight piece of $H^l\otimes
\q(M)$ is zero. Note that this piece equals  $\cok (g_l^*: H^l
(P^{c+1})\to H^l(P^c))$; hence all $g_l^*$ are surjective. Since the
category of (rational) polarizable pure Hodge structures is semisimple, we can
choose a family of splittings for $g_l^*$ (that respect the Hodge
structures). Since the splittings are Hodge, by Hodge 
conjecture (for $P^c\times P^{c+1}$) we obtain that these splittings
can be realized by a morphism  $h:P^{c+1}\to P^c$ of Chow motives.
We obtain that $H^*(h)$ gives a splitting of $H^*(g)$; the assumption (2) (applied to $h\circ g$) implies that $g$ is  a projection onto a
direct summand.
\end{proof}

\begin{rema}

1. Note that Conjecture \ref{lecon} implies the conservativity of
the singular realization, which certainly implies condition 2 of
Proposition \ref{procon}.

2. In order to get a similar result for the \'etale realization one
should replace the Hodge conjecture by the Tate conjecture.

3. One could also count the (minimal) number of non-zero terms of
$t(M)$.

\end{rema}

We also note that for any morphism $f:Y\to Z$ for $Y,Z\in \hk$ the
morphisms $H^l(f):H^{l}(E(Z))\to H^{l}(E(Y))$ for $Y,Z\in \hk$ are
strictly compatible with the weight filtration. Therefore $H^l(f)$
is zero if and only if the corresponding map of $E_1$-terms in
$S(Y)\to S(Z)$ is zero. Hence the map $\hk(Y,Z)\to A(H^l(E(X)),
H^l(E(Y)))$ factors through $t_{0*}$.

For cohomology with integral coefficients one may apply the previous
statement for rational cohomology to obtain that $f^*$ is zero on
$H^l\otimes \q$ if $t_0(f)=0$. Hence if $t_0(f)=0$ then $f^*$ is
zero on cohomology modulo torsion.

Using the results of subsection \ref{wecompl} one can compute
$\hk(Z,T)$ for $Z=m\ob(\mg(X))$, $T=m\ob(\mg(Y))$, and also compute
$$M_N(X,Y)=\imm t_{N*}(\hk(Z,T)\to \hk_N(t_N(Z),t_N(T)).$$

\begin{rema}\label{ass} 
Formally all our filtrations and spectral sequences depend on the
choices of enhancements for realizations. One can check that they
are independent in fact; yet in order to use the theory described
above it is necessary (at least) to prove that enhancements exist.
This seems to be true for all reasonable cases; yet proofs could be
difficult (see \ref{sing} below).

There is a way to avoid these difficulties completely; it is
studied in \cite{bws}. The idea is to consider a set of axioms of
so-called {\it weight structures} for a triangulated category $C$;
the axioms are somewhat similar to those for $t$-structures (yet
the consequences of the axioms are quite distinct!). One could say
that any object of $C$ has a {\it weight decomposition} which is
not unique but is 'unique up to a homotopy' (in a certain sense).
In the case when $C=\hk$ one could consider the decompositions
given by the 'stupid filtration' (see Proposition \ref{stufil}).
Then for any 
(covariant) functor $H:C\to A$, where $A$ is an abelian category,
one can consider the filtration of $H(X)$ by $H(X_{[a,+\infty]})$
for $a\in\z$. Note that while $X$ does not determine
$X_{[a,+\infty]}$ uniquely, it does determine the image of
$H(X_{[a,+\infty]})$ in $H(X)$. If $H$ is contravariant then one
should consider the image $H(X_{[-\infty,b]})\to H(X)$. The objects
$X_{[a,+\infty]}$,  $X_{[-\infty,b]}$ give a Postnikov tower for
$X$; this yields a spectral sequence $T$ converging to $H(X)$. Its
$E_1$-terms are $H_i(P^j)$; hence they are not determined by $X$.
Yet starting from $E_2$ all terms of the spectral sequence are
canonical and functorial in $X$. In fact, the filtration induced by
$G_b$ on $G$ can be obtained from the filtration for $T$ by
Deligne's decalage.  Hence the spectral sequences and the filtration
for $T$ coincide (up a certain shift of indices) with the terms
of $S$ defined in (\ref{spectr}).

This approach for constructing weight filtrations of realizations of
Hanamura's motives was described in \cite{ha3}; yet the proof of
Proposition 3.5 loc. cit. relies heavily on  (sort of) enhancements
for realizations.

The advantage of our alternative approach is that enhancements are
no longer needed; in particular, it can be applied to the stable
homotopy category for which no (differential graded) enhancements
exists. Yet in this abstract setting it is difficult to define
truncation functors (especially the 'higher' ones). The reason for
that is (as was noted by several authors) that the axiomatics of an
(abstract) triangulated category is not 'rigid enough'.

\end{rema}

\subsection{$qfh$-descent cohomology theories;
 motives of singular varieties}\label{sing}

 Some 'standard' cohomology theories are difficult to represent by a complex of
sheaves with transfers (in the way described in subsection
\ref{etdr}). One of the ways to do this is to use $qfh$-topology.

We recall that the $qfh$-topology is the topology on the set of all
varieties whose coverings are quasi-finite universal topological
coverings (see \cite{htop} for a precise definition). In particular,
the $qfh$-topology is stronger than the flat topology and the
$cdh$-toplogy. There is a natural functor from $\dme$ to the derived
category of $qfh$-sheaves with homotopy invariant cohomology (it is
denoted as $DM_{qfh}(k)$); this functor is surjective on objects.
Note also that any 'ordinary' topological sheaf restricted to $\var$
gives a $qfh$-sheaf. Moreover, $qfh$-descent follows from proper
descent combined with Zariski descent (see section 2 of
\cite{htop}).

Let $C$ be a complex of presheaves (possibly without transfers)
whose hypercohomology satisfies $qfh$-descent. Then the
$qfh$-hypercohomology of the $qfh$-sheafification of $C$ coincides
with the hypercohomology of $C$ (for example, a similar statement
was proved in the proof of Theorem 5.5 of \cite{4}). Therefore if
the cohomology of $C$ is homotopy invariant then it can be
presented by means of 'representable' functors on $\dms$ as it was
described in subsection \ref{etdr} above. In particular, this shows
that Betti and Hodge cohomology theories can be enhanced to
differential graded realizations.

Now let $C$ be a complex  of  $qfh$-sheaves. It was proved in
\cite{htop} (see Theorems 3.4.1 and 3.4.4) that the
$qfh$-hypercohomology of a variety $X$ with coefficients in $C$
coincides with the \'etale hypercohomology of $C$ in the cases when
either $C$ is a
 $\q$-vector space sheaf complex and $X$ is normal or $C$ is a
 locally constant \'etale sheaf complex. Hence in this cases
 the \'etale hypercohomology of $C$ also gives a 'representable' realization.

  Note that these realization compute
$qfh$-hypercohomology with coefficients in $C$ of any (not
necessarily smooth) variety $X$. Hence $\mg(X)_{qfh}$ (i.e., the
image of $\mg(X)$ in $DM_{qfh}(k)$) seems to be the natural choice
for the $qfh$-motif of a (possibly) singular variety. In particular,
its 'standard' realizations have  the 'right' values of
$H^i_{et}(X,\z/l\z(m))$ (at least, in the case when $k$ has finite
\'etale cohomological dimension).

\section{Concluding remarks}\label{last}

In subsection \ref{ssubcat} we  give a general description of
subcategories of $\hk$  generated by  fixed sets of objects. In
particular, this method can be used to obtain the description of the
category of effective Tate motives (i.e., the full triangulated
subcategory of $\dms$ generated by $\z(n)$ for $n>0$).

In \ref{loc} we describe the construction of 'localization of
differential graded categories' (due to Drinfeld). This gives us a
description of localizations of $\hk$. All such localizations come
from differential graded functors. As an application, we prove that
the motif of a smooth  variety is a mixed Tate one if and only if its
weight complex (as defined in \cite{gs}, cf. subsection \ref{expmgc}
and \ref{wecompl}) is.

In \ref{motcharp} we verify that over an arbitrary perfect field
one can apply our theory (at least) with rational coefficients.
One of the main tools is the Poincare duality in characteristic
$p$ proved by Beilinson and Vologodsky. Moreover over finite
fields the Beilinson-Parshin conjecture (that $H^i(P,\q(n))$ for
smooth projective $P$ could be  non-zero only for $i=2n$) holds
iff $t \q:\dmge\q\to K^b(\chow \q)$ is an equivalence (note that
here $\dmge\q$, $\chow\q$, and $t \q$ denote the corresponding
idempotent completions); see Proposition \ref{becon}. We also
describe an idea for constructing certain 'infinite integral'
weight complex in finite characteristic (see \S7.3 of \cite{bws}
for a complete proof).

In \ref{sindepl} we prove that traces of endomorphisms of cohomology
of motives induced by endomorphisms of motives do not depend on the
choice of a Weil cohomology theory. This result generalizes Theorem
3.3 of \cite{bloesn}.

In \ref{summand} we remark that one can easily add direct summands
of objects to $J$. In particular, one could include $[P][2i](i)$
into $J$.

In \ref{mn} we consider a functor $m_N:\hk_N\to \dme$ that maps
$[P]$ into the $N$-th canonical truncation of $C(P)$ (as a complex
of sheaves).

\subsection{Subcategories of $\hk$
generated by  fixed sets of objects}\label{ssubcat}

Let $B$ be a set of objects of $\hk$; we assume that $B$ is closed
with respect to direct sums.

Let $B'$ denote some full additive subcategory of $J'=\prt(J)$ such
that the corresponding objects of $\hk$ are exactly elements of $B$
(up to isomorphism).

Let $\mathfrak{B}$ denote the smallest triangulated category of
$\hk$ containing $B$.

\begin{pr}\label{subcat} $\mathfrak{B}$ is
canonically isomorphic to $\trp(B')$.
\end{pr}
\begin{proof} Follows immediately from Theorem 1 in \S4 of \cite{bk}.
\end{proof}

\begin{rema}\label{rscat}
1. It follows immediately  that in order to calculate the smallest
triangulated category of containing an {\it arbitrary} fixed set of
objects in $\hk$ it is sufficient to know morphisms between the
corresponding objects in $J'$ (i.e., certain complexes) as well as
the composition rule for those morphisms.

2. We obtain that for any triangulated subcategory of $D\subset\hk$
the embedding $D\to\hk$ is isomorphic to $\trp(E)$ for some
differential graded functor $E:F\to G$. Here $G$ is usually equal to
$J'$ (though sometimes it suffices to take $G=J$); $F$ depends on
$D$.

It follows that for any $h\in \obj \hk$ the representable
contravariant functor $h^*:D\to \ab: d\to \hk(d,h)$
 can be represented as $H^0(u)$ for some contravariant differential
 graded functor $u:F\to B^-(\ab)$. See part 2 Proposition \ref{local} below
 for a similar statement for localization functors.

3. Using this statement one can easily calculate the triangulated
category of (mixed effective) Tate motives (cf. \cite{le1}). It is
sufficient to take $B=\sum_{a_i\ge 0}[({\mathbb P}^1)^{a_i}]$, i.e.,
the additive category generated by motives of non-negative powers of
the projective line. This gives a certain extension of the
description of \cite{spit} to the case of integral coefficients.
Note that the description of the category of effective Tate motives
immediately gives a description of the whole category of Tate
motives since $\z(1)$ is quasi-invertible with respect to $\otimes$.
Besides, one could expand $J$: see \S\ref{summand} below.

\end{rema}

\subsection{Localisations of $\hk$} \label{loc}

Let $C$ be a differential graded category satisfying the homotopical
flatness condition, i.e., for any $X,Y\in\obj C$ all $C^i(X,Y)$ are
torsion-free. Note that both $J$ and $J'=\prtp(J)$ satisfy this
condition.

In \cite{drinf} V. Drinfeld has proved (modifying a preceding result
of B. Keller) that for $C$ satisfying the homotopical flatness
condition and any full differential graded subcategory $B$ of $C$
there exists a differential graded quotient $C/B$ of $C$ modulo $B$.
This means that there exists a  differential graded $g:C\to C/B$
that is surjective on objects such that $\trp(g)$ induces an
equivalence $\trp(C)/\trp(B)\to \trp(C/B)$ (i.e.,
$\trp(C)/\trp(B)\cong \trp(C/B)$, $\trp(g)$ is zero on $\trp(B)$ and
induces this equivalence).

The objects of $C/B$ are the same as for $C$, whereas for $C_1,C_2\in
\obj C=\obj (C/B),\ i\in\z$, we define
\begin{equation}\label{drf}
\begin{gathered}(C/B)^i(C_1,C_2)=C^i(C_1,C_2)\bigoplus \bigoplus_{j\ge
0}\\ \bigoplus_{B_1,\dots,B_j\in \obj B,\sum
a_i=i+j}C^{a_1}(C_1,B_1)\otimes \eps_{B_1}\otimes
C^{a_2}(B_1,B_2)\otimes \eps_{B_2}\dots  \otimes \eps_{B_j}\otimes
C^{a_j}(B_j,C_2).\end{gathered}\end{equation} Here $\eps_b\in
(C/B)^{-1}(b,b)$ for each $b\in\obj B\in \obj (C/B)$ is a 'canonical
new morphism' such that $d_b\eps_b=\id_b$; $\eps_b$ spans a
canonical direct summand $\z \eps_b\subset (C/B)^{-1}(b,b)$. From
this condition one recovers the differential on morphisms of $C/B$.

For example, this construction  (for $C=J$) gives an explicit
description of the localization of $\hk$ by the triangulated
category generated by all $[Q],Q\in\spv, \dim Q<n$, for a fixed $n$
(and hence also of the corresponding localization of $\dms$).

More generally, for localizations of $\hk$ modulo some $A\subset
\hk$ it is sufficient to know the complexes $C/B([P],[Q])$ and the
composition law for a certain $B$ and all $P,Q\in\spv$; here either
$C=J$ or $C=J'$. In the case when $A$ is not generated by objects of
(stupid) length $0$ we are forced to take $C=J'$; this makes the
direct sum in (\ref{drf}) huge.

\begin{rema}
1. One can extend $J$ as in \S\ref{summand} below. This allows one to
obtain a certain description of the (triangulated) category of {\it
birational} motives; see \cite{kabir}.

2. Explicit calculations by means of  (\ref{drf}) could be
difficult.
 We will only need (for part 3 of
Proposition \ref{local} below) the following obvious property of the
construction:
 if $C^i(-,-)=0$ for $i>0$ then the same is true for
$C/B$.  Note also that in this case $C^0(X,Y)= (C/B)^0(X,Y)$ for all
$X,Y\in\obj C$; $H(C/B)(X,Y)$ is a certain easily described factor
of  $HC(X,Y)$. These facts immediately yield the formula for
$\homm(\stackrel{\_}{M}(X), \stackrel{\_}{M}(Y)[i])$ for $i\ge 0$ as
described in Corollary 7.9 of \cite{kabir} (note that loc. cit. is
false for $i<0$); here $X,Y\in \spv$, $\stackrel{\_}{M}(X),
\stackrel{\_}{M}(Y)$ are their birational motives.

\end{rema}

\begin{pr}\label{local}
1. If $F:\hk\to T$ is a certain localization functor ($T$ is a
triangulated category) then $F\cong\trp(G)$ for a certain
differential graded functor $G$ from $J$.

 2. For any $t\in T$ the contravariant functor $t^*:\hk\to \ab:X\to T(F(x),t)$
 can be represented as $H^0(u)$ for some contravariant differential
 graded functor $u:J\to B^-(\ab)$.

3. Let $B$ be a full additive subcategory of $J$, $\obj B= T\subset
\spv$. Let $\mathfrak{B}$ denote the smallest triangulated
subcategory of $\hk$ that contains all objects of $B$ and is closed
with respect to taking direct summands (in $B$). Then for $M\in \obj
\hk$ we have $M\in \obj \mathfrak{B}$ if and only if $t_0(M)$ is  a direct
summand of a  complex all whose terms have the form $[Q],\ Q\in T$.

\end{pr}
\begin{proof}

 1. Let $A=\{X\in \obj \hk, F(X)=0\}\subset \obj J'=\obj \hk$, 
 denote the corresponding full subcategory of $J'$ by $B'$. By
Proposition \ref{subcat}, $\trp(B')$ is isomorphic to the
categorical kernel of $F$. Let $H:J'\to J'/B'$ denote the functor
given by Drinfeld's construction. Since $\trp(J')=\hk$, we obtain
$\trp(H)\cong F$. Hence by part 2 of Remark \ref{remf} we can take
$G$ being the restriction of $H$ to $J\subset J'$.

2. Let $w$ denote some element that corresponds to $t$ in
$\prtp(J'/B')$. Then it suffices to take
$u([P])=\prtp(J'/B')(\prtp(G)([P],w))$.

3. If $M\in \obj \mathfrak{B}$ then $t_0(M)$ belongs to the
triangulated subcategory of $\hk_0=K^b(\cho)$ that contains all
$[Q],\ Q\in T$ and is closed with respect to taking direct
summands. This subcategory consists exactly of complexes described
in the assertion 3.

We prove the converse implication.

We can assume that $t_0(M)$ is homotopy equivalent to a complex all
whose terms have the form $[Q],\ Q\in T$.

 Let $I$ denote $J/B$, let $K=\trp(I)$.
 Let $S:J\to J/B$ be the localization functor of \cite{drinf}; let $u=\trp(S)$.
  Since $I^i(-,-)=0$ for $i>0$, we have natural functors $I\to S(HI)$ (see \S \ref{dg}) and
  $v:K\to K_0$, where $K_0=K^b(HI)$ (the weight complex functor for this case).
Hence we have
 functors $u_0:\hk_0\to K_0$ and $v:K\to K_0$ such that
 $u_0\circ t_0= v\circ u$. We obtain $v(u(M))=0$.

Suppose that $M=j(X)$ for some $X\in\obj\hk'$. We can also assume
that $X\in \obj K$.


For any object $U\in \obj K=\obj \prtp(I)$  we consider the
 (contravariant) differential graded functor $U^*:J\to
B(\ab)$ that maps $[P]$ to $\prtp(I)([P],U)$.  Then we can consider the spectral sequence (\ref{spectr}) for
$\trp (U^*)$:
$E_1^{ij}(S)\implies \trp(U^*)(X[-i-j])= K(X,U[i+j])$.
 As in  \S \ref{hodge}, we note that its $E_1$-terms
depend only on $u_0(t_0(X))$ (they are just
$K_0(u_0(t_0(X)),v(U)[i+j])$). Since $u_0(t_0(X))= v(u(X))=0$, we have
$\trp (U^*)(X)=0$ for any $U$. Since $\trp (U^*)(X)=K(X,U)$, we
obtain that $X\cong 0$ in $K$. Hence $M\in \obj \mathfrak{B}$.
\end{proof}

Part 3 is a generalization of part 5 of Theorem \ref{ttn} (there
$T=\ns$).

\begin{coro}\label{coroloc}
Let $X$ belong to $\sv$. Then
 $\mg(X)$ is a mixed Tate motif (as described in part 3 of Remark
 \ref{rscat}) in $\dmge$ (i.e., we add direct summands) if and only if the
 complex $t_0(U)\in K^b(\chow)$ is.
\end{coro}

\begin{proof}
We apply part 3 of Proposition \ref{local} for $T=\{\sqcup_{a_i\ge
0}({\mathbb P}^1)^{a_i}\}.$

We obtain that $\mg^c(X)$ is a mixed Tate motif if and only if
$t_0(\mg^c(X))=W(X)$ is mixed Tate as an object of $K^b(\chow)$.

On the category of geometric of Voevodsky's motives $\dmgm\supset
\dmge$ (see \S\ref{stens} of this paper
 and \S4.3 of \cite{1}) we have a well-defined
duality such that $\z(n)^*=\z(-n)$. Therefore the category of mixed
Tate motives is a self-dual subcategory of $\dmgm$. Since
$\mg(X)^*=\mg^c(X)(-n)[-2n]$, $n=\dim X$ (see Theorem 4.3.7 of
\cite{1}; we can assume that $X$ is equidimensional), we obtain that
$\mg(X)$ is a Tate motif if and only if $\mg^c(X)$ is.
\end{proof}

One can translate this statement into a certain condition on the
motives $\mg(Y_{i_1}\cap Y_{i_2}\cap\dots \cap Y_{i_r})$ (in the
notation of Proposition \ref{mgc}).

\subsection{Motives in finite characteristic; $t\otimes\q$ is (conditionally)
an  isomorphism over a finite field}\label{motcharp}

In this subsection we consider Voevodsky'a motives over an arbitrary
perfect $k$.

In the proof of Theorem \ref{main} we used two facts:

(i) The cohomology of Suslin complex of a smooth projective variety
coincides with its hypercohomology.

(ii) Motives of smooth projective varieties generate $\dms$.

To the knowledge of the author, neither of these facts is known in
finite characteristic.

Yet for an covering of any variety over any perfect $k$ one has the
Mayer-Vietoris triangle (see Proposition 4.1.1 of \cite{1}); one
also has the blow-up distinguished triangle (see Proposition 3.5.3
of \cite{1} and Proposition 5.21 of \cite{3}). Besides, for a closed
embedding of smooth varieties one has the Gysin distinguished
triangle (see \cite{de2}).

\subsubsection{Certain integral arguments}\label{intmchp}

 In \cite{bev}  it was proved unconditionally
that $\dms$ has a differential graded enhancement. In fact, this
follows rather easily from the localization technique of Drinfeld
described in subsection \ref{loc} along with Voevodsky's description
of $\dms$ as a localization of $K^b(\smc)$. Moreover, Proposition
6.7 of \cite{bev} extends the Poincare duality for Voeovodsky
motives to our case. Therefore for $P,Q\in\spv$  we obtain
$$\dms(\mg(P), \mg(Q)[i])=\cho([P],[Q])\text { for }i=0;\ 0\text{
for }i>0 .$$ Hence the triangulated subcategory $DM_{pr}$ of $\dms$
generated by $[P],\ P\in\spv$ can be described as $\tr(I)$ for a
certain negative differential graded $I$ (see part 2 of Remark
\ref{prod}). Note also that the complexes $I(P,Q)$ compute motivic
cohomology of $P\times Q$.

If we define  the categories $I_N$ similarly to  $J_N$ 
(see
subsection \ref{tn}) then $I_0=\cho$. Hence there exists a
conservative weight complex functor $t_0:DM_{pr}\to K^b(\cho)$.
Moreover, for any 'enhanceable' realization of $DM_{pr}$ (again
one can easily check that these include motivic cohomology and
\'etale cohomology) and any $X\in \obj DM_{pr}$ one has the
motivic descent spectral sequence (\ref{spectr}).

The problem is that (to the knowledge of the author) at this moment
there is no way to prove that $DM_{pr}$ contains the motives of all
varieties (though it certainly contains the motives of varieties
that have 'smooth projective stratification').

Yet for any variety $X$ one could consider the functor $X^*_{pr}$
that equals $\dmge(-,\mg(X))$ restricted to $DM_{pr}$. It has a
differential graded 'enhancement'; in \S7.3 of \cite{bws} we prove
that $X^*_{pr}$ is representable by (at least) an object of a
certain infinite analogue of $DM_{pr}$.  This  gives a (possibly,
infinite) weight complex for $X$.

\subsubsection{Rational arguments}

For rational coefficients (ii) was proved in Appendix B of
\cite{hubka}: it suffices to note that if there exists an \'etale
finite morphism $U\to V$ for smooth $U,V$ then $\mg(V)\otimes \q$ is
a direct summand of $\mg(U)\otimes \q$. Hence the existence of de
Jong's alterations (see \cite{dej}) yields (ii) for motives with
rational coefficients (up to an idempotent completion of
categories).


Hence we obtain that in the characteristic $p$ case the category
$\dmge\q$ is the idempotent completion of $\tr(I)$ for a certain
negative differential graded $I$ (see above). Moreover (by the same
arguments as  in the characteristic $0$ case) there exists a
conservative  weight complex functor $t\q:\dmge\q\to K^b(\chow\q)$.

Now we show that, in contrast to the characteristic $0$ case, over
finite fields the Beilinson-Parshin conjecture (see below) implies
that $t \q$ is an equivalence.

\begin{pr}\label{becon}
Let $k$ be a finite field.

Suppose that for any $X\in \sv$ the only non-zero cohomology group
$H^i(X,\q(n))$ is $H^{2n}$.

Then $t\q$ is an equivalence of categories.
\end{pr}
\begin{proof}
By Theorem 4.2.2 of \cite{1} (note that it is valid with rational
coefficients in our case also!) for any $P,Q\in\spv$, $i\in\z$ we
have
\begin{equation}\label{fbpc}
\dme(\mg(P)\otimes \q, \mg(Q)\otimes\q[i])=H^{i+2r}(P\times
Q,\q(r)),\end{equation} where $r$ is the dimension of $Q$.
Therefore, if the Beilinson-Parshin conjecture holds then
$$\dme(\mg(P)\otimes \q,
\mg(Q)\otimes\q[i])=\cho\otimes\q([P],[Q])\text{
 for }i=0;\ 0\text{ otherwise}.$$ We obtain $$\dms\otimes
\q([P],[Q][i])\cong K^b(\cho \otimes
\q)(t_0\otimes\q([P]),t_0\otimes\q([Q][i]))$$ for all $i\in \z$.
Since $[P]$ for $P\in\spv$ generate $\dms\otimes\q$ as a
triangulated category, the same easy standard reasoning as the one
used in the proof of Theorem \ref{main} shows that $t_0\otimes \q$
induces a similar isomorphism for any two objects of
$\dms\otimes\q$; hence $t_0\otimes \q$ is a full embedding. Since
$[P],\ P\in\spv$ generate $K^b(\cho \otimes \q)$ as a triangulated
category, we obtain that $t_0\otimes \q$ is an equivalence of
categories. Lastly, since $\dmge\q$ is the idempotent completion of
$\dms\otimes \q$ and $K^b(\chow \q)$ is the idempotent completion of
$K^b(\cho\otimes\q)$,  $t\q$ is an equivalence also.
\end{proof}

\begin{rema}\label{nummotp}

1. Obviously,  the converse statement holds also: if $t \q$ is an
isomorphism then $\dmge\q (\mg(P)\otimes \q,\mg(Q)\otimes \q[i])=0$
for $i\neq 0$, $P,Q\in\sv$; hence the Beilinson-Parshin conjecture
holds.

2. As it was noted in Remark \ref{nummot}, using $t\q$ one can
obtain a functor $t_{num}\q:\dmge\q\to K^b(Mot^{eff}_{num}\q)$ (the
'wrong' homotopy category of effective numerical motives). Over a
finite field we can use Kunneth projectors for numerical motives
(cf. \S2 of \cite{mi}; note that the projectors are functorial!) to
replace $K^b(Mot^{eff}_{num}\q)$ by the 'correct' category
$K^b(\mathcal{M})$ such that $\mathcal{M}$ is tannakian, $H_i(P)$
(the $i$-th homological component) is put in degree $i$. Moreover,
the 'standard' conjectures assumed in Theorem 5.3 of \cite{mi} (Tate
conjecture $+$ numerical equivalence coincides with rational
equivalence) imply that $t_{num}$ can be extended to an equivalence
of $\dmgm\q$ with the motivic $t$-category of Theorem 5.3 of
\cite{mi} (see also Theorem 56 of \cite{skahn}). The relation
between different conjectures is discussed in 1.6.5 of \cite{skahn}.
\end{rema}

\subsection{An application: independence on $l$ for traces of open correspondences}\label{sindepl}

Our formalism easily implies that for any $X\in \obj\dmge$ and any
$f\in \dmge(X,X)$ the trace of the map $f^*$ induced on the
(rational)  $l$-adic \'etale   cohomology  $H(X)$  does not depend
on $l$. This gives a considerable generalization of Theorem 3.3 of
\cite{bloesn} (see below).

Now we formulate the main statement  more precisely. In this
subsection we do not demand the characteristic of $k$ to be $0$.

\begin{pr} \label{indepl}
Let $H$ denote any of rational $l$-adic \'etale   cohomology
theories or the singular cohomology corresponding to any embedding
of $k$ into $\mathbb{C}$ (in characteristic $0$). Denote by $\tr
f^*_H$ the sum $(-1)^i\tr f^*_{H^i(X)}$.

Then $\tr f^*_H$ 
does not depend on the choice of $H$.
\end{pr}
\begin{proof}
This statement is well known for Chow motives. Now we reduce
everything to this case.

We consider the  weight spectral sequence $S$ for $H(X)$ (see
\S\ref{hodge} and part 2 of Remark \ref{rzn}). Recall that   $S$ is
functorial in $H$ and $\dmge$-functorial in $X$.
 Hence any $f$ induces a (unique) endomorphism
$f^*_S$ of $S$. To construct $f^*_S$ one could also use the spectral
sequence $T$ described in Remark \ref{ass}. Hence to reduce the
statement to the case of Chow motives it suffices to apply the
following formula (that follows immediately from the general
spectral sequence formalism):
\begin{equation}
\tr f^*_H=\sum_{i,j}{-1}^{i+j}f^*_{H^i(P^j)};
\end{equation}
see III.7.4c of \cite{gelman}.
\end{proof}

Note that, as we have shown above, these rational arguments also
work in finite characteristic.

Now, Theorem 3.3 of \cite{bloesn} states the independence from the
choice of $H$ of traces of maps induced by {\it open
correspondences} on the cohomology of  $U\in\sv$ with compact
support. Open correspondences and the corresponding cohomology maps
were described in Definition 3.1 loc. cit. We will not recall this
definition here; we will only note that the map $\Gamma_*$ of
loc. cit. comes from a certain $f\in \dmge(\mg^c(U),\mg^c(U))$.
Indeed, $\Gamma_*$ was defined as a composition $(p_1)_*\circ
(p_2)^*$ (in the notation of \cite{bloesn}). One can define the
morphism of motives corresponding to $(p_2)^*$ using the
functoriality of $\mg^c$ with respect to proper morphisms (see \S4.1
of \cite{1}). To define the morphism corresponding to $(p_1)_*$ one
should use Corollary 4.2.4 loc. cit.

It follows that we have generalized Theorem 3.3 of \cite{bloesn}
considerably. Indeed, we do not demand $U$ to be a complement of a
smooth projective variety by a strict normal crossing divisor. This
is especially important in the characteristic $p$ case. Besides, it
seems that open correspondences of the type considered in
  Theorem 2.8 ibid.  also give
endomorphisms of $\mg^c(U)$.

\begin{rema}\label{rindepl}
1. Certainly, exactly the same argument proves the independence from
$H$ of $n_{\lambda}(H)=(-1)^i n_{\lambda} f^*_{H^i(X)}$; here
$n_{\lambda} f^*_{(H^i(X))}$ for a fixed algebraic number $\lambda$
denotes the algebraic multiplicity of the eigenvalue ${\lambda}$ for
the operator $f^*_{H^i(X)}$.

2. In fact, all these statements follow easily from the following
statement: the appropriately defined group $K_0(\en \chow)$ surjects
onto $K_0(\en \dmge)$; see \S5.4 of \cite{bws} for details.
\end{rema}

\subsection{Adding kernels of projectors to $J$}\label{summand}

The description of the derived category of Tate motives and of the
Tate twist on $\hk$ (see \S\ref{stens}) would be nicer if $\z(i),\
i\ge 0$ would be motives of length $0$ (see part 3 of Remark
\ref{rscat} and subsection \ref{motlen}). To this end we show that
one can easily add direct summands of objects to $J$. This means
that we consider a  differential graded category $J^+$ whose
objects are $\{(j,p):\ j\in \obj J,\ p\in J^0(j,j),\ j^2=j\}$. Its
morphisms are defined as usually for idempotent completions, i.e.,
$J^+ ((j,p), (j',p'))= j'\circ J(j,j')\circ p$.

Indeed, if the cohomology of a complex of sheaves coincides with its
hypercohomology, the same it true for any direct summand of this
complex. Therefore, if $D,D'$ are  direct summands (in $\cmk$) of
$C(P), C(P')$ respectively, $P,P'\in\spv$, then the natural analogue
of Proposition \ref{dmc} will be valid for $\dme(p(D),p(D'))$. It
easily follows that $m$ can be extended to a full embedding of
$\tr(J^+)$ into $\dme$ (cf.  Theorem \ref{main}).

In particular, let $P\subset Q\in\spv$ and let there exist a
section $j:Q\to P$ of the inclusion. Then  the cone of $j$
naturally lies in $J^+$ (note that the cone is isomorphic to
$m\ob(\mg^c(Q-P)$).

  For example, one can present $[\p^1]$ as $[\pt]\oplus
  \z(1)[2]$. Hence
  for $P\in\spv,\ i\ge 0$ one could consider $P[2i](i)$ as an object of  $J^+$
  (cf. the reasoning in the proof of part 1 of Theorem \ref{ttn}
and part 3 of Remark \ref{rscat}).

Yet this method (most probably) cannot give a (canonical)
differential enhancement of the whole $\dmge$ (or $\dmge{}'$). To
obtain an enhancement for it one could apply the 'infinite
diagram' method of Hanamura (see \S2 of \cite{h}).


\subsection{The functors $m_N$}\label{mn}

 We consider
the functor $J_N\to B^-(\ssc)$ that maps $[P]$ into $SC_N(P)$. Here
$SC_N^i(P)$ is the Nisnevich sheafification of the presheaf
$C_N^i(P)(-)$ (they coincide for $i\neq -N$). We consider the
corresponding functor $h_N:\hk_N\to K^-(\ssc)$ and $m_N=p\circ h_N$.
Note that for any $X\in \obj \hk_N$ we have $m_N(X)\in\dme$.

By part 1 of Proposition \ref{cqua} the natural morphisms $C_P\to
SC_N(P)$ in $\kmk$ induce a transformation of functors $\tr_N:m\to
m_N$. Besides $tr_N$ for any $X\in \hk$ is induced by a canonical
map in $\kmk$.

It seems that no nice analogue of part 3 of Theorem \ref{ttn} is
valid for $m_N$. Yet for low-dimensional varieties $m_N$ coincides
with $m$.

\begin{pr}

Suppose that the Beilinson-Soul\'e vanishing conjecture holds over
$k$. Then $m_{2n}(X)\cong m(X)$ in $\dme$ if the dimension of $X$ is
$\le n$.
\end{pr}
\begin{proof}
We check by induction on $n$ that $tr_{2n}$ is the identity for
$m(X)$. This is obviously valid for $n=0$.

Applying the same reasoning as  in the proof of part 1 of Theorem
\ref{ttn} we obtain that it is sufficient to prove the assertion for
smooth projective $X$ of dimension $\le n$.

It is sufficient to check that $\dme(\mg(Y)[N],\mg(X))=0$ for any
$Y\in \sv$, $N\ge 2n$.

By Theorem 4.3.2 of \cite{1} if the dimension of $X$ equals $n$ then
we have $\ihom_{\dme}(\mg(X), \z(n)[2n])\cong\mg(X)$. Hence
$$\begin{gathered} \dme(\mg(Y)[N],\mg(X))=\dme(\mg(Y\times X)[N],
\z(n)[2n])\\ = \dme(\mg(Y\times X)[N-2n], \z(n)). \end{gathered}$$
It remains to note that by the Beilinson-Soul\'e conjecture
$\dme(\mg(Y\times X), \z(n)[i])=0$ for $i<0$.
 \end{proof}

 \begin{rema}
1. We also see that for $k$ of characteristic $0$ and any $N$ there
exist $P,Q\in \spv$ such that $\hk(P[N],[Q])\neq 0$. Hence none of
$t_N$ and $m_N$ are full functors.

2. It can be also easily checked that $tr_{2n}$ being identical
for all $X$ of dimension $\le n$, $n\in\z$, implies the
Beilinson-Soul\'e conjecture.
 \end{rema}


\begin{thebibliography}{1}



\bibitem[BaS01]{ba} Balmer P.,  Schlichting M. Idempotent completion
of triangulated categories//  Journal of Algebra 236, no. 2 (2001),
819-834.


\bibitem[B-VRS03]{bv} Barbieri-Viale L., Rosenschon A., Saito M. Deligne's conjecture
on 1-motives// Ann. of Math. (2) 158 (2003), no. 2, 593--633.

\bibitem[BeV08]{bev} Beilinson A., Vologodsky V. A DG guide to Voevodsky
motives//  Geom. Funct. Analysis, vol. 17, no. 6, 2008,
1709--1787. 




 \bibitem[Blo94]{blo}  Bloch S.
 The moving lemma for higher Chow groups// Journ. of Alg. Geometry
3 (1994), no. 3, 537--568.


\bibitem[BlE07]{bloesn}  Bloch S., Esnault H., K\"unneth projectors for open
varieties, in: Algebraic cycles and motives, vol. 2, LMS lecture note series
344 (2007), 54--72.


\bibitem[BoK90]{bk} Bondal A. I., Kapranov M. M. Enhanced triangulated categories.
(Russian)// Mat. Sb. 181 (1990), no. 5, 669--683; translation in
Math. USSR-Sb. 70 (1991), no. 1, 93--107.


\bibitem[Bon10]{bws} Bondarko M.,
Weight structures vs.  $t$-structures; weight filtrations,
 spectral sequences, and complexes (for motives and in general)//
 J. of K-theory, v. 6, i.	03, pp. 387--504, 2010,
see also \url{http://arxiv.org/abs/0704.4003}


\bibitem[dJo06]{dej} de Jong A. J.
Smoothness, semi-stability and alterations// Inst. Hautes Etudes
Sci. Publ. Math. No. 83 (1996), 51--93.

\bibitem[Deg06]{de2} Deglise F.
 Around the Gysin triangle, preprint,
\url{http://www.math.uiuc.edu/K-theory/0764/}

\bibitem[Dri04]{drinf} Drinfeld V. DG quotients of DG categories// J. of Algebra 272
(2004), 643--691.

\bibitem[FrS02]{frisus} Friedlander E.,  Suslin A.
The spectral sequence relating algebraic $K$-theory to the motivic
cohomology// Ann. Sci. Ec. Norm. Sup. no. 35, 2002, 773--875.



\bibitem[GeM02]{gelman} Gelfand S., Manin Yu.,
  Methods of homological algebra.
  2nd ed. Springer Monographs in Mathematics.
  Springer-Verlag, Berlin, 2003. xx+372 pp.

\bibitem[GiS96]{gs} Gillet H., Soul\'e C. Descent, motives and $K$-theory//
J.  reine und angew. Math. 478, 1996, 127--176.

\bibitem[GiS07]{sg} Gillet H., Soul\'e C., Weight complexes  for arithmetic varieties, unpublished, see
\url{http://www.math.uic.edu/~henri/weightcomplexes\_slides.pdf}

\bibitem[GN-A02]{gu} Guillen F., Navarro Aznar V.
 Un critere d'extension des foncteurs definis sur les schemas lisses//
   Publ. Math. Inst. Hautes Etudes Sci. No. 95, 2002, 1--91.




\bibitem[Han95]{h} Hanamura M. Mixed motives and algebraic cycles, I//
Math. Res. Letters 2, 1995,  811--821.


\bibitem[Han00]{ha2} Hanamura M. Homological and cohomological motives of algebraic
varieties// Inv. Math 142, 2000,  319--349.

\bibitem[Han04]{hb} Hanamura M. Mixed motives and algebraic cycles, II//
Inv. Math 158, 2004,  105--179.

\bibitem[Han99]{ha3} Hanamura M. Mixed motives and algebraic cycles, III//
Math. Res. Letters 6, 1999,  61--82.

\bibitem[Hub00]{hu} Huber A.
Realizations of  Voevodsky's motives// J. Alg. Geometry,  9, 2000,
no. 4,  755--799.

\bibitem[HuK06]{hubka} Huber A., Kahn B.
The slice filtration and mixed Tate motives, Compos. Math. 142(4),
2006, 907 -- 936.


\bibitem[Jan92]{ja} Jannsen U. Motives, numerical equivalence and
semi-simplicity// Inv. Math 107, 1992,  447--452.

\bibitem[Kah05]{skahn} Kahn B., Algebraic K-theory, algebraic cycles and
arithmetic geometry, $K$-theory handbook, vol. 1, 351--428,
Springer-Verlag, Berlin, 2005.

\bibitem[KaS02]{kabir} Kahn B., Sujatha R., Birational motives, I,
preprint, 
\url{http://www.math.uiuc.edu/K-theory/0596/}


\bibitem[KrM95]{km} Kriz I., May J.  Operads, algebras, modules and motives,
Asterisque no. 233, 1995.


\bibitem[Lev09]{lesm} Levine M. Smooth motives,  
 in R. de Jeu, J.D. Lewis (eds.), Motives and Algebraic Cycles. A Celebration in Honour of Spencer J. Bloch, Fields Institute Communications 56, American Math. Soc. (2009), 175--231.


\bibitem[Lev05]{le1} Levine M. Mixed motives,  $K$-theory handbook,
vol. 1, 429--521, Springer-Verlag, Berlin, 2005.

\bibitem[Lev94]{le2} Levine M.  Bloch's higher Chow groups revisited// Asterisque
226, 10 (1994),  235--320.

\bibitem[Lev98]{le3} Levine M. Mixed motives, Math. surveys and
Monographs 57, AMS, Prov. 1998.

\bibitem[MVW06]{vbook} Mazza C., Voevodsky V.,  Weibel Ch. Lecture notes on motivic cohomology,
Clay Mathematics Monographs, vol. 2, 2006, see also
\url{http://www.math.rutgers.edu/~weibel/MVWnotes/prova-hyperlink.pdf}

\bibitem[MiR]{mi} Milne J., Ramachandran N. Motivic complexes over 
finite fields
 and the ring of correspondences at the generic point,
Pure \& App. Math. Quarterly (Tate issue), 5 (2009), 1219--1252.

\bibitem[Spi] {spit}  Spitzweck M.
 Some constructions for Voevodsky's triangulated categories of motives, preprint.

\bibitem[Sus03] {grf} Suslin A., On the Grayson spectral sequence// Tr. Mat. Inst. Steklova, 241, 2003, 218--253, translation in Proc. Steklov Inst. Math. 2003, no. 2
(241), 202--237.

 \bibitem[Voe00a]{1} Voevodsky V. Triangulated categories of motives over a field, in:
 Voevodsky V.,  Suslin A., and  Friedlander E.
 Cycles, transfers and motivic homology theories, Annals of
 Mathematical studies, vol. 143, Princeton University Press,
 2000,  188--238.

\bibitem[Voe00b]{3} Voevodsky V. Cohomological theory of presheaves with transfers,
same volume,  87--137.

\bibitem[FrV00]{4} Friedlander E.,  Voevodsky V. Bivariant cycle cohomology,
same volume,  138--187.

\bibitem[Voe96] {htop} Voevodsky, V. Homology of schemes//  Selecta Math. (N.S.) 2  (1996),  no. 1, 111--153.

\bibitem[Voe10]{voevc} Voevodsky, V. Cancellation theorem// 
Doc. Math., extra volume: Andrei Suslin's Sixtieth Birthday (2010), 671--685.


\end{thebibliography}
\end{document}